\documentclass[reqno, a4paper]{amsart}
\usepackage[left=3cm, right=3cm, top=3cm, bottom=3cm]{geometry}
\usepackage{amssymb, amsmath, amsthm}
\usepackage{amsfonts, bm, bbm, mathrsfs}
\usepackage{cases}
\usepackage{empheq}
\usepackage{enumitem}
\setlist[enumerate]{label=\upshape(\arabic*)}
\usepackage[colorlinks=true]{hyperref}
\newtheorem{theorem}{Theorem}[section]
\newtheorem{assumption}[theorem]{Assumption}
\newtheorem{lemma}[theorem]{Lemma}
\newtheorem{corollary}[theorem]{Corollary}
\newtheorem{proposition}[theorem]{Proposition}
\newtheorem{condition}[theorem]{Condition}
\theoremstyle{definition}
\newtheorem{definition}[theorem]{Definition}
\newtheorem{example}[theorem]{Example}
\theoremstyle{remark}
\newtheorem{remark}{Remark}
\allowdisplaybreaks
\begin{document}
    \title{Well-posedness and large deviations for 3D primitive equations with jumps}
    \author[Z. Shen]{Zhuoyang Shen}
    \address[Z. Shen]{School of Mathematics and Statistics, Lanzhou University, Lanzhou, 730000, P.R. China}
    \email{shenzhy2023@lzu.edu.cn}
    \author[X. Li]{Xinhua Li}
    \address[X. Li]{School of Mathematics and Statistics, Lanzhou University, Lanzhou, 730000, P.R. China}
    \email{lixinhua@lzu.edu.cn}
    \author[C. Sun]{Chunyou Sun$^*$}
    \address[C. Sun]{School of Mathematics and Statistics, Lanzhou University, Lanzhou, 730000, P.R. China;
    School of Mathematics and Statistics, Donghua University, Shanghai, 201620, P.R. China}
    \email{sunchy@lzu.edu.cn, sunchunyou@gmail.com}
    \thanks{$^*$Corresponding author}
 \begin{abstract}
We establish the global pathwise well-posedness and a Freidlin--Wentzell large deviation principle for the three-dimensional viscous primitive equations driven by multiplicative L\'evy noise represented by a compensated Poisson random measure. Under local Lipschitz and moment growth assumptions on the noise coefficient, we construct local martingale solutions and obtain a unique maximal pathwise solution with c\`adl\`ag trajectories in the $H^1$ energy space. Additional sixth-moment growth in $H^1$ and $L^6$ growth bounds on the baroclinic component of the noise yield global pathwise well-posedness. The proof combines the jump It\^o formula with baroclinic and anisotropic energy estimates to rule out finite-time blow-up. Under further exponential integrability and structural assumptions on the noise coefficient, we establish a Freidlin--Wentzell large deviation principle in a strong path topology by the weak convergence approach.
\end{abstract}
    \maketitle
    \tableofcontents
    \section{Introduction}\label{sec:intro}
\subsection{The primitive equations}

The primitive equations describe large-scale oceanic and atmospheric flows under the Boussinesq and hydrostatic approximations. In this paper, we study global pathwise well-posedness and large deviations for the three-dimensional equations driven by multiplicative L\'evy noise. We work on the cylindrical domain $\mathcal O=\mathcal D\times(-h,0)$, where $\mathcal D\subset\mathbb R^2$ is bounded and smooth. The unknowns are the horizontal velocity $\bm v=(v_1,v_2)$, the vertical velocity $w$, the pressure $p$ and the temperature $T$. The system reads
\begin{equation}\label{eq:intro-model}
    \begin{cases}
        \mathrm{d}\bm v+\big((\bm v\cdot\nabla)\bm v+w\partial_z\bm v +f\bm k\times\bm v+\nabla p+L_1\bm v\big)\mathrm{d}t =\displaystyle\int_E G_1(\bm u(t-),\xi)\,\tilde\eta(\mathrm{d}t,\mathrm{d}\xi), \\
        \partial_zp+T=0,\quad \mathrm{div}\,\bm v+\partial_zw=0, \\
        \mathrm{d}T+\big((\bm v\cdot\nabla)T+w\partial_zT+L_2T\big)\mathrm{d}t =\displaystyle\int_E G_2(\bm u(t-),\xi)\,\tilde\eta(\mathrm{d}t,\mathrm{d}\xi),
    \end{cases}
\end{equation}
where $\bm u=(\bm v,T)^\top$, $\nabla$ and $\mathrm{div}$ act in the horizontal variables, $f$ is the Coriolis parameter and $\bm k$ is the vertical unit vector. The operators $L_1,L_2$ describe full viscosity and thermal diffusion. Here $\eta$ is a Poisson random measure on $[0,\infty)\times E$ with intensity measure $\mathrm{d}t\otimes\nu$, where $E$ is a locally compact Polish space, and $\tilde\eta=\eta-\mathrm{d}t\otimes\nu$ is its compensated version. The boundary conditions are specified in Subsection~\ref{subsec:model}. In particular, the temperature satisfies the Robin condition $\partial_zT=-\alpha_uT$, $\alpha_u>0$, on the top surface and is insulated elsewhere.

For the deterministic equations with full dissipation, Lions, Temam and Wang \cite{0766.35039} established global weak solutions, and Cao and Titi \cite{1151.35074} proved global well-posedness of strong solutions. A central step in the latter proof is an $L^6$ estimate for the baroclinic mode, which leads to the $H^1$ bound. Further global regularity results were obtained by Kobelkov \cite{Kobelkov2007} and Kukavica and Ziane \cite{KukavicaZiane2007}. Hieber and Kashiwabara \cite{HieberKashiwabara2016} developed an $L^p$ approach to global strong well-posedness.

\subsection{Stochastic primitive equations}

For equations driven by Wiener noise, well-posedness has been studied in both two and three dimensions. Glatt-Holtz and Temam \cite{GlattHoltzTemam2011} constructed pathwise solutions in two dimensions. In three dimensions, Guo and Huang \cite{1179.37115} studied global well-posedness and attractors under additive white noise, and Debussche, Glatt-Holtz, Temam and Ziane \cite{1246.35194} proved global existence and regularity under multiplicative white noise. Gao and Sun \cite{1354.35160} obtained existence and uniqueness under smallness conditions on the initial data and the noise. Brzeźniak and Slavík \cite{BrzezniakSlavik2021} allowed gradient-dependent multiplicative and transport noise, proving the existence of unique maximal solutions that are strong in both the probabilistic and PDE senses. They also established global existence under additional conditions. Glatt-Holtz, Kukavica, Vicol and Ziane \cite{1338.37065} studied the long-time behavior by constructing invariant measures.

Large deviations have also been studied for the three-dimensional equations. Dong, Zhai and Zhang \cite{1370.60103} proved a Freidlin--Wentzell large deviation principle under multiplicative white noise. Slavík \cite{Slavik2022} established large and moderate deviation principles and an almost sure central limit theorem for gradient-dependent multiplicative white noise. Agresti and Theewis \cite{AgrestiTheewis2025} proved a large deviation principle for transport noise with turbulent pressure.

L\'evy noise provides a way to model intermittent perturbations through jumps. For the Navier--Stokes equations, Sakthivel and Sritharan \cite{SakthivelSritharan2012} constructed martingale solutions, and Motyl \cite{Motyl2013,Motyl2014} developed a martingale solution framework in unbounded three-dimensional domains. For the two-dimensional Navier--Stokes equations with jumps, Brzeźniak, Peng and Zhai \cite{1519.35220} proved global strong well-posedness and a large deviation principle using the Poisson control framework of \cite{1259.60065,1231.60018}.

For the primitive equations with L\'evy noise, Sun and Gao \cite{SunGao2013} established well-posedness in two dimensions, while Sun, Gao and Li \cite{SunGaoLi2018} and Ma and Zhang \cite{MaZhang2021} proved large deviation principles for additive and multiplicative L\'evy noise, respectively. In three dimensions, Su and Liu \cite{SuLiu2023} studied exponential stability, assuming the existence of global weak solutions. To develop a strong solution theory for \eqref{eq:intro-model}, one must also control the jumps in the nonlinear energy identities used for the three-dimensional estimates.

We establish global pathwise well-posedness for \eqref{eq:intro-model} at the $H^1$ level and prove a Freidlin--Wentzell large deviation principle for the resulting strong solutions. We first construct local martingale solutions that are strong in the PDE sense and have c\`adl\`ag $H^1$ paths. Pathwise uniqueness gives a unique maximal pathwise solution. Additional sixth-moment and baroclinic $L^6$ growth bounds imply that this solution is global for every square-integrable $H^1$-valued initial datum, without a smallness restriction. Under further exponential integrability and structural assumptions on the noise, the large deviation principle holds at the same spatial regularity as the well-posedness result. Both results use Lipschitz control in the $L^2$ norm, without an $H^1$-Lipschitz assumption on the noise coefficient.

\subsection{Main results}

Let $H$ and $V$ denote the $L^2$ and $H^1$ energy spaces, and let $A$ be the dissipation operator with domain $\mathscr D(A)$. The functional setting and solution notions are specified in Subsections \ref{subsec:fs} and \ref{subsec:solutions}.

The local theory uses Assumptions \ref{asm:martingale} and \ref{asm:A-est}. These impose Lipschitz continuity in the $H$ norm on bounded subsets of $V$, componentwise moment growth bounds of orders $2,4$ in $H$, and moment growth bounds of orders $2,4$ in $V$.

\begin{theorem}[Local martingale solutions]\label{thm:intro-martingale}
    Let Assumptions \ref{asm:martingale} and \ref{asm:A-est} hold, and let $\bm u_0\in L^2(\Omega,\mathcal F_0;V)$. Then, for every $\mathcal T>0$, \eqref{eq:intro-model} admits a local martingale solution $(\hat\eta,\hat{\bm u},\tau)$ whose initial value has the same law as $\bm u_0$. Its lifetime satisfies $0<\tau\leq\mathcal T$ almost surely, and
    \begin{equation*}
        \hat{\mathbb E}\Big[\sup_{t\in[0,\tau]} \|\hat{\bm u}(t)\|_V^2 +\int_0^\tau\|\hat{\bm u}(t)\|_{\mathscr D(A)}^2 \,\mathrm{d}t\Big]<\infty.
    \end{equation*}
    The solution has $V$-valued c\`adl\`ag paths on $[0,\tau]$.
\end{theorem}

\begin{theorem}[Pathwise uniqueness and maximal solution]
    \label{thm:intro-pathwise} Under the assumptions of Theorem \ref{thm:intro-martingale}, pathwise uniqueness holds among local solutions satisfying its energy bound: solutions with the same initial datum and Poisson random measure on the same stochastic basis agree on their common existence interval. Moreover, \eqref{eq:intro-model} has a unique maximal pathwise solution $(\bm u,\xi)$ on the prescribed stochastic basis.
\end{theorem}

For global well-posedness, Assumption \ref{asm:global} adds a sixth-moment growth bound in $V$ and second- and sixth-moment growth bounds for the baroclinic component of the velocity noise in $L^6$.

\begin{theorem}[Global well-posedness]\label{thm:intro-global}
    Let the assumptions of Theorem \ref{thm:intro-martingale} and Assumption \ref{asm:global} hold. Then the maximal solution of Theorem \ref{thm:intro-pathwise} satisfies $\xi=\infty$ almost surely. For every $\mathcal T>0$,
    \begin{equation*}
        \sup_{t\in[0,\mathcal T]}\|\bm u(t)\|_V^2 +\int_0^{\mathcal T}\|\bm u(t)\|_{\mathscr D(A)}^2\,\mathrm{d}t <\infty\qquad\mathbb P\text{-almost surely}.
    \end{equation*}
\end{theorem}

Pathwise well-posedness allows us to express the solution as a measurable function of the initial datum and the Poisson random measure. The large deviation argument also requires estimates for equations with controlled Poisson intensity, uniformly over controls with bounded entropy, and convergence of their solutions.

Fix $\mathcal T>0$ and a deterministic initial datum $\bm u_0\in V$. Let $\bm u^\varepsilon$ solve \eqref{eq:intro-model} with jump amplitude $\varepsilon$ and Poisson intensity $\varepsilon^{-1}\nu(\mathrm{d}\xi)\,\mathrm{d}t$, as specified in \eqref{eqn:spe}. Set
\begin{equation*}
    \Upsilon_{\mathcal T}(V) =D([0,\mathcal T];V)\cap L^2(0,\mathcal T;\mathscr D(A)),
\end{equation*}
equipped with the intersection topology of the Skorokhod $J_1$ topology on the first factor and the norm topology on the second. Assumption \ref{assump:LDP} requires exponential integrability of the Lipschitz and growth envelopes. It also controls how the noise coefficient enters the combined barotropic and vertical energy estimate. The velocity noise coefficient may depend on temperature, as illustrated in Example \ref{ex:coupled-noise}.

\begin{theorem}[Large deviation principle]\label{thm:intro-ldp}
    Under Assumption \ref{assump:LDP}, the family $\{\bm u^\varepsilon\}_{\varepsilon\in(0,1)}$ satisfies a Freidlin--Wentzell large deviation principle on $\Upsilon_{\mathcal T}(V)$ with speed $\varepsilon^{-1}$ and good rate function
    \begin{equation*}
        I(\phi)=\inf\big\{L_{\mathcal T}(g): g\in\mathbb S,\ \bm u^g=\phi\big\}, \qquad \phi\in\Upsilon_{\mathcal T}(V).
    \end{equation*}
    Here $\bm u^g$ solves the skeleton equation \eqref{eqn:skeleton}, $\mathbb S$ is the set of nonnegative measurable controls with finite entropy, and
    \begin{equation*}
        L_{\mathcal T}(g)=\int_0^{\mathcal T}\int_E \big(g\log g-g+1\big)\,\nu(\mathrm{d}\xi)\,\mathrm{d}t,
    \end{equation*}
    with $0\log0=0$ and $\inf\varnothing=\infty$.
\end{theorem}

\subsection{Proof strategy}

We first construct local martingale solutions in the $H^1$ energy space. The basic $L^2$ estimate does not control the hydrostatic convection at this level. Following Brzeźniak and Slavík \cite{BrzezniakSlavik2021}, we therefore truncate the convection outside a small $V$-neighbourhood of the linear evolution of the initial datum. Uniform energy bounds give compactness. The local $H$-Lipschitz condition, combined with localization in $V$ and fourth-moment estimates, then justifies passage to the limit in the noise coefficient. We identify the limiting martingale solution together with its Poisson driver, so that pathwise uniqueness can be applied in the Gyöngy--Krylov argument \cite{GyongyKrylov1996}. This yields the regularised solution on the original stochastic basis; its stopped version is a local pathwise solution of the original equation.

At a jump exit time, the solution may exceed the energy threshold even though its pre-exit energy is bounded. The continuation argument therefore requires an estimate of the terminal state, including the exit jump. A linear equation with the localized drift and noise coefficients gives a $V$-valued c\`adl\`ag extension whose terminal value lies in $L^2(\Omega;V)$. We apply the local existence theorem with this value as initial datum to extend the solution beyond the stopping time and prove the energy blow-up criterion.

To obtain global solutions, we adapt the baroclinic $L^6$ and derivative estimates of Cao and Titi \cite[Sections~3.2--3.3]{1151.35074}. We control the It\^o remainder by the second- and sixth-moment growth bounds on the baroclinic noise; the sixth-order contribution is absent in the Wiener case. This yields the baroclinic bound required for the derivative estimates. The random coefficients and source terms in these successive estimates have time integrals controlled by the preceding bounds. Localizing these integrals gives uniform bounds in probability for the full $H^1$ energy on the local existence intervals. The blow-up criterion therefore excludes finite-time explosion. Localization of the initial datum extends the result to square-integrable $V$-valued initial data.

The large deviation principle requires estimates that are uniform over controls with bounded entropy. For the skeleton equation, the joint growth condition controls the forcing in a combined estimate for the barotropic gradient and vertical derivative, from which the global $H^1$ bounds follow. We must also deduce convergence in $\Upsilon_{\mathcal T}(V)$ from convergence of the control measures, although the control forcing is only in $L^1(0,\mathcal T;V)$ and the noise coefficient is Lipschitz only in $H$. Following \cite{1519.35220,1370.60103}, we introduce auxiliary linear equations driven by the change in control forcing and by the jump martingale. Energy and entropy estimates show that their solutions converge to zero. Convergence in $L^2(0,\mathcal T;V)$, combined with uniform integrability of the control coefficients, makes the remaining weighted errors vanish in the $V$ estimate. We thus obtain continuity of the skeleton solution map and convergence of the controlled stochastic solutions, as required by the weak convergence criterion of \cite{1259.60065}, based on the Poisson variational representation in \cite{1231.60018}.

\subsection{Organisation of the paper}

Section \ref{sec:framework} introduces the functional setting and the notions of solution. Section \ref{sec:wp-martingale} constructs local martingale solutions through the regularised Galerkin scheme. Section \ref{sec:pathwise} proves pathwise uniqueness and the existence of a maximal pathwise solution, followed by global well-posedness. Section \ref{sec:ldp} establishes the large deviation principle.

    \section{Mathematical framework}\label{sec:framework}
\subsection{The stochastic primitive equations (SPEs)}\label{subsec:model}
Let $\mathcal D\subset\mathbb R^2$ be a bounded domain with smooth boundary $\partial \mathcal D$, and set $\mathcal O=\mathcal D\times(-h,0)$. We consider the three-dimensional stochastic primitive equations (SPEs) on $\mathcal O\times[0,\infty)$ driven by multiplicative L\'evy noise:
\begin{empheq}[left=\empheqlbrace]{align}
    & \displaystyle\,\mathrm{d}\bm v + ((\bm v\cdot\nabla)\bm v + w\partial_z\bm v+ f\bm k\times\bm v + \nabla p + L_1\bm v)\,\mathrm{d}t=\int_EG_1(\bm v(t-),T(t-),\xi)\tilde{\eta}(\mathrm{d}t,\mathrm{d}\xi),\label{eqn:o-1} \\
    & \partial_zp + T = 0,\label{eqn:o-2} \\
    & \mathrm{div}\bm v + \partial_z w = 0,\label{eqn:o-3} \\
    & \displaystyle\,\mathrm{d} T + ((\bm v\cdot \nabla) T + w\partial_z T+ L_2T)\,\mathrm{d}t=\int_EG_2(\bm v(t-),T(t-),\xi)\tilde{\eta}(\mathrm{d}t,\mathrm{d}\xi),\label{eqn:o-4}
\end{empheq}
where $\bm v=(v_1,v_2)$, $w$, $p$ and $T$ denote the horizontal velocity, vertical velocity, pressure and temperature, respectively. Here $f$ is the Coriolis parameter, $\bm k$ is the vertical unit vector, and $\tilde\eta$ is the compensated Poisson random measure specified below. We write $\nabla=(\partial_x,\partial_y)$, $\nabla_3=(\nabla,\partial_z)$, $\Delta=\partial_{xx}+\partial_{yy}$ and $\mathrm{div}\,\bm v=\partial_xv_1+\partial_yv_2$. The viscosity and thermal diffusion operators are
\begin{gather*}
    L_1\bm v:=- \frac{1}{Re_1}\Delta\bm v - \frac{1}{Re_2}\partial_{zz}\bm v, \\
    L_2T:=- \frac{1}{Rt_1}\Delta T - \frac{1}{Rt_2}\partial_{zz}T,
\end{gather*}
where $Re_1, Re_2$ are Reynolds numbers and $Rt_1, Rt_2$ are the horizontal and vertical thermal parameters, respectively. For simplicity, we take
\begin{equation*}
    Re_1 = Re_2 = Rt_1 = Rt_2 = 1.
\end{equation*}

As in \cite{1246.35194}, the boundary $\partial\mathcal{O}$ is partitioned into the top $\Gamma_u:=\mathcal{D}\times\{0\}$, the bottom $\Gamma_b:=\mathcal{D}\times\{-h\}$ and the sides $\Gamma_l:=\partial \mathcal{D}\times(-h,0)$. We denote by $\bm n$ the outward unit normal to $\partial \mathcal{D}$, and by $\alpha_u>0$ the heat-exchange coefficient at the surface. We impose the boundary conditions

\begin{align}
    \partial_z\bm v=0,\ w=0,\ \partial_zT=-\alpha_uT\ & \text{on }\Gamma_u,\label{bd_cond:o-1} \\
    \partial_z\bm v=0,\ w=0,\ \partial_zT=0\ & \text{on }\Gamma_b,\label{bd_cond:o-2} \\
    \displaystyle \bm v\cdot\bm n=0,\ \frac{\partial\bm v}{\partial\bm n}\times\bm n=0,\ \frac{\partial T}{\partial\bm n}=0\ & \text{on }\Gamma_l\label{bd_cond:o-3}.
\end{align}

Integrating \eqref{eqn:o-2} and \eqref{eqn:o-3} from $-h$ to $z$, respectively, and using \eqref{bd_cond:o-2}, we have
\begin{align}
    p(x,y,z,t) & = -\int_{-h}^zT(x,y,z^\prime,t)\,\mathrm{d}z^\prime + p_b(x,y,t), \\
    w(x,y,z,t) & = -\int_{-h}^z\mathrm{div}\bm v(x,y,z^\prime,t)\,\mathrm{d}z^\prime =: \Phi(\bm v)(x, y, z, t),
\end{align}
where $p_b$ is the pressure on $\Gamma_b$. Moreover, by \eqref{bd_cond:o-1} we have
\begin{equation}\label{id:div_v}
    \int_{-h}^0\mathrm{div}\bm v\,\mathrm{d}z = - w(x,y,0,t) = 0.
\end{equation}
Then \eqref{eqn:o-1}--\eqref{eqn:o-4} can be rewritten as
\begin{empheq}[left=\empheqlbrace]{align}
    & \begin{aligned}\label{eqn:r-1}
        \displaystyle\,\mathrm{d}\bm v + ((\bm v\cdot\nabla)\bm v + \Phi(\bm v)\partial_z\bm v + f\bm k\times\bm v + \nabla p_b - \int_{-h}^z\nabla T\,\mathrm{d}z^\prime & + L_1\bm v)\,\mathrm{d}t \\
        = & \int_EG_1(\bm v(t-),T(t-),\xi)\tilde{\eta}(\mathrm{d}t,\mathrm{d}\xi),
    \end{aligned} \\
    & \displaystyle\,\mathrm{d}T + ((\bm v\cdot \nabla) T + \Phi(\bm v)\partial_z T + L_2T)\,\mathrm{d}t=\int_EG_2(\bm v(t-),T(t-),\xi)\tilde{\eta}(\mathrm{d}t,\mathrm{d}\xi),\label{eqn:r-2} \\
    & \displaystyle \int_{-h}^0\mathrm{div}\bm v\,\mathrm{d}z = 0.\label{eqn:r-3}
\end{empheq}
The boundary conditions become
\begin{align}
    \partial_z\bm v=0,\ \partial_zT=-\alpha_uT\ & \text{on }\Gamma_u,\label{bd_cond:r-1} \\
    \partial_z\bm v=0,\ \partial_zT=0\ & \text{on }\Gamma_b,\label{bd_cond:r-2} \\
    \displaystyle \bm v\cdot\bm n=0,\ \frac{\partial\bm v}{\partial\bm n}\times\bm n=0,\ \frac{\partial T}{\partial\bm n}=0\ & \text{on }\Gamma_l.\label{bd_cond:r-3}
\end{align}
Writing $\bm u=(\bm v,T)^\top$, we prescribe the initial condition
\begin{equation}
    \bm u(0) = \bm u_0 = (\bm v_0,T_0)^\top.
\end{equation}

\subsection{Functional spaces and abstract functionals}\label{subsec:fs}
Let $L^p(\mathcal{O})$, $L^p(\mathcal{D})$ and $H^m(\mathcal{O}),H^m(\mathcal{D})$ be the classical Lebesgue spaces and Sobolev spaces, respectively. For $1\leq p<\infty$, set
\begin{equation*}
    |\phi|_p= \begin{cases}
        \displaystyle\left(\int_{\mathcal{O}}|\phi(x,y,z)|^p\,\mathrm{d}x\,\mathrm{d}y\,\mathrm{d}z\right)^{1/p}, & \forall\phi\in L^p(\mathcal{O}), \\
        \displaystyle\left(\int_{\mathcal{D}}|\phi(x,y)|^p\,\mathrm{d}x\,\mathrm{d}y\right)^{1/p}, & \forall\phi\in L^p(\mathcal{D}).
    \end{cases}
\end{equation*}
For $k\in\mathbb{N}^+$, set
\begin{equation*}
    \|\psi\|_{H^k}^2=\sum_{0\leq|\alpha|\leq k}|\partial^\alpha\psi|_2^2,\quad \forall \psi\in H^k(\mathcal{O})\ \text{or}\ H^k(\mathcal{D}).
\end{equation*}

Define the working spaces for \eqref{eqn:r-1}--\eqref{eqn:r-3}
\begin{align*}
    \mathcal{V}_1 & :=\left\{\bm v\in (C^\infty(\overline{\mathcal{O}}))^2:\frac{\partial\bm v}{\partial z}\big|_{\Gamma_u,\Gamma_b}=0,\ \bm v\cdot\bm n=0,\ \frac{\partial\bm v}{\partial\bm n}\times\bm n=0\ \text{on }\Gamma_l,\ \int_{-h}^0\mathrm{div}\bm v\,\mathrm{d}z=0\right\}, \\
    \mathcal{V}_2 & :=\left\{T\in C^\infty(\overline{\mathcal{O}}):\partial_zT=-\alpha_uT\ \text{on }\Gamma_u,\ \partial_zT=0\ \text{on }\Gamma_b,\ \frac{\partial T}{\partial\bm n}\big|_{\Gamma_l}=0\right\}.
\end{align*}
Let $H_1$ and $H_2$ be the closures of $\mathcal{V}_1$ and $\mathcal{V}_2$ in $L^2(\mathcal{O})\times L^2(\mathcal{O})$ and $L^2(\mathcal{O})$, respectively. Similarly, let $V_1$ and $V_2$ be the closures of $\mathcal{V}_1$ and $\mathcal{V}_2$ in $H^1(\mathcal{O})\times H^1(\mathcal{O})$ and $H^1(\mathcal{O})$, respectively. We denote the norms and inner products on $V_i$ by $\|\cdot\|_{V_i}$ and $\langle\cdot,\cdot\rangle_{V_i}$, and on $H_i$ by $|\cdot|_{H_i}$ and $\langle\cdot,\cdot\rangle_{H_i}$, for $i=1,2$, where the $H_i$ inner products are induced by $L^2$. We equip $V_i$ with the energy inner products $\langle\cdot,\cdot\rangle_{V_i}:=a_i(\cdot,\cdot)$ defined below. These induce the same topology as $H^1$, by Lemma \ref{lem:Stokes_op}. In particular, $V_2=H^1(\mathcal O)$ and $H_2=L^2(\mathcal O)$.

Set $V=V_1\times V_2$ and $H=H_1\times H_2$, with inner products and norms
\begin{gather*}
    \langle \bm u,\hat{\bm u}\rangle_V=\langle\bm v,\hat{\bm v}\rangle_{V_1}+\langle T,\hat{T}\rangle_{V_2},\quad \|\bm u\|_V=\langle \bm u,\bm u\rangle_V^{1/2}, \\
    \langle \bm u,\hat{\bm u}\rangle_H=\langle\bm v,\hat{\bm v}\rangle_{H_1}+\langle T,\hat{T}\rangle_{H_2},\quad |\bm u|_H=\langle \bm u,\bm u\rangle_H^{1/2},
\end{gather*}
where $\bm u=(\bm v,T)^\top$ and $\hat{\bm u}=(\hat{\bm v},\hat T)^\top$.

Denoting the duals of $V$ and $H$ by $V'$ and $H'$, we obtain the Gelfand triple
\begin{equation*}
    V\subset H\cong H^\prime\subset V^\prime.
\end{equation*}
Define the functionals $a:V\times V\to\mathbb{R}$, $a_1:V_1\times V_1\to\mathbb{R}$, $a_2:V_2\times V_2\to\mathbb{R}$ and their corresponding linear operators $A:V\to V^\prime$, $A_1:V_1\to V_1^\prime$, $A_2:V_2\to V_2^\prime$ by
\begin{equation*}
    a(\bm u,\hat{\bm u}):={}_{V^\prime}\langle A\bm u,\hat{\bm u}\rangle_V=a_1(\bm v,\hat{\bm v})+a_2(T,\hat{T}),
\end{equation*}
where
\begin{align*}
    a_1(\bm v,\hat{\bm v}) & :={}_{V^\prime_1}\langle A_1\bm v,\hat{\bm v}\rangle_{V_1}=\int_{\mathcal{O}}\left(\nabla\bm v\cdot\nabla\hat{\bm v}+\frac{\partial\bm v}{\partial z}\cdot\frac{\partial\hat{\bm{v}}}{\partial z}\right)\,\mathrm{d}x\,\mathrm{d}y\,\mathrm{d}z, \\
    a_2(T,\hat{T}) & :={}_{V^\prime_2}\langle A_2T,\hat{T}\rangle_{V_2}=\int_{\mathcal{O}}\left(\nabla T\cdot\nabla\hat{T}+\frac{\partial T}{\partial z}\cdot\frac{\partial\hat{T}}{\partial z}\right)\,\mathrm{d}x\,\mathrm{d}y\,\mathrm{d}z+\alpha_u\int_{\Gamma_u}T\,\hat{T}\,\mathrm{d}x\,\mathrm{d}y.
\end{align*}
We also denote by $A_i$ and $A$ their realizations on $H_i$ and $H$, respectively, with
\begin{equation}\label{def:operator-domain}
    \mathscr D(A_i)=\{u\in V_i:A_i u\in H_i\},\qquad \mathscr D(A)=\{\bm u\in V:A\bm u\in H\}.
\end{equation}
Let $P_{H_1}$ be the orthogonal projection from $L^2(\mathcal O;\mathbb R^2)$ onto $H_1$, and set $P_H=\operatorname{diag}(P_{H_1},I)$.

\begin{lemma}\label{lem:Stokes_op}
    The form $a$ is continuous and coercive on $V$, and
    \begin{equation}\label{est:a}
        |a(\bm u,\hat{\bm u})|\leq\|\bm u\|_V\|\hat{\bm u}\|_V, \qquad a(\bm u,\bm u)=\|\bm u\|_V^2\asymp\|\bm u\|_{H^1}^2.
    \end{equation}
    The realization $A$ on $H$ is positive and self-adjoint with compact inverse. Its domain is
    \begin{equation}\label{eq:operator-domain-bc}
        \mathscr D(A)=\{\bm u\in V\cap H^2(\mathcal O;\mathbb R^3): \bm u\text{ satisfies }\eqref{bd_cond:r-1}\text{--}\eqref{bd_cond:r-3}\}.
    \end{equation}
    For $\bm u=(\bm v,T)\in\mathscr D(A)$,
    \begin{equation*}
        A\bm u=(P_{H_1}L_1\bm v,L_2T),\qquad |A\bm u|_H\asymp\|\bm u\|_{H^2}.
    \end{equation*}
    Moreover,
    \begin{equation*}
        \mathscr D(A^{1/2})=V,\qquad |A^{1/2}\bm u|_H=\|\bm u\|_V.
    \end{equation*}
\end{lemma}

\begin{proof}
    For temperature, coercivity follows from
    \begin{equation*}
        |T|_2^2\leq2h|T(\cdot,0)|_{L^2(\mathcal D)}^2+2h^2|\partial_zT|_2^2.
    \end{equation*}
    For velocity, the Poincar\'e inequality on each horizontal slice gives $|\bm v|_2\leq C|\nabla\bm v|_2$, since the normal trace vanishes and the only constant field tangent to $\partial\mathcal D$ is zero. The representation theorem for closed symmetric forms and the compact embedding $V\hookrightarrow H$ give positivity, self-adjointness, compactness of the inverse and the square-root identities.

    To identify the domain, set $\overline{\bm v}=h^{-1}\int_{-h}^0\bm v\,\mathrm{d}z$ and decompose $\bm v=\overline{\bm v}+\widetilde{\bm v}$. The horizontal operator is the two-dimensional Stokes operator on the barotropic component and the vector Laplacian on the zero-average component, both with $\bm v\cdot\bm n=0$ and $(\partial_{\bm n}\bm v)\times\bm n=0$. Elliptic regularity on the smooth domain $\mathcal D$ gives their $H^2$ domains. Combining these with the vertical Neumann operator gives the velocity domain in \eqref{eq:operator-domain-bc}. For temperature, combine the horizontal Neumann operator with the vertical Robin--Neumann operator. Products of the corresponding eigenfunctions, with eigenvalues $\lambda,\mu$, give the domain of their sum; the bound $\lambda\mu\leq(\lambda+\mu)^2/4$ controls the mixed derivatives. Integration by parts identifies $A$ with the stated differential operator. The resulting elliptic estimate and coercivity give $\|\bm u\|_{H^2}\asymp|A\bm u|_H$.
\end{proof}

We equip $\mathscr D(A)$ with the norm $\|\bm u\|_{\mathscr D(A)}:=|A\bm u|_H$.

Set $W=V\cap H^2(\mathcal O;\mathbb R^3)$, with its $H^2$ norm. For smooth fields define
\begin{equation*}
    b(\bm u,\hat{\bm u},\check{\bm u}) =\int_{\mathcal O}\big[(\bm v\cdot\nabla)\hat{\bm u} +\Phi(\bm v)\partial_z\hat{\bm u}\big]\cdot\check{\bm u} \,\mathrm{d}x\,\mathrm{d}y\,\mathrm{d}z.
\end{equation*}
This form extends continuously to $V\times W\times V$ and $V\times V\times W$. It defines $B:V\times V\to W'$ by duality, and its restrictions satisfy
\begin{equation*}
    B:V\times W\to V',\qquad B:\mathscr D(A)\times\mathscr D(A)\to H.
\end{equation*}
In the latter case, $B(\bm u,\hat{\bm u})=P_H[(\bm v\cdot\nabla)\hat{\bm u}+\Phi(\bm v)\partial_z\hat{\bm u}]$. We write $B(\bm u)=B(\bm u,\bm u)$ and denote its velocity and temperature components by $B_1,B_2$.

\begin{lemma}\label{lem:Euler_bi_op}
    There is $C>0$ with the following properties.
    \begin{enumerate}
        \item For $\bm u\in V$, $\hat{\bm u}\in\mathscr D(A)$ and $\check{\bm u}\in V$,
            \begin{equation}\label{id:B-skew}
                b(\bm u,\hat{\bm u},\check{\bm u})=-b(\bm u,\check{\bm u},\hat{\bm u}),\qquad b(\bm u,\hat{\bm u},\hat{\bm u})=0.
            \end{equation}

        \item For the same arguments,
            \begin{equation}\label{est:B_op-prod-norm}
                \begin{aligned}
                    |b(\bm u,\hat{\bm u},\check{\bm u})|\leq{} & C|\bm u|_H^{1/2}\|\bm u\|_V^{1/2}\|\hat{\bm u}\|_V\|\check{\bm u}\|_V \\
                    & +C\|\bm u\|_V|\partial_z\hat{\bm u}|_2^{1/2}\|\partial_z\hat{\bm u}\|_{H^1}^{1/2}|\check{\bm u}|_H^{1/2}\|\check{\bm u}\|_V^{1/2}.
                \end{aligned}
            \end{equation}
            In particular, for $\bm u\in\mathscr D(A)$,
            \begin{equation}\label{est:B_op-dual-norm}
                \|B(\bm u)\|_{V'}\leq C|\bm u|_H^{1/2}\|\bm u\|_V^{3/2}+C\|\bm u\|_V^{3/2}\|\bm u\|_{\mathscr D(A)}^{1/2}.
            \end{equation}
        \item For $\bm u,\hat{\bm u}\in\mathscr D(A)$ and $\check{\bm u}\in H$,
            \begin{equation}\label{est:B-H-norm}
                \begin{aligned}
                    |\langle B(\bm u,\hat{\bm u}),\check{\bm u}\rangle_H| \leq C\Big( & |\bm v|_6\|\hat{\bm u}\|_V^{1/2}|A\hat{\bm u}|_H^{1/2} \\
                    & +\|\bm u\|_V^{1/2}|A\bm u|_H^{1/2} |\partial_z\hat{\bm u}|_2^{1/2}\|\partial_z\hat{\bm u}\|_{H^1}^{1/2}\Big)|\check{\bm u}|_H.
                \end{aligned}
            \end{equation}
            In particular,
            \begin{equation}\label{est:B-H-symmetric}
                |B(\bm u,\hat{\bm u})|_H \leq c_b\|\bm u\|_V^{1/2}|A\bm u|_H^{1/2} \|\hat{\bm u}\|_V^{1/2}|A\hat{\bm u}|_H^{1/2}.
            \end{equation}
            Also, for $\bm u,\check{\bm u}\in V$ and $\hat{\bm u}\in\mathscr D(A)$,
            \begin{equation}\label{est:B-VDA}
                |b(\bm u,\hat{\bm u},\check{\bm u})| \leq c_b\|\bm u\|_V|A\hat{\bm u}|_H\|\check{\bm u}\|_V.
            \end{equation}
        \item For $\bm u,\hat{\bm u}\in\mathscr D(A)$, $\bm d=\bm u-\hat{\bm u}$ and every $\varepsilon>0$,
            \begin{equation}\label{est:B-diff_prod}
                |\langle B(\bm u)-B(\hat{\bm u}),\bm d\rangle_H| \leq\varepsilon\|\bm d\|_V^2+C_\varepsilon\big(\|\hat{\bm u}\|_V^4+\|\hat{\bm u}\|_V^2\|\hat{\bm u}\|_{\mathscr D(A)}^2\big)|\bm d|_H^2.
            \end{equation}
    \end{enumerate}
\end{lemma}
\begin{proof}
    Since $\partial_z\Phi(\bm v)=-\operatorname{div}\bm v$, the full velocity $(\bm v,\Phi(\bm v))$ is divergence free and has zero normal flux. Integration by parts proves the skew identity. For horizontal convection, H\"older's inequality and \eqref{est:GN-d3-3}--\eqref{est:GN-d3-6} give
    \begin{equation*}
        |b_{\rm hor}(\bm u,\hat{\bm u},\check{\bm u})|\leq |\bm v|_3|\nabla\hat{\bm u}|_2|\check{\bm u}|_6 \leq C|\bm u|_H^{1/2}\|\bm u\|_V^{1/2}\|\hat{\bm u}\|_V\|\check{\bm u}\|_V.
    \end{equation*}
    For any scalar or vector field $f$, set $F_f(x)=(\int_{-h}^0|f(x,z)|^2\,\mathrm{d}z)^{1/2}$. The two-dimensional Ladyzhenskaya inequality and $|\nabla F_f|\leq F_{\nabla f}$ imply
    \begin{equation*}
        \|F_f\|_{L^4(\mathcal D)}\leq C|f|_2^{1/2}(|f|_2+|\nabla f|_2)^{1/2}.
    \end{equation*}
    Since $\sup_z|\Phi(\bm v)|\leq\sqrt h F_{\mathrm{div}\bm v}$,
    \begin{equation*}
        |b_{\rm vert}(\bm u,\hat{\bm u},\check{\bm u})|\leq C\|F_{\mathrm{div}\bm v}\|_2\|F_{\partial_z\hat{\bm u}}\|_4\|F_{\check{\bm u}}\|_4,
    \end{equation*}
    which proves \eqref{est:B_op-prod-norm}, the analogue of \cite[Lemma 2.1, (11)]{1354.35160}. Taking the supremum over $\|\check{\bm u}\|_V\leq1$ proves \eqref{est:B_op-dual-norm}. Alternatively, using exponents $4,4,2$ on these three factors, and estimating horizontal convection by $|\bm v|_6|\nabla\hat{\bm u}|_3|\check{\bm u}|_2$, gives \eqref{est:B-H-norm}, as in \cite[Lemma~2.1, (12)]{1354.35160} and \cite[(2.16)]{BrzezniakSlavik2021}. The hydrostatic projection is a contraction in $L^2$. Since $\|\bm u\|_V\leq C|A\bm u|_H$ and $\|\partial_z\hat{\bm u}\|_{H^1}\leq C|A\hat{\bm u}|_H$, \eqref{est:B-H-norm} gives \eqref{est:B-H-symmetric}. Estimate \eqref{est:B_op-prod-norm} similarly gives \eqref{est:B-VDA}. These are the bounds in \cite[(2.15) and (2.14)]{BrzezniakSlavik2021}, respectively.

    For the difference estimate, skew symmetry gives $\langle B(\bm u)-B(\hat{\bm u}),\bm d\rangle_H=b(\bm d,\hat{\bm u},\bm d)$. The horizontal term is bounded by $C\|\hat{\bm u}\|_V|\bm d|_H^{1/2}\|\bm d\|_V^{3/2}$. The vertical term satisfies
    \begin{equation*}
        \begin{aligned}
            |b_{\rm vert}(\bm d,\hat{\bm u},\bm d)| & \leq C|\partial_z\hat{\bm u}|_2^{1/2}\|\partial_z\hat{\bm u}\|_{H^1}^{1/2}|\bm d|_H^{1/2}\|\bm d\|_V^{3/2} \\
            & \leq\varepsilon\|\bm d\|_V^2+C_\varepsilon|\partial_z\hat{\bm u}|_2^2\|\partial_z\hat{\bm u}\|_{H^1}^2|\bm d|_H^2.
        \end{aligned}
    \end{equation*}
    Young's inequality and $|\partial_z\hat{\bm u}|_2\leq C\|\hat{\bm u}\|_V$, $\|\partial_z\hat{\bm u}\|_{H^1}\leq C\|\hat{\bm u}\|_{\mathscr D(A)}$ prove \eqref{est:B-diff_prod}.
\end{proof}

Define $R:V\to H$ by
\begin{equation*}
    R(\bm u)=\left(P_{H_1}\left[f\bm k\times\bm v-\int_{-h}^z\nabla T\,\mathrm{d}z'\right],0\right).
\end{equation*}
Then $|R(\bm u)|_H\leq C\|\bm u\|_V$, and its associated form is
\begin{align*}
    r(\bm u,\hat{\bm u}) & :={}_{V^\prime}\langle R(\bm u),\hat{\bm u}\rangle_V \\
    & =\int_{\mathcal{O}}(f\bm k\times\bm v)\cdot\hat{\bm v}+(-\int_{-h}^z\nabla T\,\mathrm{d}z^\prime)\cdot\hat{\bm v}\,\mathrm{d}x\,\mathrm{d}y\,\mathrm{d}z.
\end{align*}

\begin{samepage}
    \begin{lemma}
        If $\bm u=(\bm v,T),\hat{\bm u}=(\hat{\bm v},\hat{T})\in V$, then
        \begin{enumerate}
            \item
                \begin{equation}\label{id:r}
                    r(\bm u,\bm u)={}_{V^\prime}\langle R(\bm u),\bm u\rangle_V=-\int_{\mathcal{O}}(\int_{-h}^z\nabla T\,\mathrm{d}z^\prime)\cdot\bm v\,\mathrm{d}x\,\mathrm{d}y\,\mathrm{d}z.
                \end{equation}
            \item
                \begin{equation}\label{est:r}
                    |r(\bm u,\hat{\bm u})|\leq C|\bm v|_{H_1}|\hat{\bm v}|_{H_1} + C(|T|_{H_2}\|\hat{\bm v}\|_{V_1}\land\|T\|_{V_2}|\hat{\bm v}|_{H_1}),
                \end{equation}
                in particular,
                \begin{equation}
                    \begin{aligned}
                        |r(\bm u,\bm u)| & \leq C(|T|_{H_2}\|\bm v\|_{V_1}\land\|T\|_{V_2}|\bm v|_{H_1}) \\
                        & \leq C|\bm u|_H\|\bm u\|_V.
                    \end{aligned}
                \end{equation}
        \end{enumerate}
    \end{lemma}
\end{samepage}

\begin{remark}\label{rk:r}
    The identity \eqref{id:r} follows from $(f\bm k\times\bm v)\cdot\bm v=0$. Moreover, \eqref{id:div_v} and the lateral boundary condition give
    \begin{align*}
        \int_{\mathcal{O}}\nabla p_b\cdot\bm v\,\mathrm{d}x\,\mathrm{d}y\,\mathrm{d}z & =\int_{\mathcal{D}}\nabla p_b\cdot(\int_{-h}^0\bm v\,\mathrm{d}z)\,\mathrm{d}x\,\mathrm{d}y \\
        & =-\int_{\mathcal{D}}p_b(\int_{-h}^0\mathrm{div}\bm v\,\mathrm{d}z)\,\mathrm{d}x\,\mathrm{d}y=0.
    \end{align*}
\end{remark}

Let $(\Omega,\mathcal F,\mathbb F,\mathbb P)$ be a stochastic basis satisfying the usual conditions, where $\mathbb F:=\{\mathcal F_t\}_{t\geq0}$. Let $E$ be a locally compact Polish space with Borel $\sigma$-field $\mathcal E$ and a $\sigma$-finite measure $\nu$. We take $\eta$ to be an $\mathbb F$-Poisson random measure on $[0,\infty)\times E$ with intensity $\mathrm{d}t\otimes\nu$. Its compensated version satisfies
\begin{equation}\label{id:cPoisson}
    \tilde\eta((0,t]\times U)=\eta((0,t]\times U)-t\nu(U), \qquad t\geq0,\quad U\in\mathcal E,\quad\nu(U)<\infty.
\end{equation}

Set $G(\bm u,\xi):=(G_1(\bm u,\xi),G_2(\bm u,\xi))^\top$. The system \eqref{eqn:r-1}--\eqref{eqn:r-3} can be rewritten in the following abstract form
\begin{equation}\label{eqn:abstract}
    \begin{cases}
        \displaystyle\,\mathrm{d}\bm u + (A\bm u + B(\bm u) + R(\bm u))\,\mathrm{d}t = \int_EG(\bm u(t-),\xi)\tilde{\eta}(\mathrm{d}t,\mathrm{d}\xi), \\
        \bm u(0)=\bm u_0\in V.
    \end{cases}
\end{equation}

\subsection{Notions of solution}\label{subsec:solutions}

Fix a stochastic basis and a Poisson random measure as above. We denote by $D(I;X)$ the space of c\`adl\`ag paths from an interval $I$ into a Banach space $X$. When $I$ is clear from the context, we abbreviate $L^p(I;X)$ and $C(I;X)$ as $L^pX$ and $CX$, respectively.

\begin{definition}[Martingale solution]\label{def:martingale-sol}
    Let $\mathcal{T}>0$ and let $\bm u_0$ be an $\mathcal{F}_0$-measurable $V$-valued random variable.  A pair $(\eta,\bm u)$ is called a \emph{martingale solution} of \eqref{eqn:abstract} on $[0,\mathcal{T}]$ if
    \begin{enumerate}
        \item $\bm u$ is an $H$-valued, $\mathbb F$-adapted process with paths in $D([0,\mathcal T];V')\cap L^2(0,\mathcal T;V)$ almost surely; moreover $B(\bm u)\in L^1(0,\mathcal T;V')$ almost surely. The integrand $G(\bm u(s-),\xi)$ is well defined, predictable and locally square integrable in $H$ with respect to $\mathrm{d}s\,\nu(\mathrm{d}\xi)$;
        \item for every $\psi\in V$ the identity
            \begin{equation}\label{eqn:weak-form}
                \begin{aligned}
                    \langle\bm u(t),\psi\rangle_H & +\int_0^t{}_{V^\prime}\langle A\bm u(s)+B(\bm u(s)) +R(\bm u(s)),\psi\rangle_V\,\mathrm{d}s \\
                    & =\langle\bm u_0,\psi\rangle_H +\int_0^t\int_E {}_{V^\prime}\langle G(\bm u(s-),\xi),\psi\rangle_V\, \tilde{\eta}(\mathrm{d}s,\mathrm{d}\xi)
                \end{aligned}
            \end{equation}
            holds for all $t\in[0,\mathcal{T}]$, $\mathbb{P}$-almost surely.
    \end{enumerate}
\end{definition}

\begin{definition}[Uniform bound and the class \texorpdfstring{$\mathfrak C$}{C}]\label{def:unif-bound}
    A martingale solution $(\eta,\bm u)$ on $[0,\mathcal{T}]$ is said to satisfy the \emph{uniform bound} if
    \begin{equation}\label{eq:class-bound}
        \mathbb E\Big[\sup_{t\in[0,\mathcal{T}]}\|\bm u(t)\|^2_V +\int_0^{\mathcal{T}}\|\bm u\|^2_{\mathscr D(A)}\,\mathrm{d}t\Big] <\infty .
    \end{equation}
    We write $\mathfrak{C}([0,\mathcal{T}])$ for the collection of martingale solutions on $[0,\mathcal{T}]$ satisfying \eqref{eq:class-bound}.  The same definition applies verbatim to the solutions of the regularised equation \eqref{eqn:u-hat-main}, that is, of \eqref{eqn:abstract} with $B$ replaced by $\theta(\|\cdot-\bm u^*\|_V)B(\cdot)$; the class so obtained is still denoted by $\mathfrak{C}([0,\mathcal{T}])$.
\end{definition}

\begin{lemma}\label{lem:strong-cadlag}
    Under \eqref{asm:g_moment-V} with $p=1$, a solution in $\mathfrak C([0,\mathcal T])$ has a modification in $D([0,\mathcal T];V)$. The same assertion holds for the regularised equation.
\end{lemma}
\begin{proof}
    By \eqref{est:B-H-symmetric},
    \begin{equation*}
        \int_0^{\mathcal T}|B(\bm u)|_H^2\,\mathrm{d}t \leq C\sup_{t\leq\mathcal T}\|\bm u(t)\|_V^2 \int_0^{\mathcal T}|A\bm u|_H^2\,\mathrm{d}t<\infty \quad\text{almost surely}.
    \end{equation*}
    Write the equation as $\mathrm{d}\bm u+A\bm u\,\mathrm{d}t =\bm f\,\mathrm{d}t+\int_EG\,\tilde\eta$, where $\bm f=-B(\bm u)-R(\bm u)$ (with the cutoff when appropriate). Set $\sigma_m=\inf\{t:\int_0^t|\bm f|_H^2\,\mathrm{d}s\geq m\}\wedge\mathcal T$. The linear equation with forcing $\mathbbm{1}_{[0,\sigma_m]}\bm f$ and noise $\mathbbm{1}_{[0,\sigma_m]}G$ has square-integrable data in $V$, $L^2(0,\mathcal T;H)$ and $L^2(\Omega\times[0,\mathcal T]\times E;V)$. For its spectral approximations, the $V$ energy identity and the Burkholder--Davis--Gundy inequality give
    \begin{equation*}
        \mathbb E\sup_{t\leq\mathcal T}\|\bm z(t)\|_V^2 +\mathbb E\int_0^{\mathcal T}|A\bm z|_H^2\,\mathrm{d}t \leq C\mathbb E\left[\|\bm u_0\|_V^2 +\int_0^{\sigma_m}|\bm f|_H^2\,\mathrm{d}t +\int_0^{\sigma_m}\int_E\|G\|_V^2\,\nu(\mathrm{d}\xi)\,\mathrm{d}t\right].
    \end{equation*}
    Applied to differences of spectral approximations, this estimate gives convergence in the mean-square uniform $V$ norm and hence a $V$-valued c\`adl\`ag limit. Uniqueness for the linear equation identifies this limit with $\bm u$ up to $\sigma_m$. Since $\int_0^{\mathcal T}|\bm f|_H^2<\infty$ almost surely, the assertion follows by letting $m\to\infty$.
\end{proof}

\begin{definition}[Pathwise solutions]\label{def:pathwise-sol}
    Let the stochastic basis, the Poisson random measure $\eta$ and the $\mathcal F_0$-measurable initial datum $\bm u_0$ be fixed. A \emph{pathwise solution} of \eqref{eqn:abstract} on $[0,\mathcal{T}]$ is a martingale solution $(\eta,\bm u)$ on this basis whose solution component $\bm u$ is measurable with respect to $\sigma(\eta,\bm u_0)$. A pair $(\bm u,\tau)$, where $\tau>0$ is an $\mathbb F$-stopping time, is called a \emph{local pathwise solution} if $\bm u$ is $\mathbb F$-adapted and $\sigma(\eta,\bm u_0)$-measurable, $\bm u(\cdot\wedge\tau)$ has $V$-valued c\`adl\`ag paths, $\bm u(0)=\bm u_0$, the integrand $\mathbbm{1}_{(0,\tau]}(s)G(\bm u(s-),\xi)$ is predictable and locally square integrable in $H$ with respect to $\mathrm{d}s\,\nu(\mathrm{d}\xi)$, the weak form \eqref{eqn:weak-form} holds $\mathbb P$-almost surely for all $t\in[0,\tau]$, and
    \begin{equation}\label{eq:local-bound}
        \mathbb E\Big[\sup_{t\in[0,\tau]}\|\bm u(t)\|^2_V +\int_0^\tau\|\bm u\|^2_{\mathscr D(A)}\,\mathrm{d}t\Big] <\infty .
    \end{equation}
    A \emph{maximal pathwise solution} is a pair $(\bm u,\xi)$, where $\xi$ is an $\mathbb F$-stopping time, such that
    \begin{enumerate}
        \item there exists an increasing sequence
            $\{\tau_N\}_{N\geq1}$ of stopping times with $\tau_N\uparrow\xi$ almost surely such that $(\bm u,\tau_N)$ is a local pathwise solution for every $N\in\mathbb N$, and
        \item $\hat\tau\leq\xi$ almost surely for every local pathwise solution
            $(\hat{\bm u},\hat\tau)$ on this stochastic basis, driven by the same $\eta$ and starting from the same $\bm u_0$;
    \end{enumerate}
    it is called \emph{global} if $\xi=+\infty$ almost surely.
\end{definition}

The moment bound in this definition is imposed on each localizing interval. A global maximal solution need not belong to $\mathfrak C([0,\mathcal T])$ on every deterministic interval; uniqueness follows by comparison on the common localizing intervals.

\subsection{Properties of \texorpdfstring{$\overline{\bm v}$}{v} and \texorpdfstring{$\widetilde{\bm v}$}{v}}

Denote the barotropic mode $\overline{\bm v}$ by
\begin{equation*}
    \overline{\bm v}(x,y)=\frac{1}{h}\int_{-h}^0\bm v(x,y,z)\,\mathrm{d}z,\quad \text{in } \mathcal{D},
\end{equation*}
and the baroclinic mode $\widetilde{\bm v}$ by
\begin{equation*}
    \widetilde{\bm v}=\bm v - \overline{\bm v},\quad \text{in } \mathcal{O}.
\end{equation*}
For a scalar or vector field $\psi$, set
\begin{equation}\label{eq:z-average}
    \overline\psi=\frac{1}{h}\int_{-h}^0\psi\,\mathrm{d}z, \qquad\widetilde\psi=\psi-\overline\psi.
\end{equation}
Averaging the convection term and integrating by parts in $z$ gives
\begin{equation}\label{cal:convection}
    \overline{(\bm v\cdot\nabla)\bm v + \Phi(\bm v)\partial_z\bm v}=(\overline{\bm v}\cdot\nabla)\overline{\bm v}+\overline{(\widetilde{\bm v}\cdot\nabla)\widetilde{\bm v}+(\nabla\cdot\widetilde{\bm v})\widetilde{\bm v}},
\end{equation}
which, together with \eqref{bd_cond:o-1} and \eqref{bd_cond:o-2}, gives the vertically averaged form of \eqref{eqn:r-1} on $\mathcal D\times[0,\infty)$:
\begin{empheq}[left=\empheqlbrace]{align}
    & \begin{aligned}\label{eq:barotropic}
        \mathrm{d}\overline{\bm v} + \Big[(\overline{\bm v}\cdot\nabla)\overline{\bm v} & + \overline{(\widetilde{\bm v}\cdot\nabla)\widetilde{\bm v}+(\nabla\cdot\widetilde{\bm v})\widetilde{\bm v}} + f\bm k\times\overline{\bm v} + \nabla p_b \\
        & - \frac{1}{h}\int_{-h}^0\int_{-h}^z\nabla T\,\mathrm{d}z^\prime\,\mathrm{d}z - \Delta\overline{\bm v}\Big]\,\mathrm{d}t= \int_E\overline{G_1}(\bm v(t-),T(t-),\xi)\tilde{\eta}(\mathrm{d}t,\mathrm{d}\xi),
    \end{aligned} \\
    & \mathrm{div}\overline{\bm v}=0,\label{eq:div-v-barotropic}
\end{empheq}
and boundary conditions
\begin{equation}\label{bd_cond-barotropic}
    \overline{\bm v}\cdot\bm n=0,\ \frac{\partial \overline{\bm v}}{\partial\bm n}\times \bm n=0,\quad \text{on} \ \partial\mathcal{D}.
\end{equation}
Subtracting \eqref{eq:barotropic} from \eqref{eqn:r-1} gives the baroclinic equation on $\mathcal O\times[0,\infty)$:
\begin{equation}\label{eq:baroclinic}
    \begin{aligned}
        \displaystyle\,\mathrm{d}\widetilde{\bm v} + ((\widetilde{\bm v}\cdot\nabla)\widetilde{\bm v} & + \Phi(\widetilde{\bm v})\partial_z\widetilde{\bm v} + (\overline{\bm v}\cdot\nabla)\widetilde{\bm v} + (\widetilde{\bm v}\cdot\nabla)\overline{\bm v} - \overline{(\widetilde{\bm v}\cdot\nabla)\widetilde{\bm v}+(\nabla\cdot\widetilde{\bm v})\widetilde{\bm v}} \\
        & + f\bm k\times\widetilde{\bm v} -\int_{-h}^z\nabla T\,\mathrm{d}z^\prime+\frac{1}{h}\int_{-h}^0\int_{-h}^z\nabla T\,\mathrm{d}z^\prime\,\mathrm{d}z + L_1\widetilde{\bm v})\,\mathrm{d}t \\
        = & \int_E(\widetilde{G_1}(\bm v(t-),T(t-),\xi))\tilde{\eta}(\mathrm{d}t,\mathrm{d}\xi),
    \end{aligned}
\end{equation}
and boundary conditions
\begin{align}
    \partial_z\widetilde{\bm v}=0,\quad & \text{on} \ \Gamma_u, \\
    \partial_z\widetilde{\bm v}=0,\quad & \text{on} \ \Gamma_b, \\
    \widetilde{\bm v}\cdot\bm n=0,\ \frac{\partial \widetilde{\bm v}}{\partial\bm n}\times\bm n=0\quad & \text{on} \ \Gamma_l\label{bd:baroclinic-3}.
\end{align}

\subsection{Some inequalities}
We use the following interpolation inequalities \cite{0687.35071,0588.35003}. For $\phi\in H^1(\mathcal D)$,
\begin{align}
    |\phi|_4 & \leq C|\phi|_2^{1/2}\|\phi\|_{H^1(\mathcal D)}^{1/2}, \label{est:GN-d2-4} \\
    |\phi|_8 & \leq C|\phi|_6^{3/4}\|\phi\|_{H^1(\mathcal D)}^{1/4}. \label{est:GN-d2-8}
\end{align}
For $\phi\in H^1(\mathcal O)$,
\begin{align}
    |\phi|_3 & \leq C|\phi|_2^{1/2}\|\phi\|_{H^1(\mathcal O)}^{1/2}, \label{est:GN-d3-3} \\
    |\phi|_6 & \leq C\|\phi\|_{H^1(\mathcal O)}. \label{est:GN-d3-6}
\end{align}
For $\phi\in H^2(\mathcal O)$, we also have
\begin{equation}\label{est:GN-d3-10}
    |\phi|_{10}^5\leq C\|\phi\|_{H^1(\mathcal O)}^4 \|\phi\|_{H^2(\mathcal O)}.
\end{equation}
For $\phi\in H^1(\mathcal D)$ with $|\phi|^3\in H^1(\mathcal D)$, applying \eqref{est:GN-d2-4} to $|\phi|^3$ gives
\begin{equation}\label{est:L12}
    |\phi|_{L^{12}(\mathcal D)}^{12} \leq C|\phi|_{L^6(\mathcal D)}^6 \int_{\mathcal D}|\phi|^4|\nabla\phi|^2\,\mathrm{d}x\,\mathrm{d}y +C|\phi|_{L^6(\mathcal D)}^{12}.
\end{equation}

We use the following localized Gr\"onwall lemma; see \cite[Proposition~2.1]{BrzezniakSlavik2021}.

\begin{lemma}[Stochastic Gr\"onwall inequality]\label{lem:sto_gron}
    Let $X$ be a nonnegative adapted c\`adl\`ag process and let $Y,Z,R$ be nonnegative progressively measurable processes. Suppose $\tau\leq\mathcal T$ is a stopping time and, for a deterministic $K<\infty$,
    \begin{equation}\label{eq:gron-hyp}
        \mathbb EX(0)+\mathbb E\int_0^\tau(RX+Z)\,\mathrm{d}s<\infty, \qquad\int_0^\tau R\,\mathrm{d}s\leq K\quad\text{almost surely}.
    \end{equation}
    Assume that, for all stopping times $0\leq\tau_a\leq\tau_b\leq\tau$,
    \begin{equation}\label{eq:gron-local}
        \mathbb E\left[\sup_{\tau_a\leq t\leq\tau_b}X(t)+\int_{\tau_a}^{\tau_b}Y\,\mathrm{d}s\right] \leq C_0\mathbb E\left[X(\tau_a)+\int_{\tau_a}^{\tau_b}(RX+Z)\,\mathrm{d}s\right].
    \end{equation}
    Then, for a constant $C=C(C_0,K)$,
    \begin{equation}\label{eq:gron-concl}
        \mathbb E\left[\sup_{t\leq\tau}X(t)+\int_0^\tau Y\,\mathrm{d}s\right] \leq C\mathbb E\left[X(0)+\int_0^\tau Z\,\mathrm{d}s\right].
    \end{equation}
\end{lemma}
\begin{proof}
    Starting from $\tau_0=0$, partition $[0,\tau]$ at successive increments $(2C_0)^{-1}$ of $\int_0^tR\,\mathrm{d}s$. Each interval satisfies $\int_{\tau_j}^{\tau_{j+1}}R\,\mathrm{d}s\leq(2C_0)^{-1}$, so absorption in \eqref{eq:gron-local} gives
    \begin{equation*}
        \mathbb E\left[\sup_{\tau_j\leq t\leq\tau_{j+1}}X(t) +\int_{\tau_j}^{\tau_{j+1}}Y\,\mathrm{d}s\right] \leq2C_0\mathbb E\left[X(\tau_j)+\int_{\tau_j}^{\tau_{j+1}}Z\,\mathrm{d}s\right].
    \end{equation*}
    At most $\lceil2C_0K\rceil+1$ intervals are required. Iteration proves \eqref{eq:gron-concl}, with a constant depending on $C_0$ and $K$.
\end{proof}

    \section{Existence of martingale solutions}\label{sec:wp-martingale}
We use a Galerkin approximation with a cutoff around the linear evolution of the initial datum, following \cite{BrzezniakSlavik2021}. We first obtain uniform estimates, then prove tightness and identify the limiting martingale problem.
\begin{assumption}\label{asm:martingale}
    For $i=1,2$, assume that $G_i:V\times E\to V_i$ is a measurable mapping.
    \begin{enumerate}
        \item (Linear growth in $H$) For $p=1,2$ and $(\bm v,T)\in V$,
            \begin{equation}\label{asm:g_linear-H}
                \begin{aligned}
                    \int_E|G_1(\bm v,T,\xi)|_{H_1}^{2p}\nu(\mathrm{d}\xi) & \leq C(1+|\bm v|_{H_1}^{2p}), \\
                    \int_E|G_2(\bm v,T,\xi)|_{H_2}^{2p}\nu(\mathrm{d}\xi) & \leq C(1+|T|_{H_2}^{2p}).
                \end{aligned}
            \end{equation}
        \item (Local Lipschitz in $H$ on bounded subsets of $V$) For every $r>0$, there exists $C_r>0$ such that, for $\bm u,\hat{\bm u}\in V$ with $\|\bm u\|_V\vee\|\hat{\bm u}\|_V\leq r$,
            \begin{equation}\label{asm:g_Lip-local}
                \int_E|G(\bm u,\xi)-G(\hat{\bm u},\xi)|_H^2\,\nu(\mathrm{d}\xi) \leq C_r|\bm u-\hat{\bm u}|_H^2.
            \end{equation}
    \end{enumerate}
\end{assumption}

For the strong PDE estimates we also impose the following moment bounds.

\begin{assumption}\label{asm:A-est}
    In addition to Assumption \ref{asm:martingale}, assume that, for $p=1,2$ and $\bm u\in V$,
    \begin{equation}\label{asm:g_moment-V}
        \int_E\|G(\bm u,\xi)\|_V^{2p}\,\nu(\mathrm{d}\xi)\leq C_p\big(1+\|\bm u\|_V^{2p}\big).
    \end{equation}
\end{assumption}

\begin{lemma}\label{lem:local-noise-convergence}
    Suppose that $G:V\times E\to V$ is measurable and satisfies \eqref{asm:g_linear-H} with $p=1$ and \eqref{asm:g_Lip-local}. Let $v_n,v$ be jointly measurable $V$-valued processes such that
    \begin{equation*}
        \mathbb E\int_0^{\mathcal T}|v_n-v|_H^2\,\mathrm{d}s\to0,\qquad
        \sup_n\mathbb E\int_0^{\mathcal T}\|v_n\|_V^4\,\mathrm{d}s+\mathbb E\int_0^{\mathcal T}\|v\|_V^4\,\mathrm{d}s<\infty.
    \end{equation*}
    Then $G(v_n,\cdot)\to G(v,\cdot)$ in $L^2(\Omega\times(0,\mathcal T)\times E;H)$.
\end{lemma}
\begin{proof}
    Fix $r\geq1$ and set $B_{n,r}=\{\|v_n\|_V\vee\|v\|_V\leq r\}$. On this set, \eqref{asm:g_Lip-local} bounds the integral of the squared coefficient difference over $E$ by $C_r|v_n-v|_H^2$. On its complement, \eqref{asm:g_linear-H} with $p=1$ and the embedding $V\hookrightarrow H$ give a bound by $C(1+\|v_n\|_V^2+\|v\|_V^2)$. The elementary inequality
    \begin{equation*}
        (1+a^2+b^2)\mathbbm1_{\{a\vee b>r\}}\leq\frac{3}{r^2}(a^4+b^4),\qquad a,b\geq0,\quad r\geq1,
    \end{equation*}
    therefore yields
    \begin{equation}\label{est:local-noise-tail}
        \begin{split}
            \mathbb E\int_0^{\mathcal T}\int_E|G(v_n,\xi)-G(v,\xi)|_H^2\,\nu(\mathrm{d}\xi)\,\mathrm{d}s
            &\leq C_r\mathbb E\int_0^{\mathcal T}|v_n-v|_H^2\,\mathrm{d}s\\
            &\quad+\frac{C}{r^2}\mathbb E\int_0^{\mathcal T}(\|v_n\|_V^4+\|v\|_V^4)\,\mathrm{d}s.
        \end{split}
    \end{equation}
    For fixed $r$, the first term on the right tends to zero by the assumed convergence in $L^2(\Omega\times(0,\mathcal T);H)$. The uniform fourth-moment bound controls the second term by $C/r^2$, with $C$ independent of $n$ and $r$. Taking the upper limit as $n\to\infty$ in \eqref{est:local-noise-tail} and then letting $r\to\infty$ proves the claimed convergence.
\end{proof}

\begin{theorem}\label{thm:exi-martingale}
    Assume that Assumptions \ref{asm:martingale} and \ref{asm:A-est} hold and that $\bm u_0\in L^4(\Omega,\mathcal F_0;V)$. Then equation \eqref{eqn:abstract} admits a local martingale solution, and the regularised equation \eqref{eqn:u-hat-main}, in which the convection term carries the cutoff, admits a global martingale solution; both satisfy the bounds of Proposition \ref{prop:uniform-DA}. More precisely, there exist a filtered probability space $(\hat\Omega,\hat{\mathcal F},\hat{\mathbb F},\hat{\mathbb P})$ with $\hat{\mathbb F}:=\{\hat{\mathcal F}_t\}_{t\geq0}$, an $\hat{\mathbb F}$-Poisson random measure $\hat\eta$ on $[0,\infty)\times E$ with intensity measure $\mathrm{d}t\otimes\nu$, and an $\hat{\mathcal F}_0$-measurable $V$-valued initial datum $\hat{\bm u}_0$ such that
    \begin{equation*}
        \mathcal L_{\hat{\mathbb P}}(\hat{\bm u}_0)=\mathcal L_{\mathbb P}(\bm u_0).
    \end{equation*}
    Write $\tilde{\hat\eta}=\hat\eta-\mathrm{d}t\otimes\nu$ and $\hat{\bm u}^*(t)=e^{-tA}\hat{\bm u}_0$. On this basis there exists an adapted process $\hat{\bm u}$ with the following properties:
    \begin{enumerate}
        \item The process satisfies $\hat{\bm u}(0)=\hat{\bm u}_0$ and has paths in $D([0,\infty);V)\cap L_{\mathrm{loc}}^2([0,\infty);\mathscr D(A))$ almost surely. For every $\mathcal T>0$,
            \begin{equation}\label{est:u-hat-strong}
                \hat{\mathbb E}\Big[\sup_{t\in[0,\mathcal{T}]}\|\hat{\bm u}(t)\|_V^2+\int_0^{\mathcal{T}}\|\hat{\bm u}\|^2_{\mathscr D(A)}\,\mathrm{d}t\Big]<\infty .
            \end{equation}
        \item For any $\psi\in V$, the following equation holds $\hat{\mathbb{P}}$\text{-}a.s. for all $t\geq0$:
            \begin{equation}\label{eqn:u-hat-main}
                \begin{aligned}
                    \langle\hat{\bm{u}}(t),\psi\rangle_H & +\int_0^t{}_{V^\prime}\langle A\hat{\bm{u}}(s)+\theta\big(\|\hat{\bm u}(s)-\hat{\bm u}^*(s)\|_V\big)B(\hat{\bm u}(s))+R(\hat{\bm u}(s)),\psi\rangle_V\,\mathrm{d}s \\
                    = & \langle\hat{\bm u}_0,\psi\rangle_H+\int_0^t\int_E{}_{V^\prime}\langle G(\hat{\bm u}(s-),\xi),\psi\rangle_V\,\tilde{\hat{\eta}}(\mathrm{d}s,\mathrm{d}\xi).
                \end{aligned}
            \end{equation}
        \item With $\tau$ as in \eqref{def:tau-exit}, the stopped process $\hat{\bm u}(\cdot\wedge\tau)$ satisfies the weak form \eqref{eqn:weak-form} of \eqref{eqn:abstract} on the random interval $[0,\tau]$; that is, \eqref{eqn:abstract} admits a local martingale solution.
    \end{enumerate}
\end{theorem}

The local conclusion extends to $\bm u_0\in L^2(\Omega,\mathcal F_0;V)$ by Proposition \ref{prop:initial-localization}. Throughout this section, we assume $\bm u_0\in L^4(\Omega,\mathcal F_0;V)$.

\subsection{Galerkin approximation}\label{subsec:galerkin}

\subsubsection{The regularised scheme}

We truncate the convection term outside a fixed $V$-neighbourhood of the linear evolution of the initial datum, as in \cite[Section~3]{BrzezniakSlavik2021}.

Let $\{e_i\}_{i\geq1}\subset\mathscr D(A)$ be an orthonormal basis of $H$ consisting of eigenvectors of $A$, $Ae_i=\lambda_ie_i$, with $0<\lambda_1\leq\lambda_2\leq\cdots$ and $\lambda_i\to\infty$; such a basis exists by Lemma \ref{lem:Stokes_op}. We choose the basis from the component eigenbases of $A_1$ and $A_2$, so its spectral projections preserve the velocity and temperature components. For $n\geq1$ put $H_n:=\mathrm{span}\{e_1,\dots,e_n\}\subset\mathscr D(A)$ and let $P_n:V^\prime\to H_n$ be the projection defined by
\begin{equation}
    P_n(\phi)=\sum_{i=1}^n{}_{V^\prime}\langle\phi,e_i\rangle_Ve_i, \qquad\phi\in V^\prime.
\end{equation}
The restriction of $P_n$ to $H$ is the $H$-orthogonal projection onto $H_n$, and, since the $e_i$ are eigenvectors of $A$,
\begin{equation}
    A P_n\phi=P_nA\phi\quad\text{for every }\phi\in\mathscr D(A).
\end{equation}
Let $\kappa>0$ be a constant, to be fixed in the proof of Proposition \ref{prop:uniform-DA}, and let $\theta\in C^\infty(\mathbb{R};[0,1])$ be such that
\begin{equation}\label{eq:cutoff}
    \mathbbm{1}_{[-\kappa/2,\kappa/2]}\leq\theta\leq\mathbbm{1}_{[-\kappa,\kappa]}.
\end{equation}
For every $n\geq1$ let $\bm u_n^*$ denote the solution of the linear problem
\begin{equation}\label{eqn:background}
    \frac{\mathrm{d}}{\mathrm{d}t}\bm u_n^*(t)+A\bm u_n^*(t)=0, \qquad \bm u_n^*(0)=P_n\bm u_0.
\end{equation}
Thus $\bm u_n^*(t)=P_ne^{-tA}\bm u_0$, and for every $p\geq2$ and $\mathcal T>0$,
\begin{equation}\label{est:background}
    \sup_{t\in[0,\mathcal T]}\|\bm u_n^*(t)\|_V^p +\int_0^{\mathcal T}\|\bm u_n^*\|_V^{p-2} \|\bm u_n^*\|_{\mathscr D(A)}^2\,\mathrm{d}t +\Big(\int_0^{\mathcal T}\|\bm u_n^*\|^2_{\mathscr D(A)}\,\mathrm{d}t\Big)^{p/2} \leq C_p\|\bm u_0\|_V^p,
\end{equation}
with $C_p$ independent of $n$.

For $\psi\in H_n$ we consider the finite-dimensional stochastic system
\begin{equation}\label{eqn:galerkin}
    \begin{aligned}
        \langle\bm u_n(t),\psi\rangle_H & +\int_0^t\Big\langle A\bm u_n +\theta\big(\|\bm u_n-\bm u_n^*\|_V\big)B_{(n)}(\bm u_n) +R_{(n)}(\bm u_n),\psi\Big\rangle_H\mathrm{d}s \\
        & =\langle P_n\bm u_0,\psi\rangle_H +\int_0^t\int_E \langle G_{(n)}(\bm u_n(s-),\xi),\psi\rangle_H \tilde{\eta}(\mathrm{d}s,\mathrm{d}\xi),
    \end{aligned}
\end{equation}
where the projected operators are defined by
\begin{equation}
    B_{(n)}(\bm u_n):=P_nB(\bm u_n),\quad R_{(n)}(\bm u_n):=P_nR(\bm u_n),\quad G_{(n)}(\bm u_n,\xi):=P_nG(\bm u_n,\xi).
\end{equation}

For fixed $n$, equivalence of the $H$ and $V$ norms on $H_n$ and \eqref{asm:g_Lip-local} give local Lipschitz continuity of the noise coefficient. After localization of the initial datum, the cutoff drift is also locally Lipschitz, uniformly on bounded time intervals. The finite-dimensional local construction therefore gives a unique maximal solution; see \cite[Theorem~117]{1070.60002}. Let $\zeta_{n,R}$ be its first exit time from the $H$-ball of radius $R$. The stopped second-moment estimate proved below gives
\begin{equation*}
    \mathbb E\sup_{t\leq\mathcal T\wedge\zeta_{n,R}}|\bm u_n(t)|_H^2\leq C_{\mathcal T}(1+\mathbb E|P_n\bm u_0|_H^2),
\end{equation*}
with a constant independent of $R$, including the terminal jump. Consequently $\mathbb P(\zeta_{n,R}\leq\mathcal T)\leq C_{\mathcal T}(1+\mathbb E|P_n\bm u_0|_H^2)/R^2\to0$. Thus the solution is global and belongs to $D([0,\mathcal T];H_n)$.

\subsubsection{Uniform \texorpdfstring{$L^2$}{L2} estimates}

\begin{lemma}
    Assume Assumption \ref{asm:martingale} holds. For $n\geq 1$, let $\bm u_n$ be the solution of the equation \eqref{eqn:galerkin}. Then, for $\mathcal{T}\geq 0 $ and $p=1,2$, the following uniform estimate holds:
    \begin{equation}\label{est:galerkin-uniform}
        \sup_{n\geq1}\mathbb{E}\left(\sup_{t\in[0,\mathcal{T}]}|\bm u_n(t)|_H^{2p}+\int_0^{\mathcal{T}}|\bm u_n(t)|_H^{2(p-1)}\|\bm u_n(t)\|_V^2\,\mathrm{d}t\right)\leq C(\mathcal{T}).
    \end{equation}
\end{lemma}
\begin{proof}
    We suppress the Galerkin index. All stochastic integrands are evaluated at left limits, and the estimates are first applied after localization. For $m=2,4$, set $F_m(z)=|z|^m$ and
    \begin{equation*}
        Q_m(z,g)=|z+g|^m-|z|^m-m|z|^{m-2}(z,g).
    \end{equation*}
    Convexity and Taylor's formula give
    \begin{equation*}
        0\leq Q_m(z,g)\leq C_m\big(|z|^{m-2}|g|^2+|g|^m\big).
    \end{equation*}
    The spectral projections preserve the two components and are contractions in $H_i$. Thus \eqref{asm:g_linear-H} applies to the projected noise. The convection cancels in each component energy, even with the cutoff. Moreover, by \eqref{id:r}--\eqref{est:r},
    \begin{equation*}
        m|\bm u|_H^{m-2}|r(\bm u,\bm u)| \leq\frac m2|\bm u|_H^{m-2}\|\bm u\|_V^2+C_m|\bm u|_H^m.
    \end{equation*}
    It\^o's formula for $|\bm u|_H^m$ therefore yields
    \begin{equation}\label{eq:low-energy-ito}
        \begin{split}
            |\bm u(t)|_H^m+\frac m2\int_0^t |\bm u|_H^{m-2}\|\bm u\|_V^2\,\mathrm{d}s \leq |\bm u_0|_H^m+C_m\int_0^t|\bm u|_H^m\,\mathrm{d}s +M_m(t)+\int_0^t\int_E Q_m(\bm u,G)\,\eta(\mathrm{d}s,\mathrm{d}\xi),
        \end{split}
    \end{equation}
    where $M_m=m\int_0^\cdot\int_E|\bm u|_H^{m-2} \langle\bm u,G\rangle_H\,\tilde\eta(\mathrm{d}s,\mathrm{d}\xi)$. The Burkholder--Davis--Gundy inequality and compensation imply, for every $\delta>0$,
    \begin{align*}
        \mathbb E\sup_{s\leq t}|M_m(s)| & \leq C_m\mathbb E\left[ \sup_{s\leq t}|\bm u(s)|_H^{m/2} \left(\int_0^t\int_E|\bm u(s-)|_H^{m-2}|G|_H^2 \,\eta(\mathrm{d}s,\mathrm{d}\xi)\right)^{1/2}\right] \\
        & \leq\delta\mathbb E\sup_{s\leq t}|\bm u(s)|_H^m +C_{m,\delta}\int_0^t\big(1+\mathbb E|\bm u(s)|_H^m\big)\,\mathrm{d}s, \\
        \mathbb E\int_0^t\int_E Q_m(\bm u,G)\,\eta(\mathrm{d}s,\mathrm{d}\xi) & \leq C_m\int_0^t\big(1+\mathbb E|\bm u(s)|_H^m\big)\,\mathrm{d}s.
    \end{align*}
    Estimate the supremum and the terminal dissipation separately in \eqref{eq:low-energy-ito}, then add the bounds and choose $\delta$ small. This gives
    \begin{equation*}
        \mathbb E\left[\sup_{s\leq t}|\bm u(s)|_H^m +\int_0^t|\bm u|_H^{m-2}\|\bm u\|_V^2\,\mathrm{d}s\right] \leq C_m\mathbb E|\bm u_0|_H^m +C_m\int_0^t\left(1+\mathbb E\sup_{r\leq s}|\bm u(r)|_H^m\right)\,\mathrm{d}s.
    \end{equation*}
    Gr\"onwall's inequality and Fatou's lemma remove the localization and prove \eqref{est:galerkin-uniform}, with $m=2p$.

    For later use we record the component estimates. The temperature equation has no $R$ term, and its linear martingale satisfies
    \begin{equation}\label{est:G2-compensated}
        \mathbb E\sup_{s\leq t}\left| \int_0^s\int_E\langle T(r-),G_2(\bm u(r-),\xi)\rangle_{H_2} \,\tilde\eta(\mathrm{d}r,\mathrm{d}\xi)\right| \leq\delta\mathbb E\sup_{s\leq t}|T(s)|_{H_2}^2 +C_\delta\int_0^t(1+\mathbb E|T(s)|_{H_2}^2)\,\mathrm{d}s.
    \end{equation}
    The same argument for $|T|_{H_2}^m$ yields
    \begin{align}
        \mathbb E\left[\sup_{t\leq\mathcal T}|T(t)|_{H_2}^2 +\int_0^{\mathcal T}\|T\|_{V_2}^2\,\mathrm{d}t\right] & \leq C_{\mathcal T}(1+\mathbb E|T_0|_{H_2}^2),\label{est:T-L2} \\
        \mathbb E\left[\sup_{t\leq\mathcal T}|T(t)|_{H_2}^4 +\int_0^{\mathcal T}|T|_{H_2}^2\|T\|_{V_2}^2\,\mathrm{d}t\right] & \leq C_{\mathcal T}(1+\mathbb E|T_0|_{H_2}^4).\label{est:T-L4L2}
    \end{align}
    For the velocity equation, the drift term is $-m|\bm v|_{H_1}^{m-2} r(\bm u,\bm u)$, and
    \begin{equation*}
        m|\bm v|_{H_1}^{m-2}|r(\bm u,\bm u)| \leq\frac m2|\bm v|_{H_1}^{m-2}\|\bm v\|_{V_1}^2 +C_m\big(|\bm v|_{H_1}^m+|T|_{H_2}^m\big).
    \end{equation*}
    Using the temperature bounds, compensation and Gr\"onwall's inequality gives
    \begin{align}
        \mathbb E\left[\sup_{t\leq\mathcal T}|\bm v(t)|_{H_1}^2 +\int_0^{\mathcal T}\|\bm v\|_{V_1}^2\,\mathrm{d}t\right] & \leq C_{\mathcal T}(1+\mathbb E|\bm u_0|_H^2),\label{est:v-L2} \\
        \mathbb E\left[\sup_{t\leq\mathcal T}|\bm v(t)|_{H_1}^4 +\int_0^{\mathcal T}|\bm v|_{H_1}^2\|\bm v\|_{V_1}^2\,\mathrm{d}t\right] & \leq C_{\mathcal T}(1+\mathbb E|\bm u_0|_H^4).\label{est:v-L4L2}
    \end{align}
    All constants are independent of the Galerkin dimension.
\end{proof}

\subsubsection{Uniform estimates in \texorpdfstring{$\mathscr D(A)$}{D(A)}}\label{ssubsec:DA}

Set $\bar{\bm u}_n=\bm u_n-\bm u_n^*$. We estimate the weighted dissipation
\begin{equation}\label{eq:weighted-dissipation}
    \int_0^{\mathcal T}\|\bar{\bm u}_n\|_V^{p-2}|A\bar{\bm u}_n|_H^2\,\mathrm{d}t.
\end{equation}
The nonlinear estimates follow the decomposition in \cite[(3.13)--(3.15)]{BrzezniakSlavik2021}.

\begin{proposition}[Uniform $\mathscr D(A)$ estimates]\label{prop:uniform-DA}
    Let Assumptions \ref{asm:martingale} and \ref{asm:A-est} hold, and let $\bm u_0\in L^q(\Omega,\mathcal F_0;V)$ with $q\geq4$. There is $\kappa>0$ such that, for every $p\in[2,4]$,
    \begin{equation}\label{est:uniform-DA-p}
        \mathbb E\bigg[\sup_{t\leq\mathcal T}\|\bar{\bm u}_n(t)\|_V^p +\int_0^{\mathcal T}\|\bar{\bm u}_n\|_V^{p-2}|A\bar{\bm u}_n|_H^2\,\mathrm{d}t\bigg] \leq C_p\big(1+\mathbb E\|\bm u_0\|_V^4\big),
    \end{equation}
    uniformly in $n$. If also $2p\leq q$, then
    \begin{equation}\label{est:uniform-DA-dissipation}
        \mathbb E\bigg(\int_0^{\mathcal T}|A\bar{\bm u}_n|_H^2\,\mathrm{d}t\bigg)^{p/2} \leq C_p\big(1+\mathbb E\|\bm u_0\|_V^{2p}\big).
    \end{equation}
    In particular,
    \begin{equation}\label{est:uniform-DA}
        \sup_n\mathbb E\bigg[\sup_{t\leq\mathcal T}\|\bar{\bm u}_n(t)\|_V^2 +\int_0^{\mathcal T}|A\bar{\bm u}_n|_H^2\,\mathrm{d}t\bigg]<\infty.
    \end{equation}
    Writing $M_n(t)=\int_0^t\int_EG_{(n)}(\bm u_n(s-),\xi)\,\tilde\eta(\mathrm{d}s,\mathrm{d}\xi)$, we also have
    \begin{equation}\label{est:uniform-DA-moments}
        \sup_n\mathbb E\|M_n\|_{W^{\alpha,p}(0,\mathcal T;H)}^p<\infty, \qquad 0<\alpha<1/p.
    \end{equation}
\end{proposition}

\begin{proof}
    We suppress $n$ and write $\bm u=\bar{\bm u}+\bm u^*$, $x=\|\bar{\bm u}\|_V$, $y=|A\bar{\bm u}|_H$, $a=\|\bm u^*\|_V$, and $d=|A\bm u^*|_H$. All stochastic integrands below are evaluated at left limits. The spectral projections are contractions on $V$ and $H$. For $F_p(\bm w)=\|\bm w\|_V^p$, let
    \begin{equation*}
        Q_p(\bm w,\bm g)=F_p(\bm w+\bm g)-F_p(\bm w)-DF_p(\bm w)\bm g.
    \end{equation*}
    Convexity and Taylor's formula give
    \begin{equation}\label{est:Psi}
        0\leq Q_p(\bm w,\bm g) \leq C_p\big(\|\bm w\|_V^{p-2}\|\bm g\|_V^2+\|\bm g\|_V^p\big).
    \end{equation}
    The finite-dimensional jump It\^o formula yields
    \begin{equation}\label{eq:ito-DA}
        x(t)^p+p\int_0^t x^{p-2}y^2\,\mathrm{d}s =J_1(t)+J_2(t)+\mathcal M_p(t)+\mathcal Q_p(t),
    \end{equation}
    where
    \begin{align*}
        J_1(t) & =-p\int_0^t\theta(x)x^{p-2}\langle B(\bm u),A\bar{\bm u}\rangle_H\,\mathrm{d}s, \\
        J_2(t) & =-p\int_0^t x^{p-2}\langle R(\bm u),A\bar{\bm u}\rangle_H\,\mathrm{d}s, \\
        \mathcal M_p(t) & =p\int_0^t\int_E x(s-)^{p-2} \langle\bar{\bm u}(s-),G_{(n)}(\bm u(s-),\xi)\rangle_V \,\tilde\eta(\mathrm{d}s,\mathrm{d}\xi), \\
        \mathcal Q_p(t) & =\int_0^t\int_E Q_p\big(\bar{\bm u}(s-),G_{(n)}(\bm u(s-),\xi)\big) \,\eta(\mathrm{d}s,\mathrm{d}\xi).
    \end{align*}
    Here $\bar{\bm u}(0)=0$; for $p=2$ the factor $x^{p-2}$ is understood to be $1$.

    By bilinearity and \eqref{est:B-H-symmetric},
    \begin{equation*}
        \begin{aligned}
            p\theta(x)x^{p-2}|\langle B(\bar{\bm u}),A\bar{\bm u}\rangle_H| & \leq pc_b\kappa x^{p-2}y^2, \\
            p\theta(x)x^{p-2}|\langle B(\bm u^*),A\bar{\bm u}\rangle_H| & \leq\varepsilon x^{p-2}y^2+C_{p,\varepsilon}\kappa^{p-2}a^2d^2, \\
            p\theta(x)x^{p-2}|\langle B(\bar{\bm u},\bm u^*)+B(\bm u^*,\bar{\bm u}),A\bar{\bm u}\rangle_H| & \leq\varepsilon x^{p-2}y^2+C_{p,\varepsilon}\kappa^pa^2d^2.
        \end{aligned}
    \end{equation*}
    Indeed, before Young's inequality the mixed terms are bounded by $C_p\theta(x)x^{p-3/2}y^{3/2}a^{1/2}d^{1/2}$. Choose $c_b\kappa\leq1/8$ and then $\varepsilon$ sufficiently small. Since $|R(\bm u)|_H\leq C\|\bm u\|_V$, we obtain
    \begin{equation}\label{est:cutoff-B-weighted}
        |J_1(t)|+|J_2(t)| \leq\frac p2\int_0^t x^{p-2}y^2\,\mathrm{d}s +C_{p,\kappa}\int_0^t\big(1+x^p+a^p+a^2d^2\big)\,\mathrm{d}s.
    \end{equation}
    The background estimate implies
    \begin{equation}\label{eq:int-b}
        \int_0^{\mathcal T}a^2d^2\,\mathrm{d}s\leq C\|\bm u_0\|_V^4.
    \end{equation}

    H\"older interpolation between the second- and fourth-moment bounds in \eqref{asm:g_moment-V} yields
    \begin{equation*}
        \int_E\|G(\bm u,\xi)\|_V^p\,\nu(\mathrm{d}\xi)\leq C_p(1+\|\bm u\|_V^p),\qquad 2\leq p\leq4.
    \end{equation*}
    Combining this bound with compensation and \eqref{est:Psi}, we obtain
    \begin{equation*}
        \mathbb E\mathcal Q_p(t) \leq C_p\mathbb E\int_0^t \big[x^{p-2}(1+x^2+a^2)+1+x^p+a^p\big]\,\mathrm{d}s.
    \end{equation*}
    For $p>2$, $x^{p-2}a^2\leq (p-2)x^p/p+2a^p/p$ and $x^{p-2}\leq1+x^p$; the case $p=2$ is immediate. Hence
    \begin{equation}\label{est:J4}
        \mathbb E\mathcal Q_p(t) \leq C_p\mathbb E\int_0^t\big(1+x^p+a^p\big)\,\mathrm{d}s.
    \end{equation}
    The Burkholder--Davis--Gundy inequality, followed by Young's inequality and compensation, gives, for any $\delta>0$,
    \begin{equation}\label{est:J3-p2}
        \begin{aligned}
            \mathbb E\sup_{s\leq t}|\mathcal M_p(s)| & \leq C_p\mathbb E\left[\sup_{s\leq t}x(s)^{p/2} \left(\int_0^t\int_E x(s-)^{p-2}\|G_{(n)}\|_V^2 \,\eta(\mathrm{d}s,\mathrm{d}\xi)\right)^{1/2}\right] \\
            & \leq\delta\mathbb E\sup_{s\leq t}x(s)^p +C_{p,\delta}\mathbb E\int_0^t(1+x^p+a^p)\,\mathrm{d}s.
        \end{aligned}
    \end{equation}
    First stop when $x$ or $\int_0^t y^2\,\mathrm{d}s$ exceeds an integer $m$, and restrict the initial datum to $\{\|\bm u_0\|_V\leq m\}$. Before this stop the pre-jump integrands have the required moments; compensation controls the terminal jump. All constants in the following estimate are independent of $m$. For sufficiently small $\delta$,
    \begin{equation}\label{est:before-gronwall-DA}
        \begin{aligned}
            \mathbb E\left[\sup_{s\leq t}x(s)^p+\int_0^t x^{p-2}y^2\,\mathrm{d}s\right] & \leq C_{p,\kappa}\big(1+\mathbb E\|\bm u_0\|_V^4\big) \\
            & \quad+C_{p,\kappa}\int_0^t\mathbb E\sup_{r\leq s}x(r)^p\,\mathrm{d}s.
        \end{aligned}
    \end{equation}
    Gr\"onwall's inequality gives bounds independent of the localization level. Fatou's lemma then proves the first two bounds in \eqref{est:uniform-DA-p}.

    For \eqref{est:uniform-DA-dissipation}, use \eqref{eq:ito-DA} with $p=2$ and raise the resulting inequality to the power $p/2$. The Poisson moment inequality gives
    \begin{equation*}
        \begin{aligned}
            \mathbb E\left(\int_0^{\mathcal T}\int_E\|G_{(n)}\|_V^2\,\eta(\mathrm{d}s,\mathrm{d}\xi)\right)^{p/2} & \leq C_p\mathbb E\left[ \left(\int_0^{\mathcal T}\int_E\|G_{(n)}\|_V^2\,\nu(\mathrm{d}\xi)\,\mathrm{d}s\right)^{p/2} \right. \\
            & \qquad\left.+\int_0^{\mathcal T}\int_E\|G_{(n)}\|_V^p\,\nu(\mathrm{d}\xi)\,\mathrm{d}s\right] \\
            & \leq C_p\left(1+\mathbb E\sup_{t\leq\mathcal T}\|\bm u(t)\|_V^p\right).
        \end{aligned}
    \end{equation*}
    The Burkholder--Davis--Gundy inequality bounds $\mathbb E\sup|\mathcal M_2|^{p/2}$ by the same quantity. The remaining background contribution is bounded by $C_p\mathbb E\|\bm u_0\|_V^{2p}$, by \eqref{eq:int-b}. This proves \eqref{est:uniform-DA-dissipation} when $2p\leq q$.

    Finally, the Poisson moment inequality and the preceding bounds imply
    \begin{equation*}
        \mathbb E|M_n(t)-M_n(s)|_H^p\leq C_p\big(|t-s|^{p/2}+|t-s|\big).
    \end{equation*}
    Integration against $|t-s|^{-1-\alpha p}\,\mathrm{d}s\,\mathrm{d}t$ proves \eqref{est:uniform-DA-moments} for $\alpha p<1$.
\end{proof}

\subsection{Tightness and the Skorokhod--Jakubowski representation theorem}\label{subsec:S-J-theorem}

\subsubsection{Tightness}
To pass to the limit in the Galerkin scheme, we work on the path space
\begin{equation*}
    \mathcal Z_{\mathcal{T}}=D([0,\mathcal{T}];H_{\mathrm w})\cap D([0,\mathcal{T}];V')\cap L^\infty_{\mathrm{w}^\ast}(0,\mathcal{T};V)\cap L^2_{\mathrm w}(0,\mathcal{T};\mathscr D(A))\cap L^2(0,\mathcal{T};H)
\end{equation*}
endowed with the supremum of the five topologies. The first factor consists of weakly c\`adl\`ag $H$-valued paths and carries the weakest topology making $u\mapsto(u(\cdot),h)_H$ continuous into $D([0,\mathcal T];\mathbb R)$ with the Skorokhod $J_1$ topology for every $h\in H$. The remaining factors carry the $J_1$, weak$^\ast$, weak and norm topologies, respectively. The weak$^\ast$ topology is induced by the duality with $L^1(0,\mathcal T;V')$.

\begin{lemma}\label{lem:tight-u_n}
    $\{\mathcal L(\bm u_n)\}_{n\geq1}$ is tight on $(\mathcal{Z}_{\mathcal{T}}, \mathcal{B}(\mathcal{Z}_{\mathcal{T}}))$.
\end{lemma}
\begin{proof}
    Write $F_n=A\bm u_n+\theta(\|\bm u_n-\bm u_n^*\|_V)B_{(n)}(\bm u_n)+R_{(n)}(\bm u_n)$. Write $\bar{\bm u}_n=\bm u_n-\bm u_n^*$. On the support of the cutoff, $\|\bar{\bm u}_n\|_V\leq\kappa$. Expanding $B(\bar{\bm u}_n+\bm u_n^*)$ and using \eqref{est:B-H-symmetric} and Young's inequality gives
    \begin{equation}\label{est:cutoff-drift}
        \theta^2|B_{(n)}(\bm u_n)|_H^2 \leq C_\kappa\big(|A\bar{\bm u}_n|_H^2 +\|\bm u_n^*\|_V^2|A\bm u_n^*|_H^2\big).
    \end{equation}
    Indeed, each mixed term is bounded in square by $C\kappa|A\bar{\bm u}_n|_H\|\bm u_n^*\|_V|A\bm u_n^*|_H$. By \eqref{est:uniform-DA}, \eqref{est:background} and \eqref{eq:int-b},
    \begin{equation}\label{est:galerkin-drift-L2}
        \sup_n\mathbb E\int_0^{\mathcal T}|F_n|_H^2\,\mathrm{d}t \leq C_{\kappa,\mathcal T}\big(1+\mathbb E\|\bm u_0\|_V^4\big)<\infty.
    \end{equation}
    For a stopping time $\tau\leq\mathcal T$, put $\tau_\delta=(\tau+\delta)\wedge\mathcal T$. Cauchy--Schwarz and the Poisson isometry give
    \begin{equation*}
        \mathbb E\left|\int_\tau^{\tau_\delta}F_n\,\mathrm{d}s\right|_H\leq C\delta^{1/2}, \qquad \mathbb E|M_n(\tau_\delta)-M_n(\tau)|_H^2\leq C\delta.
    \end{equation*}
    Consequently the Aldous condition holds:
    \begin{equation}\label{cond:aldous}
        \lim_{\delta\downarrow0}\sup_n\sup_{\tau\leq\mathcal T} \mathbb P\big(|\bm u_n(\tau_\delta)-\bm u_n(\tau)|_H>\varepsilon\big)=0, \qquad\varepsilon>0.
    \end{equation}
    The uniform $V$ bound and the compact embedding $V\hookrightarrow H$ give compact containment in $H$.

    By \eqref{est:galerkin-drift-L2} and \eqref{est:uniform-DA-moments} with $p=2$,
    \begin{equation*}
        \sup_n\mathbb E\|\bm u_n\|_{W^{\alpha,2}(0,\mathcal T;H)}^2<\infty, \qquad 0<\alpha<1/2.
    \end{equation*}
    The compact embedding
    \begin{equation}\label{eq:compact-DA-V}
        L^2(0,\mathcal T;\mathscr D(A))\cap W^{\alpha,2}(0,\mathcal T;H) \hookrightarrow\hookrightarrow L^2(0,\mathcal T;V)
    \end{equation}
    therefore yields tightness in $L^2(0,\mathcal T;V)$. Finally, bounded balls in $L^\infty(0,\mathcal T;V)$ and $L^2(0,\mathcal T;\mathscr D(A))$ are compact in their weak-star and weak topologies, respectively. The uniform estimates and the preceding bounds give tightness on $\mathcal Z_{\mathcal T}$; see \cite[Corollary~4.2]{Motyl2014} for the c\`adl\`ag tightness criterion.
\end{proof}

\subsubsection{Joint Skorokhod representation}
We refer to \cite{0923.60001,8138795} for the following Skorokhod--Jakubowski representation theorem.

\begin{lemma}[Skorokhod--Jakubowski representation theorem]\label{lem:s-j}
    Suppose $X$ is a topological space, and there exists a sequence of continuous functions $g_n:X\to\mathbb{R}$ that separates points of $X$. If $\{\mu_n\}_{n\geq1}$ is a tight sequence of probability measures on $(X,\mathcal{B}(X))$, then there exists a subsequence $\{\mu_{n_k}\}_{k\geq1}$ of $\{\mu_n\}_{n\geq1}$, a probability space $(\hat{\Omega},\hat{\mathcal{F}},\hat{\mathbb{P}})$, $X$\text{-}valued random variables $\nu_k$ and $\nu$ such that the law of $\nu_k$ is $\mu_{n_k}$ for all $k\in\mathbb{N}$ and $\nu_k\to\nu$ $\hat{\mathbb{P}}$\text{-}a.s. as $k\to\infty$.
\end{lemma}

For the noise coordinate, choose a Borel partition $E=\bigcup_{j\geq1}E_j$ with $\nu(E_j)<\infty$, and equip each $E_j$ with a Polish topology generating its Borel $\sigma$-field. Let
\begin{equation}\label{def:joint-noise-space}
    \mathcal N=\prod_{j\geq1}\mathcal N_f([0,\mathcal T]\times E_j),
\end{equation}
where $\mathcal N_f$ is the space of finite point measures with the weak topology. This is a Polish space, and $\eta$ is an $\mathcal N$-valued random variable. Bounded continuous functions on finitely many strips are dense in $L^2([0,\mathcal T]\times E,\mathrm{d}t\otimes\nu)$.

By Lemma \ref{lem:tight-u_n}, the joint laws of $(\bm u_0,\bm u_n,\eta)$ are tight on $V\times\mathcal Z_{\mathcal T}\times\mathcal N$; the initial and noise marginals have fixed laws on Polish spaces. Lemma \ref{lem:s-j} gives a subsequence, relabelled by $n$, and random variables $(a_n,\hat{\bm u}_n,\hat\eta_n)$ and $(\hat{\bm u}_0,\hat{\bm u},\hat\eta)$ on a probability space $(\hat\Omega,\hat{\mathcal F},\hat{\mathbb P})$ such that
\begin{equation}\label{sk-rep-u-law}
    \begin{gathered}
        \mathcal L_{\hat{\mathbb P}}(a_n,\hat{\bm u}_n,\hat\eta_n) =\mathcal L(\bm u_0,\bm u_n,\eta), \\
        (a_n,\hat{\bm u}_n,\hat\eta_n) \to(\hat{\bm u}_0,\hat{\bm u},\hat\eta) \quad\hat{\mathbb P}\text{-a.s. in }V\times\mathcal Z_{\mathcal T}\times\mathcal N.
    \end{gathered}
\end{equation}
In particular, $\mathcal L(\hat{\bm u}_0)=\mathcal L(\bm u_0)$, $\hat{\bm u}_n(0)=P_na_n$, and $\hat{\bm u}(0)=\hat{\bm u}_0$. Set $\hat{\bm u}_n^*(t)=e^{-tA}P_na_n$ and $\hat{\bm u}^*(t)=e^{-tA}\hat{\bm u}_0$. Then $\hat{\bm u}_n^*\to\hat{\bm u}^*$ in $C([0,\mathcal T];V)$ almost surely, by the semigroup contraction and $a_n\to\hat{\bm u}_0$ in $V$.
\begin{corollary}\label{cor:limit-bounds}
    The limit $\hat{\bm u}$ inherits the estimates of Proposition \ref{prop:uniform-DA}: for every $\mathcal{T}>0$,
    \begin{equation}\label{est:u-hat-DA}
        \hat{\mathbb E}\Big[\sup_{t\in[0,\mathcal{T}]} \|\hat{\bm u}(t)-\hat{\bm u}^*(t)\|^2_V +\int_0^{\mathcal{T}}\|\hat{\bm u}-\hat{\bm u}^*\|^2_{\mathscr D(A)} \,\mathrm{d}t\Big]\leq K,
    \end{equation}
    where $K$ is the constant in \eqref{est:uniform-DA}. Together with \eqref{est:background}, this gives \eqref{est:u-hat-strong}.
\end{corollary}
\begin{proof}
    By \eqref{sk-rep-u-law}, estimate \eqref{est:uniform-DA} holds for $\hat{\bm u}_n-\hat{\bm u}_n^*$ with a constant independent of $n$. Combining the background convergence established above with $\hat{\bm u}_n\to\hat{\bm u}$ almost surely in $\mathcal Z_{\mathcal T}$, we find that $(\hat{\bm u}_n-\hat{\bm u}_n^*)$ converges almost surely to $\hat{\bm u}-\hat{\bm u}^*$, is bounded in $L^\infty(0,\mathcal{T};V)$ and in $L^2(0,\mathcal{T};\mathscr D(A))$ almost surely, and converges weakly in $L^2(0,\mathcal T;\mathscr D(A))$ and weakly$^\ast$ in $L^\infty(0,\mathcal T;V)$. The functional $x\mapsto\sup_{t\leq\mathcal T}\|x(t)\|_V^2$, with $\|z\|_V=\infty$ for $z\in V'\setminus V$, is lower semicontinuous under $J_1$ convergence in $V'$. The $L^2(0,\mathcal T;\mathscr D(A))$ norm is lower semicontinuous under weak convergence. Fatou's lemma on $\hat\Omega$ therefore gives \eqref{est:u-hat-DA}. The final assertion follows from \eqref{est:background} and the triangle inequality.
\end{proof}

\begin{lemma}\label{lem:uni-est}
    For $\mathcal{T}>0$ and $p=1,2$, we have
    \begin{align}
        \sup_{n\geq1}\hat{\mathbb{E}}\left(\sup_{t\in[0,\mathcal{T}]}|\hat{\bm u}_n(t)|_H^{2p}+\int_0^{\mathcal{T}}|\hat{\bm u}_n(t)|_H^{2(p-1)}\|\hat{\bm u}_n(t)\|_V^2\,\mathrm{d}t\right) & \leq C(\mathcal{T})\label{est:u-hat-n}, \\
        \hat{\mathbb{E}}\left(\sup_{t\in[0,\mathcal{T}]}|\hat{\bm u}(t)|_H^{2p}+\int_0^{\mathcal{T}}|\hat{\bm u}(t)|_H^{2(p-1)}\|\hat{\bm u}(t)\|_V^2\,\mathrm{d}t\right) & \leq C(\mathcal{T})\label{est:u-hat}.
    \end{align}
\end{lemma}
\begin{proof}
    The estimate \eqref{est:u-hat-n} is directly provided by \eqref{est:galerkin-uniform} and \eqref{sk-rep-u-law}. We prove \eqref{est:u-hat} in two steps.

    \textbf{Step 1: the supremum term.} Since $\hat{\bm u}_n\to \hat{\bm u}$, $\hat{\mathbb{P}}$-a.s.\ in $D([0,\mathcal{T}];H_{\mathrm{w}})$, and convergence in the Skorokhod topology implies pointwise convergence at every continuity point of the limit, we have, for $\hat{\mathbb{P}}$-a.e.\ $\omega$ and every continuity point $t$ of $\hat{\bm u}$,
    \begin{equation*}
        \hat{\bm u}_n(t)\rightharpoonup \hat{\bm u}(t)\quad\text{weakly in }H.
    \end{equation*}
    Since $\hat{\bm u}\in D([0,\mathcal{T}];H_{\mathrm{w}})$, the set of its continuity points is dense in $[0,\mathcal{T}]$. For such a point $t$, the weak lower semicontinuity of $|\cdot|_H$ yields
    \begin{equation*}
        |\hat{\bm u}(t)|_H\leq\liminf_{n\to\infty}|\hat{\bm u}_n(t)|_H\leq\liminf_{n\to\infty}\sup_{s\in[0,\mathcal{T}]}|\hat{\bm u}_n(s)|_H.
    \end{equation*}
    For an arbitrary $t\in[0,\mathcal{T})$, choose continuity points $t_m\downarrow t$; by the weak right continuity of $\hat{\bm u}$, $\hat{\bm u}(t_m)\rightharpoonup \hat{\bm u}(t)$ weakly in $H$, and hence
    \begin{equation*}
        |\hat{\bm u}(t)|_H\leq\liminf_{m\to\infty}|\hat{\bm u}(t_m)|_H\leq\liminf_{n\to\infty}\sup_{s\in[0,\mathcal{T}]}|\hat{\bm u}_n(s)|_H.
    \end{equation*}
    At the endpoint $t=\mathcal{T}$, since the time reparametrizations in the Skorokhod topology fix the endpoints, the convergence $\hat{\bm u}_n\to\hat{\bm u}$ in $D([0,\mathcal{T}];H_{\mathrm{w}})$ implies $\hat{\bm u}_n(\mathcal{T})\rightharpoonup\hat{\bm u}(\mathcal{T})$ weakly in $H$, and the same estimate follows from the weak lower semicontinuity of $|\cdot|_H$.

    Taking the supremum over $t\in[0,\mathcal{T}]$ and then the $2p$-th power, we obtain, $\hat{\mathbb{P}}$-a.s.,
    \begin{equation*}
        \sup_{t\in[0,\mathcal{T}]}|\hat{\bm u}(t)|_H^{2p}\leq\liminf_{n\to\infty}\sup_{t\in[0,\mathcal{T}]}|\hat{\bm u}_n(t)|_H^{2p}.
    \end{equation*}
    Therefore, by Fatou's lemma and \eqref{est:u-hat-n},
    \begin{equation*}
        \begin{aligned}
            \hat{\mathbb{E}}\sup_{t\in[0,\mathcal{T}]}|\hat{\bm u}(t)|_H^{2p} & \leq\hat{\mathbb{E}}\liminf_{n\to\infty}\sup_{t\in[0,\mathcal{T}]}|\hat{\bm u}_n(t)|_H^{2p} \\
            & \leq\liminf_{n\to\infty}\hat{\mathbb{E}}\sup_{t\in[0,\mathcal{T}]}|\hat{\bm u}_n(t)|_H^{2p} \leq C(\mathcal{T}).
        \end{aligned}
    \end{equation*}

    \textbf{Step 2: the integral term.} Set $d_n(t,\omega):=|\hat{\bm u}_n(t,\omega)-\hat{\bm u}(t,\omega)|_H$. Since $\hat{\bm u}_n\to\hat{\bm u}$, $\hat{\mathbb{P}}$-a.s.\ in $L^2(0,\mathcal{T};H)$, Chebyshev's inequality yields, for $\hat{\mathbb{P}}$-a.e.\ $\omega$ and every $\varepsilon>0$,
    \begin{equation*}
        m_n(\omega):=\big|\{t\in[0,\mathcal{T}]:d_n(t,\omega)>\varepsilon\}\big| \leq\varepsilon^{-2}\int_0^{\mathcal{T}}d_n(t,\omega)^2\,\mathrm{d}t\to 0,
    \end{equation*}
    i.e., $d_n(\cdot,\omega)\to0$ in Lebesgue measure. Since $0\leq m_n\leq\mathcal{T}$ and $m_n\to0$ $\hat{\mathbb{P}}$-a.s., the bounded convergence theorem gives
    \begin{equation*}
        (\mathrm{d}t\otimes\hat{\mathbb{P}})\big(\{d_n>\varepsilon\}\big)=\hat{\mathbb{E}}\,m_n\to 0,
    \end{equation*}
    i.e., $d_n\to0$ in measure on $[0,\mathcal{T}]\times\hat{\Omega}$. Hence there exists a subsequence $\{\hat{\bm u}_{n_k}\}_{k\geq1}$, which is independent of $\omega$, such that for $\mathrm{d}t\otimes\hat{\mathbb{P}}$-a.e. $(t,\omega)$, we have
    \begin{equation*}
        |\hat{\bm u}_{n_k}(t,\omega)-\hat{\bm u}(t,\omega)|_H\to0,\quad\text{as }k\to\infty.
    \end{equation*}
    Let $A\subset[0,\mathcal{T}]\times\hat{\Omega}$ be the set of full $\mathrm{d}t\otimes\hat{\mathbb{P}}$-measure on which the above convergence holds, and fix $(t,\omega)\in A$. With the convention $\|\bm u\|_V:=+\infty$ for $\bm u\in H\setminus V$, we claim that
    \begin{equation*}
        \|\hat{\bm u}(t,\omega)\|_V^2\leq\liminf_{k\to\infty}\|\hat{\bm u}_{n_k}(t,\omega)\|_V^2.
    \end{equation*}
    Indeed, if the right-hand side is infinite, there is nothing to prove. Otherwise, choose a further subsequence, still denoted by $\{n_k\}$, realizing the $\liminf$; then $\{\hat{\bm u}_{n_k}(t,\omega)\}$ is bounded in $V$, so it admits a subsequence converging weakly in $V$. Since $V\hookrightarrow H$ and $\hat{\bm u}_{n_k}(t,\omega)\to\hat{\bm u}(t,\omega)$ in $H$, the weak limit in $V$ must coincide with $\hat{\bm u}(t,\omega)$, and the claim follows from the weak lower semicontinuity of $\|\cdot\|_V$.

    Since also $|\hat{\bm u}_{n_k}(t,\omega)|_H^{2(p-1)}\to|\hat{\bm u}(t,\omega)|_H^{2(p-1)}$, and since for non-negative sequences $\{\alpha_k\}_{k\geq1}$ and $\{\beta_k\}_{k\geq1}$ with $\alpha_k\to\alpha$ and $\beta\leq\liminf_{k\to\infty}\beta_k$ one has
    \begin{equation*}
        \alpha\beta\leq\liminf_{k\to\infty}\alpha_k\beta_k,
    \end{equation*}
    we infer that, for $\mathrm{d}t\otimes\hat{\mathbb{P}}$-a.e.\ $(t,\omega)$,
    \begin{equation*}
        |\hat{\bm u}(t,\omega)|_H^{2(p-1)}\|\hat{\bm u}(t,\omega)\|_V^2 \leq\liminf_{k\to\infty}(|\hat{\bm u}_{n_k}(t,\omega)|_H^{2(p-1)}\|\hat{\bm u}_{n_k}(t,\omega)\|_V^2).
    \end{equation*}
    Applying Fatou's lemma with respect to the measure $\mathrm{d}t\otimes\hat{\mathbb{P}}$ and using \eqref{est:u-hat-n}, we obtain
    \begin{equation*}
        \begin{aligned}
            \hat{\mathbb{E}}\int_0^{\mathcal{T}}|\hat{\bm u}(t)|_H^{2(p-1)}\|\hat{\bm u}(t)\|_V^2\,\mathrm{d}t & \leq\hat{\mathbb{E}}\int_0^{\mathcal{T}}\liminf_{k\to\infty}(|\hat{\bm u}_{n_k}(t)|_H^{2(p-1)}\|\hat{\bm u}_{n_k}(t)\|_V^2)\,\mathrm{d}t \\
            & \leq\liminf_{k\to\infty}\hat{\mathbb{E}}\int_0^{\mathcal{T}}|\hat{\bm u}_{n_k}(t)|_H^{2(p-1)}\|\hat{\bm u}_{n_k}(t)\|_V^2\,\mathrm{d}t \\
            & \leq\sup_{n\geq1}\hat{\mathbb{E}}\int_0^{\mathcal{T}}|\hat{\bm u}_n(t)|_H^{2(p-1)}\|\hat{\bm u}_n(t)\|_V^2\,\mathrm{d}t \\
            & \leq C(\mathcal{T}).
        \end{aligned}
    \end{equation*}

    Combining the estimates from Step 1 and Step 2, and relabelling the constant $C(\mathcal{T})$ if necessary, we obtain \eqref{est:u-hat}.
\end{proof}

\subsection{Identification of the limiting equation}

\subsubsection{Square integrability of the limit martingale}
For the truncated convection terms, set
\begin{align*}
    \mathcal B_n(s) & =\theta(\|\bm u_n(s)-\bm u_n^*(s)\|_V)B_{(n)}(\bm u_n(s)), \\
    \widehat{\mathcal B}_n(s) & =\theta(\|\hat{\bm u}_n(s)-\hat{\bm u}_n^*(s)\|_V)B_{(n)}(\hat{\bm u}_n(s)), \\
    \widehat{\mathcal B}(s) & =\theta(\|\hat{\bm u}(s)-\hat{\bm u}^*(s)\|_V)B(\hat{\bm u}(s)).
\end{align*}

For $n\in\mathbb N$, $\psi\in V$ and $t\in[0,\mathcal T]$, define the processes
\begin{equation}\label{def:hat-M-n}
    M_n(t):= \bm u_n(t)-\bm u_n(0)+\int_0^tA\bm u_n(s)+\mathcal B_n(s)+R_{(n)}(\bm u_n(s))\,\mathrm{d}s,
\end{equation}
and
\begin{equation}\label{def:hat-N-n-psi}
    \begin{aligned}
        N_{n,\psi}(t):= & {}_{V^\prime}\langle\bm u_n(t),\psi\rangle_V^2-{}_{V^\prime}\langle\bm u_n(0),\psi\rangle_V^2-\int_0^t\int_E{}_{V^\prime}\langle G_{(n)}(\bm u_n(s-),\xi),\psi\rangle_V^2\,\nu(\mathrm{d}\xi)\,\mathrm{d}s \\
        & +2\int_0^t{}_{V^\prime}\langle A\bm u_n(s)+\mathcal B_n(s)+R_{(n)}(\bm u_n(s)),\psi\rangle_V\,{}_{V^\prime}\langle\bm u_n(s),\psi\rangle_V\,\mathrm{d}s.
    \end{aligned}
\end{equation}
Similarly, $\hat M_n(t)$ and $\hat{N}_{n,\psi}(t)$ are defined by replacing $\bm u_n$ with $\hat{\bm u}_n$ in \eqref{def:hat-M-n} and \eqref{def:hat-N-n-psi}, respectively.

For the filtration definitions, paths are extended constantly beyond $\mathcal T$. For every $t\geq 0$, let
\begin{equation*}
    \hat{\mathcal H}_{n,t}:=\sigma\{\hat{\bm u}_n(s):0\leq s\leq t\}.
\end{equation*}
We then define the right-continuous modification by
\begin{equation*}
    \hat{\mathcal H}_{n,t}^+:=\bigcap_{s>t}\hat{\mathcal H}_{n,s},
\end{equation*}
and denote by $\hat{\mathcal N}$ the collection of all $\hat{\mathbb P}$-null sets in $\hat{\mathcal F}$. Finally, we define the filtration
\begin{equation*}
    \hat{\mathcal F}_{n,t}:=\sigma(\hat{\mathcal H}_{n,t}^+\cup\hat{\mathcal N}).
\end{equation*}
The filtration $\hat{\mathbb F}_n:=\{\hat{\mathcal F}_{n,t}\}_{t\geq0}$ is the usual augmentation of $\{\hat{\mathcal H}_{n,t}\}_{t\geq0}$ and satisfies
\begin{equation*}
    \hat{\mathcal F}_{n,t}=\bigcap_{s>t}\sigma(\hat{\mathcal H}_{n,s}\cup \hat{\mathcal N}).
\end{equation*}
See \cite[Proposition~5.3 and Lemma~E.2]{8220046}.

\begin{lemma}\label{lem:martingale-1}
    For every $n\in\mathbb{N}$, the process $\hat M_n$ is a square-integrable $V^\prime$-valued $\hat{\mathbb F}_n$-martingale. Moreover, for $p=1,2$, we have
    \begin{equation}\label{est:uni-hat-M-n}
        \sup_{n\geq1}\hat{\mathbb{E}}\sup_{t\in[0,\mathcal{T}]}\|\hat M_n(t)\|_{V^\prime}^{2p}<\infty.
    \end{equation}
\end{lemma}

\begin{proof}
    On $D([0,\mathcal T];V')$, let $\pi_r(x)=x(r)$ and define
    \begin{equation}\label{eq:canonical-past-filtration}
        \begin{gathered}
            \mathscr D_t^0=\sigma(\pi_r:0\leq r\leq t), \\
            \mathscr D_t=\bigcap_{t<r\leq\mathcal T}\mathscr D_r^0 \quad(t<\mathcal T),\qquad \mathscr D_{\mathcal T}=\mathscr D_{\mathcal T}^0.
        \end{gathered}
    \end{equation}
    Fix $0\leq s<t\leq\mathcal T$, $\psi\in V$, and a bounded $\mathscr D_s^0$-measurable function $h$. Then $h(\bm u_n)$ is $\mathcal F_s$-measurable. Since $M_n$ is a martingale,
    \begin{equation*}
        \mathbb E[{}_{V'}\langle M_n(t)-M_n(s),\psi\rangle_V h(\bm u_n)]=0.
    \end{equation*}
    The drift integrals are measurable functions of the path up to their upper integration limit; the background evolution depends only on its initial value. Thus $\hat M_n$ is adapted to $\hat{\mathcal H}_{n,t}$, and equality of path laws gives
    \begin{equation}\label{eq:martingale-past-transfer}
        \hat{\mathbb E}[{}_{V'}\langle\hat M_n(t)-\hat M_n(s),\psi\rangle_V h(\hat{\bm u}_n)]=0.
    \end{equation}
    Integrability also transfers from $M_n$, so this proves the martingale property for the raw natural filtration.

    We record the augmentation argument used below. If $X$ is a c\`adl\`ag martingale for a raw filtration $\{\mathcal H_t\}_{t\geq0}$ and $\mathbb E\sup_{t\leq\mathcal T}|X(t)|<\infty$, it remains a martingale for its usual augmentation. Indeed, for $s<t$ choose $s_j\downarrow s$ with $s_j<t$. For every event $F$ in the augmented $\sigma$-field at $s$, $F$ belongs to the completion of $\mathcal H_{s_j}$, and hence $\mathbb E[\mathbbm{1}_F(X(t)-X(s_j))]=0$. Dominated convergence gives the same identity at $s$. Applied to the scalar pairings of $\hat M_n$, this proves the $\hat{\mathbb F}_n$-martingale property. The required integrable supremum follows from equality in law with the square-integrable original martingale and Doob's inequality.

    To prove \eqref{est:uni-hat-M-n}, set $F_n(s,\xi)=P_nG(\bm u_n(s-),\xi)$. The Galerkin equation gives the $H$-valued identity
    \begin{equation}\label{eq:vector-martingale-integral}
        M_n(t)=\int_0^t\int_E F_n(s,\xi)\, \tilde\eta(\mathrm{d}s,\mathrm{d}\xi).
    \end{equation}
    For $p=1,2$, the Hilbert-space martingale maximal inequality yields
    \begin{equation}\label{est:vector-martingale-maximal}
        \mathbb E\sup_{t\leq\mathcal T}|M_n(t)|_H^{2p} \leq C_p\mathbb E\left(\int_0^{\mathcal T}\int_E |F_n|_H^2\,\nu(\mathrm{d}\xi)\,\mathrm{d}s\right)^p +C_p\mathbb E\int_0^{\mathcal T}\int_E |F_n|_H^{2p}\,\nu(\mathrm{d}\xi)\,\mathrm{d}s.
    \end{equation}
    The constants are independent of $n$. For $p=1$, this follows from Doob's inequality and the Poisson isometry. For $p=2$, the Hilbert-space Burkholder--Davis--Gundy inequality bounds the left-hand side by $C\mathbb E(\int_0^{\mathcal T}\int_E|F_n|_H^2\,\mathrm{d}\eta)^2$. Decompose this integral into its compensated integral and compensator, and use the scalar Poisson isometry to obtain \eqref{est:vector-martingale-maximal}.

    Since $P_n$ is an $H$-contraction, \eqref{asm:g_linear-H} gives
    \begin{equation}\label{est:projected-noise-moments}
        \int_E|F_n(s,\xi)|_H^{2p}\,\nu(\mathrm{d}\xi) \leq C_p(1+|\bm u_n(s-)|_H^{2p}),\qquad p=1,2.
    \end{equation}
    H\"older's inequality in time and \eqref{est:galerkin-uniform} therefore imply
    \begin{equation*}
        \sup_n\mathbb E\sup_{t\leq\mathcal T}|M_n(t)|_H^{2p} \leq C_{p,\mathcal T}\left(1+ \sup_n\mathbb E\sup_{t\leq\mathcal T}|\bm u_n(t)|_H^{2p}\right) <\infty.
    \end{equation*}
    Finally, let $c_{HV'}$ be the norm of the embedding $H\hookrightarrow V'$. Equality of path laws and this embedding yield
    \begin{equation}\label{est:vector-martingale-transfer}
        \begin{aligned}
            \sup_n\hat{\mathbb E}\sup_{t\leq\mathcal T}|\hat M_n(t)|_{V'}^{2p} & =\sup_n\mathbb E\sup_{t\leq\mathcal T}|M_n(t)|_{V'}^{2p} \\
            & \leq c_{HV'}^{2p}\sup_n\mathbb E\sup_{t\leq\mathcal T}|M_n(t)|_H^{2p} <\infty.
        \end{aligned}
    \end{equation}
\end{proof}

\begin{lemma}\label{lem:square-martingale-galerkin}
    For every $n\in\mathbb N$ and $\psi\in V$, the process $\hat N_{n,\psi}$ is a square-integrable real-valued $\hat{\mathbb F}_n$-martingale. Moreover,
    \begin{equation}\label{est:N_n,psi}
        \sup_n\hat{\mathbb E}\sup_{t\leq\mathcal T}|\hat N_{n,\psi}(t)|^2 =\sup_n\mathbb E\sup_{t\leq\mathcal T}|N_{n,\psi}(t)|^2 \leq C_{\mathcal T}\|\psi\|_V^4.
    \end{equation}
\end{lemma}
\begin{proof}
    Set
    \begin{equation*}
        Z_n(t)={}_{V'}\langle\bm u_n(t),\psi\rangle_V,\qquad \gamma_n(s,\xi)={}_{V'}\langle G_{(n)}(\bm u_n(s-),\xi),\psi\rangle_V,
    \end{equation*}
    and $d_n(s)={}_{V'}\langle A\bm u_n(s)+\mathcal B_n(s) +R_{(n)}(\bm u_n(s)),\psi\rangle_V$. The scalar equation is $\mathrm{d}Z_n=-d_n\,\mathrm{d}t+\int_E\gamma_n\,\tilde\eta(\mathrm{d}t,\mathrm{d}\xi)$. The jump It\^o formula gives
    \begin{equation}\label{eq:square-ito-galerkin}
        \begin{aligned}
            Z_n(t)^2={} & Z_n(0)^2-2\int_0^tZ_n(s)d_n(s)\,\mathrm{d}s \\
            & +2\int_0^t\int_E Z_n(s-)\gamma_n(s,\xi)\, \tilde\eta(\mathrm{d}s,\mathrm{d}\xi) +\int_0^t\int_E\gamma_n(s,\xi)^2\,\eta(\mathrm{d}s,\mathrm{d}\xi).
        \end{aligned}
    \end{equation}
    Consequently, compensating the last integral and using \eqref{def:hat-N-n-psi}, we obtain
    \begin{equation}\label{eq:square-martingale-integrand}
        N_{n,\psi}(t)=\int_0^t\int_E \big[2Z_n(s-)\gamma_n(s,\xi)+\gamma_n(s,\xi)^2\big] \,\tilde\eta(\mathrm{d}s,\mathrm{d}\xi).
    \end{equation}
    The integrand is predictable. Since $P_n$ is an $H$-contraction, \eqref{asm:g_linear-H} and \eqref{est:galerkin-uniform} with $p=2$ imply
    \begin{equation}\label{est:square-martingale-integrand}
        \begin{aligned}
            & \mathbb E\int_0^{\mathcal T}\int_E |2Z_n(s-)\gamma_n+\gamma_n^2|^2\,\nu(\mathrm{d}\xi)\,\mathrm{d}s \\
            & \quad\leq C\|\psi\|_V^4\mathbb E\int_0^{\mathcal T} \left[|\bm u_n(s-)|_H^2\int_E|G_{(n)}(\bm u_n(s-),\xi)|_H^2\,\nu(\mathrm{d}\xi) +\int_E|G_{(n)}(\bm u_n(s-),\xi)|_H^4\,\nu(\mathrm{d}\xi)\right]\mathrm{d}s \\
            & \quad\leq C\|\psi\|_V^4\mathbb E\int_0^{\mathcal T} (1+|\bm u_n(s)|_H^4)\,\mathrm{d}s \leq C_{\mathcal T}\|\psi\|_V^4.
        \end{aligned}
    \end{equation}
    Thus \eqref{eq:square-martingale-integrand} is a square-integrable martingale. Doob's $L^2$ inequality and the compensated Poisson isometry bound its expected squared supremum by four times the left-hand side of \eqref{est:square-martingale-integrand}. Equality of path laws transfers this bound to $\hat N_{n,\psi}$; the past-measurable testing and augmentation argument of Lemma \ref{lem:martingale-1} transfers the martingale property. This proves \eqref{est:N_n,psi}.
\end{proof}

Let us introduce the following limiting filtration. For every $t\geq 0$, we define
\begin{equation*}
    \mathcal{H}_t:=\sigma\{\hat{\bm u}(s):0\leq s\leq t\},
\end{equation*}
and the right-continuous modification by
\begin{equation*}
    \mathcal{H}_t^+:=\bigcap_{s>t}\mathcal{H}_s.
\end{equation*}
Recall that $\hat{\mathcal N}$ is the collection of all $\hat{\mathbb P}$-null sets in $\hat{\mathcal F}$. We then define the filtration
\begin{equation*}
    \hat{\mathcal F}_t:=\sigma(\mathcal{H}_t^+\cup\hat{\mathcal N}).
\end{equation*}

For $t\in [0,\mathcal{T}]$, let us define a $V^\prime$-valued process $\hat M$ by
\begin{equation}\label{eqn:def-hat-M}
    \hat M(t):=\hat{\bm u}(t)-\hat{\bm u}(0)+\int_0^tA\hat{\bm u}(s)+\widehat{\mathcal B}(s)+R(\hat{\bm u}(s))\,\mathrm{d}s.
\end{equation}
Also, for every $\psi\in V$, we define a real-valued process $\hat N_\psi$ by
\begin{equation}\label{eqn:def-hat-N-psi}
    \begin{aligned}
        \hat N_\psi(t):= & {}_{V^\prime}\langle\hat{\bm u}(t),\psi\rangle_V^2-{}_{V^\prime}\langle\hat{\bm u}(0),\psi\rangle_V^2-\int_0^t\int_E{}_{V^\prime}\langle G(\hat{\bm u}(s-),\xi),\psi\rangle_V^2\,\nu(\mathrm{d}\xi)\,\mathrm{d}s \\
        & +2\int_0^t{}_{V^\prime}\langle A\hat{\bm u}(s)+\widehat{\mathcal B}(s)+R(\hat{\bm u}(s)),\psi\rangle_V\,{}_{V^\prime}\langle\hat{\bm u}(s),\psi\rangle_V\,\mathrm{d}s.
    \end{aligned}
\end{equation}

\begin{lemma}\label{lem:martingale-3}
    The process $\hat M$ is a square-integrable $V^\prime$-valued $\hat{\mathbb F}$-martingale. Moreover, for $p=1,2$, we have
    \begin{equation}\label{est:uni-hat-M}
        \hat{\mathbb{E}}\sup_{t\in[0,\mathcal{T}]}\|\hat M(t)\|_{V^\prime}^{2p}<\infty.
    \end{equation}
\end{lemma}
\begin{proof}
    \textbf{Step 1: square integrability.} First, we claim that, for $\hat{\mathbb{P}}$-a.e.\ $\omega$, $\hat M_n(t)\rightharpoonup\hat M(t)$ in $V^\prime$ at every continuity point $t$ of $\hat{\bm u}$. Since $\hat{\bm u}_n\to\hat{\bm u}$ in $D([0,\mathcal{T}];V^\prime)$, we have, for $\hat{\mathbb{P}}$-a.e.\ $\omega$ and every continuity point $t$ of $\hat{\bm u}$,
    \begin{equation*}
        \hat{\bm u}_n(t)\to\hat{\bm u}(t)\quad\text{in }V^\prime.
    \end{equation*}
    For $\psi\in V$, continuity of the form $a$ and \eqref{eq:strong-L2-V} give
    \begin{equation*}
        \sup_{t\leq\mathcal T}\left| \int_0^t\langle A(\hat{\bm u}_n-\hat{\bm u}),\psi\rangle\,\mathrm{d}s\right| =\sup_{t\leq\mathcal T}\left| \int_0^t a(\hat{\bm u}_n-\hat{\bm u},\psi)\,\mathrm{d}s\right| \leq C\mathcal T^{1/2}\|\psi\|_V \|\hat{\bm u}_n-\hat{\bm u}\|_{L^2(0,\mathcal T;V)}\to0.
    \end{equation*}
    The background processes converge in $CV$, and \eqref{eq:strong-L2-V} implies convergence of the cutoffs in $L^2(0,\mathcal T)$. Together with the drift convergence \eqref{eq:convergence-B}, this yields
    \begin{equation}
        \int_0^t{}_{V^\prime}\Big\langle \theta\big(\|\hat{\bm u}_n(s)-\hat{\bm u}_n^*(s)\|_V\big) B_{(n)}(\hat{\bm u}_n(s)) -\theta\big(\|\hat{\bm u}(s)-\hat{\bm u}^*(s)\|_V\big) B(\hat{\bm u}(s)),\psi\Big\rangle_V\,\mathrm{d}s\to 0.
    \end{equation}
    This identifies the convection term in the drift residual. The stochastic integral is identified in Lemma \ref{lem:retained-noise}. After obtaining the $V$-valued c\`adl\`ag version, we use the exit time
    \begin{equation}\label{def:tau-exit}
        \tau=\inf\{t\in[0,\mathcal T]: \|\hat{\bm u}(t)-\hat{\bm u}^*(t)\|_V\geq\kappa/2\} \wedge\mathcal T,\qquad\inf\varnothing=\infty.
    \end{equation}
    To prove the drift convergence, decompose
    \begin{equation*}
        \begin{aligned}
            & \lim_{n\to\infty}\int_0^t{}_{V^\prime}\langle B_{(n)}(\hat{\bm u}_n(s))-B(\hat{\bm u}(s)),\psi\rangle_V\,\mathrm{d}s \\
            = & \lim_{n\to\infty}\int_0^t{}_{V^\prime}\langle P_n B(\hat{\bm u}_n(s))-B(\hat{\bm u}(s)),\psi\rangle_V\,\mathrm{d}s \\
            = & \lim_{n\to\infty}\int_0^t{}_{V^\prime}\langle B(\hat{\bm u}_n),(P_n-I)\psi\rangle_V\,\mathrm{d}s \\
            & +\lim_{n\to\infty}\int_0^tb(\hat{\bm u}_n-\hat{\bm u},\hat{\bm u}_n,\psi)\,\mathrm{d}s+\lim_{n\to\infty}\int_0^tb(\hat{\bm u},\hat{\bm u}_n-\hat{\bm u},\psi)\,\mathrm{d}s \\
            := & \lim (I_{1,n}+I_{2,n}+I_{3,n}).
        \end{aligned}
    \end{equation*}
    We first prove the convergence for $\psi\in V\cap\mathscr D(A)$; the general case follows by density, because on the full measure set where the almost sure bounds displayed below hold all the terms are continuous in $\psi\in V$. On the full probability set of convergence in $\mathcal Z_{\mathcal T}$, weak convergence in $L^2(0,\mathcal T;\mathscr D(A))$ and weak-star convergence in $L^\infty(0,\mathcal T;V)$ imply, by uniform boundedness,
    \begin{equation}\label{eq:unif-a.s.}
        K(\omega):=\sup_n\|\hat{\bm u}_n\|_{L^\infty(0,\mathcal T;V)} +\sup_n\|\hat{\bm u}_n\|_{L^2(0,\mathcal T;\mathscr D(A))}<\infty.
    \end{equation}
    Strong convergence in $L^2(0,\mathcal T;H)$ and the spectral interpolation inequality give
    \begin{equation}\label{eq:strong-L2-V}
        \|\hat{\bm u}_n-\hat{\bm u}\|_{L^2V}^2 \leq\|\hat{\bm u}_n-\hat{\bm u}\|_{L^2H} \|\hat{\bm u}_n-\hat{\bm u}\|_{L^2\mathscr D(A)}\to0.
    \end{equation}
    On this set, by \eqref{est:B_op-dual-norm},
    \begin{equation*}
        \int_0^t\|B(\hat{\bm u}_n)\|_{V^\prime}\,\mathrm{d}s \leq C\int_0^t\Big(|\hat{\bm u}_n|_H^{1/2}\|\hat{\bm u}_n\|_V^{3/2} +\|\hat{\bm u}_n\|_V^{3/2}\|\hat{\bm u}_n\|_{\mathscr D(A)}^{1/2}\Big)\mathrm{d}s \leq C(\omega),
    \end{equation*}
    uniformly in $n$, and the fact that $\hat{\bm u}_n\rightharpoonup\hat{\bm u}$ in $L^2(0,\mathcal{T};V)$ shows that $\{\hat{\bm u}_n\}_{n\geq1}$ is bounded in $L^2(0,\mathcal{T};V)$; hence
    \begin{equation*}
        |I_{1,n}|\leq C(\omega)\|(P_n-I)\psi\|_V\to 0,\quad\text{as}\ n\to\infty.
    \end{equation*}
    Set $\delta_n=\hat{\bm u}_n-\hat{\bm u}$. For $\psi\in\mathscr D(A)$, skew symmetry and \eqref{est:B_op-prod-norm} imply
    \begin{equation*}
        |b(a,b,\psi)|=|b(a,\psi,b)|\leq C_\psi\|a\|_V\|b\|_V.
    \end{equation*}
    Hence
    \begin{equation*}
        \int_0^{\mathcal T}|b(\delta_n,\hat{\bm u}_n,\psi)|+|b(\hat{\bm u},\delta_n,\psi)|\,\mathrm{d}s \leq C_\psi(\|\hat{\bm u}_n\|_{L^2V}+\|\hat{\bm u}\|_{L^2V})\|\delta_n\|_{L^2V}\to0.
    \end{equation*}
    Together with $I_{1,n}\to0$ and density in $V$, this proves, uniformly in $t\leq\mathcal T$,
    \begin{equation}\label{eq:convergence-B}
        \int_0^t\langle B_{(n)}(\hat{\bm u}_n)-B(\hat{\bm u}),\psi\rangle\,\mathrm{d}s\to0, \qquad\psi\in V.
    \end{equation}
    The bound $\|B(\hat{\bm u}_n)\|_{L^2V'}\leq C(\omega)$ follows from \eqref{est:B_op-dual-norm} and \eqref{eq:unif-a.s.}. The cutoffs converge strongly in $L^2$ by \eqref{eq:strong-L2-V} and convergence of the backgrounds. For smooth $\psi$, the preceding difference estimate also holds after multiplication by a bounded cutoff; the cutoff difference paired with $B(\hat{\bm u})$ tends to zero by Cauchy--Schwarz. Density and the uniform $L^1V'$ bound give the claimed cutoff drift convergence for every $\psi\in V$. Similarly, since \eqref{est:r} shows that $R$ is a bounded linear operator, we have
    \begin{equation*}
        \begin{aligned}
            & \lim_{n\to\infty}\int_0^t{}_{V^\prime}\langle R_{(n)}(\hat{\bm u}_n(s))-R(\hat{\bm u}(s)),\psi\rangle_V\,\mathrm{d}s \\
            = & \lim_{n\to\infty}\int_0^t{}_{V^\prime}\langle R(\hat{\bm u}_n),(P_n-I)\psi\rangle_V\,\mathrm{d}s \\
            & +\lim_{n\to\infty}\int_0^t{}_{V^\prime}\langle R(\hat{\bm u}_n(s)-\hat{\bm u}(s)),\psi\rangle_V\,\mathrm{d}s \\
            = & 0.
        \end{aligned}
    \end{equation*}
    Therefore, we conclude that for $\hat{\mathbb{P}}$-a.e.\,$\omega$ and every continuity point $t$ of $\hat{\bm u}$,
    \begin{equation}\label{conver:M_n-hat}
        \hat M_n(t)\rightharpoonup \hat M(t)\quad\text{in }V^\prime.
    \end{equation}

    The rest of the proof in this step is similar to the proof of Lemma \ref{lem:uni-est}. \eqref{conver:M_n-hat} shows that for $\hat{\mathbb{P}}$-a.e.\,$\omega$ and every continuity point $t$ of $\hat{\bm u}$,
    \begin{equation*}
        \|\hat M(t)\|_{V^\prime}\leq \liminf_{n\to\infty}\|\hat M_n(t)\|_{V^\prime}\leq\liminf_{n\to\infty}\sup_{t\in[0,\mathcal{T}]}\|\hat M_n(t)\|_{V^\prime}.
    \end{equation*}
    For an arbitrary $t\in[0,\mathcal{T})$, we can choose continuity points $t_m\downarrow t$ as $m\to\infty$. Right continuity of $\hat M$ gives
    \begin{equation*}
        \|\hat M(t)\|_{V^\prime}=\lim_{m\to\infty}\|\hat M(t_m)\|_{V^\prime}\leq \liminf_{n\to\infty}\sup_{t\in[0,\mathcal{T}]}\|\hat M_n(t)\|_{V^\prime}.
    \end{equation*}
    At the endpoint $t=\mathcal{T}$, since the time reparametrizations in the Skorokhod topology fix the endpoints, we directly have the same estimate.

    Therefore, for $p=1,2$, by Fatou's lemma and \eqref{est:uni-hat-M-n}, we have
    \begin{equation}\label{est:martingale-lem-3-step1-1}
        \hat{\mathbb{E}}\sup_{t\in[0,\mathcal{T}]}\|\hat M(t)\|_{V^\prime}^{2p}\leq \liminf_{n\to\infty}\hat{\mathbb{E}}\sup_{t\in[0,\mathcal{T}]}\|\hat M_n(t)\|_{V^\prime}^{2p}\leq C(\mathcal{T})<\infty.
    \end{equation}

    \textbf{Step 2: the martingale property.} Choose a deterministic dense set $\mathbb C\subset[0,\mathcal T]$ containing $0$ and $\mathcal T$, such that $\hat{\bm u}$ is almost surely continuous at every point of $\mathbb C\cap(0,\mathcal T)$. Such a set exists because a c\`adl\`ag process has at most countably many fixed times with positive probability of a jump. To see this, for each jump threshold restrict to paths having at most $k$ jumps above that threshold; the expected counting measure is finite. Taking the union over $k$ and over reciprocal integer thresholds proves the assertion. Endpoint evaluation is continuous in the Skorokhod topology, including when there is a jump at $\mathcal T$. Fix $s,t\in\mathbb C$, $s<t$, and $\psi\in V$. By \eqref{conver:M_n-hat},
    \begin{equation}\label{est:martingale-lem-3-step2-1}
        \lim_{n\to\infty}{}_{V^\prime}\langle\hat M_n(t)-\hat M_n(s),\psi\big\rangle_V={}_{V^\prime}\langle\hat M(t)-\hat M(s),\psi\rangle_V, \quad \hat{\mathbb P}\text{-a.s.}
    \end{equation}

    Let $h(x)=H(x|_{[0,s]})$, where $H\in C_b(D([0,s];V'))$. This is $\mathscr D_s^0$-measurable. Restriction to $[0,s]$ is continuous at paths continuous at $s$; for $s=0$ use endpoint evaluation. Therefore
    \begin{equation}\label{est:martingale-lem-3-step2-2}
        \lim_{n\to\infty}h(\hat{\bm u}_n)=h(\hat{\bm u}),\quad \hat{\mathbb P}\text{-a.s.}
    \end{equation}
    Set
    \begin{equation*}
        \begin{aligned}
            Y_n(\omega):= & {}_{V^\prime}\langle \hat M_n(t,\omega)-\hat M_n(s,\omega),\psi\rangle_V\,h(\hat{\bm u}_n), \\
            Y(\omega):= & {}_{V^\prime}\langle \hat M(t,\omega)-\hat M(s,\omega),\psi\rangle_V\,h(\hat{\bm u}),\quad \omega\in\hat{\Omega}.
        \end{aligned}
    \end{equation*}
    By the triangle and the Cauchy--Schwarz inequalities, together with \eqref{est:uni-hat-M-n}, we have
    \begin{equation*}
        \sup_{n\geq1}\hat{\mathbb{E}}|Y_n|^2\leq 4|h|_{L^\infty}^2\|\psi\|_V^2\sup_{n\geq1}\hat{\mathbb{E}}\sup_{s\in[0,\mathcal{T}]}\|\hat M_n(s)\|_{V^\prime}^2<\infty,
    \end{equation*}
    so $\{Y_n\}_{n\geq1}$ is uniformly integrable; hence, by \eqref{est:martingale-lem-3-step2-1}, \eqref{est:martingale-lem-3-step2-2} and Vitali's convergence theorem, we have
    \begin{equation*}
        \hat{\mathbb{E}}Y =\lim_{n\to\infty}\hat{\mathbb{E}}Y_n=0,
    \end{equation*}
    i.e., for all $s,t\in\mathbb C$ with $s<t$, $\psi\in V$ and past-path tests $h(x)=H(x|_{[0,s]})$ as above, we have
    \begin{equation}
        \hat{\mathbb{E}}({}_{V^\prime}\langle \hat M(t,\omega)-\hat M(s,\omega),\psi\rangle_V\,h(\hat{\bm u}))=0.
    \end{equation}
    Bounded continuous functions on $D([0,s];V')$ generate its Borel $\sigma$-field. The functional monotone class theorem therefore extends this identity to every bounded $\mathscr D_s^0$-measurable $h$, for $s,t\in\mathbb C$.

    Next, we extend the martingale property to the whole interval $[0,\mathcal{T}]$ by right continuity.

    Now let $0\leq s<t\leq\mathcal{T}$ be arbitrary. Since $\mathbb C$ is dense in $[0,\mathcal{T}]$, we may choose sequences $\{s_k\}_{k\geq1},\{t_k\}_{k\geq1}\subset\mathbb C$ such that
    \begin{equation*}
        0\leq s\leq s_k<t\leq t_k\leq \mathcal{T},\quad\forall\,k\geq 1,
    \end{equation*}
    taking $t_k=\mathcal T$ when $t=\mathcal T$. Moreover,
    \begin{equation*}
        s_k\downarrow s,\quad t_k\downarrow t,\quad\text{as }k\to\infty.
    \end{equation*}
    We also choose and fix $h:D([0,\mathcal{T}];V^\prime)\to\mathbb{R}$, a bounded $\mathscr D_s^0$-measurable function. By the right continuity of the paths of the process ${}_{V^\prime}\langle\hat{M},\psi\rangle_V$, we have, $\hat{\mathbb{P}}$-a.s.,
    \begin{align}
        \lim_{k\to\infty}{}_{V^\prime}\langle \hat M(s_k),\psi\rangle_V\,h(\hat{\bm u}) & ={}_{V^\prime}\langle \hat M(s),\psi\rangle_V\,h(\hat{\bm u})\label{est:martingale-lem-3-step2-4}, \\
        \lim_{k\to\infty}{}_{V^\prime}\langle \hat M(t_k),\psi\rangle_V\,h(\hat{\bm u}) & ={}_{V^\prime}\langle \hat M(t),\psi\rangle_V\,h(\hat{\bm u})\label{est:martingale-lem-3-step2-5}.
    \end{align}
    Also, by \eqref{est:martingale-lem-3-step1-1} and the boundedness of $h$, we have
    \begin{equation*}
        \hat{\mathbb{E}}\sup_{t\in [0,\mathcal{T}]}|{}_{V^\prime}\langle \hat M(t),\psi\rangle_V\,h(\hat{\bm u})|^2 \leq |h|_{L^\infty}^2\|\psi\|_V^2\,\hat{\mathbb{E}}\sup_{t\in[0,\mathcal{T}]}\|\hat M(t)\|_{V^\prime}^2<\infty,
    \end{equation*}
    which implies that
    \begin{equation*}
        \sup_{k\geq1}\hat{\mathbb{E}}|{}_{V^\prime}\langle \hat M(s_k),\psi\rangle_V\,h(\hat{\bm u})|^2+\sup_{k\geq1}\hat{\mathbb{E}}|{}_{V^\prime}\langle \hat M(t_k),\psi\rangle_V\,h(\hat{\bm u})|^2<\infty.
    \end{equation*}
    Thus the sequences $\{{}_{V^\prime}\langle \hat M(s_k),\psi\rangle_V\,h(\hat{\bm u})\}_{k\geq1}$ and $\{{}_{V^\prime}\langle \hat M(t_k),\psi\rangle_V\,h(\hat{\bm u})\}_{k\geq1}$ are uniformly integrable. Combining this with \eqref{est:martingale-lem-3-step2-4} and \eqref{est:martingale-lem-3-step2-5} and applying Vitali's convergence theorem, we obtain
    \begin{equation*}
        \hat{\mathbb{E}}({}_{V^\prime}\langle \hat M(t)-\hat M(s),\psi\rangle_V\,h(\hat{\bm u}))=\lim_{k\to\infty}\hat{\mathbb{E}}({}_{V^\prime}\langle \hat M(t_k)-\hat M(s_k),\psi\rangle_V\,h(\hat{\bm u}))=0.
    \end{equation*}
    Estimate \eqref{est:uni-hat-M} follows from \eqref{est:martingale-lem-3-step1-1}.

    The identities for all bounded past-measurable tests establish the martingale property for $\{\mathcal H_t\}_{t\geq0}$. The integrable supremum established above permits the augmentation argument in Lemma \ref{lem:martingale-1}, yielding the stated $\hat{\mathbb F}$-martingale property.
\end{proof}

\begin{lemma}\label{lem:martingale-4}
    For every $\psi\in V$, the process $\hat N_\psi$ is a square-integrable real-valued $\hat{\mathbb F}$-martingale.
\end{lemma}

\begin{proof}
    Use the deterministic set $\mathbb C$ chosen in the proof of Lemma \ref{lem:martingale-3}; evaluation at $\mathcal T$ is handled by endpoint convergence. The proof is divided into two steps.

    \medskip

    \textbf{Step 1: convergence and integrability of $\hat N_\psi$.}

    We first prove the following claim.

    \medskip

    \textbf{Claim 1.} \emph{For every $t\in\mathbb C$,
    \begin{equation}\label{conv:N-pointwise}
        \hat N_{n,\psi}(t)\to\hat N_\psi(t),\quad\hat{\mathbb P}\text{-a.s.}
    \end{equation}
    }

    \medskip

    \noindent \emph{Proof of Claim 1.} Fix $t\in\mathbb C$. Since $0,t\in\mathbb C$ and
    \begin{equation*}
        \hat{\bm u}_n\to\hat{\bm u}\quad\text{in }D([0,\mathcal{T}];V^\prime),
    \end{equation*}

    we have
    \begin{equation*}
        \hat{\bm u}_n(t)\to\hat{\bm u}(t),\quad\hat{\bm u}_n(0)\to\hat{\bm u}(0) \quad\text{in }V^\prime.
    \end{equation*}
    Hence
    \begin{equation}\label{conv:N-quad}
        {}_{V^\prime}\langle\hat{\bm u}_n(t),\psi\rangle_V^2- {}_{V^\prime}\langle\hat{\bm u}_n(0),\psi\rangle_V^2 \to{}_{V^\prime}\langle\hat{\bm u}(t),\psi\rangle_V^2-{}_{V^\prime}\langle\hat{\bm u}(0),\psi\rangle_V^2 .
    \end{equation}

    We next consider the compensator term. By \eqref{eq:unif-a.s.}, for almost every $\omega$ there is $R(\omega)<\infty$ such that
    \begin{equation*}
        \sup_n\|\hat{\bm u}_n\|_{L^\infty(0,\mathcal T;V)}+\|\hat{\bm u}\|_{L^\infty(0,\mathcal T;V)}\leq R(\omega).
    \end{equation*}
    For this path, \eqref{asm:g_Lip-local} gives
    \begin{equation*}
        \int_0^{\mathcal T}\int_E|G(\hat{\bm u}_n,\xi)-G(\hat{\bm u},\xi)|_H^2\,\nu(\mathrm{d}\xi)\,\mathrm{d}s
        \leq C_{R(\omega)}\|\hat{\bm u}_n-\hat{\bm u}\|_{L^2(0,\mathcal T;H)}^2\to0.
    \end{equation*}
    Values may be replaced by left limits in these time integrals. Decompose the projected difference as
    \begin{equation*}
        P_nG(\hat{\bm u}_n)-G(\hat{\bm u})=P_n[G(\hat{\bm u}_n)-G(\hat{\bm u})]+(P_n-I)G(\hat{\bm u}).
    \end{equation*}
    The first term converges in $L^2((0,\mathcal T)\times E;H)$ by the contraction property. The second converges in the same space by dominated convergence and \eqref{asm:g_linear-H}. Cauchy--Schwarz therefore gives, almost surely,
    \begin{equation}\label{conv:N-compensator}
        \begin{aligned}
            & \int_0^t\int_E{}_{V^\prime}\langle G_{(n)}(\hat{\bm u}_n(s-),\xi),\psi\rangle_V^2\,\nu(\mathrm{d}\xi)\,\mathrm{d}s \\
            \to & \int_0^t\int_E{}_{V^\prime}\langle G(\hat{\bm u}(s-),\xi),\psi\rangle_V^2\,\nu(\mathrm{d}\xi)\,\mathrm{d}s,\quad\text{as}\ n\to\infty.
        \end{aligned}
    \end{equation}

    We finally consider the cross term. We decompose the difference as
    \begin{equation*}
        \begin{aligned}
            & \int_0^t{}_{V^\prime}\langle A\hat{\bm u}_n+\widehat{\mathcal B}_n(s)+R_{(n)}(\hat{\bm u}_n),\psi \rangle_V\,{}_{V^\prime}\langle\hat{\bm u}_n,\psi\rangle_V\,\mathrm{d}s \\
            & -\int_0^t{}_{V^\prime}\langle A\hat{\bm u}+\widehat{\mathcal B}(s)+R(\hat{\bm u}),\psi\rangle_V\,{}_{V^\prime}\langle\hat{\bm u},\psi\rangle_V\,\mathrm{d}s \\
            = & \int_0^t{}_{V^\prime}\langle A(\hat{\bm u}_n-\hat{\bm u}),\psi\rangle_V\,{}_{V^\prime}\langle\hat{\bm u},\psi\rangle_V \mathrm{d}s+\int_0^t{}_{V^\prime}\langle \widehat{\mathcal B}_n(s)-\widehat{\mathcal B}(s),\psi\rangle_V\,{}_{V^\prime}\langle\hat{\bm u},\psi\rangle_V\,\mathrm{d}s \\
            & +\int_0^t{}_{V^\prime}\langle R_{(n)}(\hat{\bm u}_n)-R(\hat{\bm u}),\psi\rangle_V\,{}_{V^\prime}\langle\hat{\bm u},\psi\rangle_V \mathrm{d}s \\
            & +\int_0^t{}_{V^\prime}\langle A\hat{\bm u}_n+\widehat{\mathcal B}_n(s)+R_{(n)}(\hat{\bm u}_n),\psi\rangle_V\,{}_{V^\prime}\langle\hat{\bm u}_n-\hat{\bm u},\psi\rangle_V\,\mathrm{d}s \\
            := & J_{1,n}+J_{2,n}+J_{3,n}+J_{4,n}.
        \end{aligned}
    \end{equation*}
    The scalar drift differences in $J_{1,n}+J_{2,n}+J_{3,n}$ are bounded in $L^2(0,\mathcal T)$ by \eqref{eq:unif-a.s.}, \eqref{est:B_op-dual-norm} and the bounds for $A$ and $R$. Their primitives converge to zero by the drift convergence used in \eqref{conver:M_n-hat}. Testing first against step functions and then using their density in $L^2(0,\mathcal T)$ proves weak $L^2$ convergence to zero. Since $\mathbbm1_{[0,t]}\langle\hat{\bm u},\psi\rangle$ belongs to $L^2(0,\mathcal T)$, it follows that
    \begin{equation}\label{conv:N-cross-1-2-3}
        J_{1,n}+J_{2,n}+J_{3,n}\to 0,\quad\text{as}\ n\to\infty.
    \end{equation}

    It remains to show that
    \begin{equation*}
        J_{4,n}:=\int_0^t{}_{V^\prime}\langle A\hat{\bm u}_n+\widehat{\mathcal B}_n(s)+R_{(n)}(\hat{\bm u}_n),\psi\rangle_V\,{}_{V^\prime}\langle\hat{\bm u}_n-\hat{\bm u},\psi\rangle_V\,\mathrm{d}s\to 0.
    \end{equation*}
    By \eqref{est:a} and \eqref{est:r}, and the strong convergence of $\hat{\bm u}_n$ to $\hat{\bm u}$ in $L^2(0,\mathcal{T};H)$, we have, $\hat{\mathbb P}$-a.s.,
    \begin{equation}\label{conv:N-cross-4-AR}
        \begin{aligned}
            & \left|\int_0^t{}_{V^\prime}\langle A\hat{\bm u}_n+R_{(n)}(\hat{\bm u}_n),\psi\rangle_V\,{}_{V^\prime}\langle\hat{\bm u}_n-\hat{\bm u},\psi\rangle_V\,\mathrm{d}s\right| \\
            \leq & C\|\psi\|_V^2\int_0^t\|\hat{\bm u}_n\|_V|\hat{\bm u}_n-\hat{\bm u}|_H\,\mathrm{d}s \\
            \leq & C\|\psi\|_V^2\|\hat{\bm u}_n\|_{L^2(0,\mathcal{T};V)}\|\hat{\bm u}_n-\hat{\bm u}\|_{L^2(0,\mathcal{T};H)}\to 0.
        \end{aligned}
    \end{equation}

    By \eqref{eq:unif-a.s.} and \eqref{est:B_op-dual-norm}, $\sup_n\|B(\hat{\bm u}_n)\|_{L^2(0,\mathcal T;V')}\leq C(\omega)$. Consequently,
    \begin{equation}\label{conv:N-cross-4-B}
        \begin{aligned}
            & \left|\int_0^t{}_{V^\prime}\langle \widehat{\mathcal B}_n(s),\psi\rangle_V\,{}_{V^\prime}\langle\hat{\bm u}_n-\hat{\bm u},\psi\rangle_V\,\mathrm{d}s\right| \\
            & \leq C\|\psi\|_V^2\int_0^t\Big(|\hat{\bm u}_n|_H^{1/2}\|\hat{\bm u}_n\|_V^{3/2}+\|\hat{\bm u}_n\|_V^{3/2}\|\hat{\bm u}_n\|_{\mathscr D(A)}^{1/2}\Big)|\hat{\bm u}_n-\hat{\bm u}|_H\,\mathrm{d}s \\
            & \leq C(\omega)\|\psi\|_V^2\|\hat{\bm u}_n-\hat{\bm u}\|_{L^2(0,\mathcal T;H)}\to0.
        \end{aligned}
    \end{equation}
    Therefore, by \eqref{conv:N-cross-1-2-3}, \eqref{conv:N-cross-4-AR} and \eqref{conv:N-cross-4-B}, we conclude that $\hat{\mathbb P}\text{-a.s.}$,
    \begin{equation}\label{conv:N-cross}
        \begin{aligned}
            & \int_0^t{}_{V^\prime}\langle A\hat{\bm u}_n+\widehat{\mathcal B}_n(s)+R_{(n)}(\hat{\bm u}_n),\psi\rangle_V\,{}_{V^\prime}\langle\hat{\bm u}_n,\psi\rangle_V\,\mathrm{d}s \\
            \to & \int_0^t{}_{V^\prime}\langle A\hat{\bm u}+\widehat{\mathcal B}(s)+R(\hat{\bm u}),\psi\rangle_V\,{}_{V^\prime}\langle\hat{\bm u},\psi\rangle_V \mathrm{d}s,\quad\text{as }n\to\infty.
        \end{aligned}
    \end{equation}
    Hence, by \eqref{conv:N-quad}, \eqref{conv:N-compensator} and \eqref{conv:N-cross}, we finish the proof of Claim 1.

    \medskip

    \textbf{Claim 2.} \emph{
    \begin{equation}
        \hat{\mathbb E}\sup_{t\in[0,\mathcal{T}]}|\hat N_\psi(t)|^2<\infty.
    \end{equation}
    }
    \medskip

    \noindent \emph{Proof of Claim 2.}

    Choose a countable dense subset $\mathbb C_0\subset\mathbb C$ containing $\mathcal T$. Claim 1 holds simultaneously on $\mathbb C_0$ outside a single null set. Hence Fatou's lemma and \eqref{est:N_n,psi} give
    \begin{equation}\label{est:N-fixed-time}
        \hat{\mathbb E}\sup_{t\in\mathbb C_0}|\hat N_\psi(t)|^2 \leq\liminf_n\hat{\mathbb E}\sup_{t\in\mathbb C_0}|\hat N_{n,\psi}(t)|^2 \leq C_{\mathcal T}\|\psi\|_V^4.
    \end{equation}
    Right continuity and $\mathcal T\in\mathbb C_0$ identify the supremum over $\mathbb C_0$ with that over $[0,\mathcal T]$. Therefore
    \begin{equation}\label{est:N-hat-quad}
        \hat{\mathbb E}\sup_{t\in[0,\mathcal T]}|\hat N_\psi(t)|^2 \leq C_{\mathcal T}\|\psi\|_V^4.
    \end{equation}

    \medskip

    \textbf{Step 2: the martingale property.} For $s,t\in\mathbb C$ with $s<t$, \eqref{conv:N-pointwise} gives
    \begin{equation}\label{est:martingale-lem-4-step2-1}
        \lim_{n\to\infty}(\hat N_{n,\psi}(t)-\hat N_{n,\psi}(s))=\hat N_\psi(t)-\hat N_\psi(s), \quad \hat{\mathbb P}\text{-a.s.}
    \end{equation}

    Let $h(x)=H(x|_{[0,s]})$ with $H\in C_b(D([0,s];V'))$. Then $h$ is $\mathscr D_s^0$-measurable and converges along the represented paths as in \eqref{est:martingale-lem-3-step2-2}. Denote
    \begin{equation*}
        \begin{aligned}
            Z_n(\omega) & :=(\hat N_{n,\psi}(t,\omega)-\hat N_{n,\psi}(s,\omega))h(\hat{\bm u}_n), \\
            Z(\omega) & :=(\hat N_\psi(t,\omega)-\hat N_\psi(s,\omega))h(\hat{\bm u}),\quad\omega\in\hat{\Omega}.
        \end{aligned}
    \end{equation*}
    With the same argument as in the proof of Lemma \ref{lem:martingale-3}, we can show that $\{Z_n\}_{n\geq1}$ is uniformly integrable. Thus, by \eqref{est:martingale-lem-4-step2-1}, \eqref{est:martingale-lem-3-step2-2} and Vitali's convergence theorem, we have, for all $s,t\in\mathbb C$ with $s<t$, $\psi\in V$ and past-path tests $h(x)=H(x|_{[0,s]})$ as above,
    \begin{equation}
        \hat{\mathbb{E}}((\hat N_\psi(t)-\hat N_\psi(s))h(\hat{\bm u}))=\hat{\mathbb{E}}Z=\lim_{n\to\infty}\hat{\mathbb{E}}Z_n=0.
    \end{equation}

    The functional monotone class theorem extends this identity to every bounded $\mathscr D_s^0$-measurable $h$, since the preceding tests generate the past-path Borel $\sigma$-field.

    Next, we extend the martingale property to the whole interval $[0,\mathcal{T}]$ by right continuity. Let $0\leq s<t\leq\mathcal{T}$ be arbitrary. Since $\mathbb C$ is dense in $[0,\mathcal{T}]$, we may choose sequences $\{s_k\}_{k\geq1},\{t_k\}_{k\geq1}\subset\mathbb C$ such that
    \begin{equation*}
        0\leq s\leq s_k<t\leq t_k\leq \mathcal{T},\quad\forall\,k\geq 1.
    \end{equation*}
    Moreover,
    \begin{equation*}
        s_k\downarrow s,\quad t_k\downarrow t,\quad\text{as }k\to\infty.
    \end{equation*}
    We also choose and fix $h:D([0,\mathcal{T}];V^\prime)\to\mathbb{R}$, a bounded $\mathscr D_s^0$-measurable function. By the right continuity of the paths of the process $\hat N_\psi$, we have, $\hat{\mathbb{P}}$-a.s.,
    \begin{align}
        \lim_{k\to\infty}\hat N_\psi(s_k)\,h(\hat{\bm u}) & =\hat N_\psi(s)\,h(\hat{\bm u})\label{est:martingale-lem-4-step2-4}, \\
        \lim_{k\to\infty}\hat N_\psi(t_k)\,h(\hat{\bm u}) & =\hat N_\psi(t)\,h(\hat{\bm u})\label{est:martingale-lem-4-step2-5}.
    \end{align}
    Also, by \eqref{est:N-hat-quad} and the boundedness of $h$, we have
    \begin{equation*}
        \hat{\mathbb{E}}\sup_{t\in [0,\mathcal{T}]}|\hat N_\psi(t)\,h(\hat{\bm u})|^2 \leq |h|_{L^\infty}^2\,\hat{\mathbb{E}}\sup_{t\in[0,\mathcal{T}]}|\hat N_\psi(t)|^2<\infty,
    \end{equation*}
    which implies that
    \begin{equation*}
        \sup_{k\geq1}\hat{\mathbb{E}}|\hat N_\psi(s_k)\,h(\hat{\bm u})|^2+\sup_{k\geq1}\hat{\mathbb{E}}|\hat N_\psi(t_k)\,h(\hat{\bm u})|^2<\infty.
    \end{equation*}
    The two sequences $\{\hat N_\psi(s_k)\,h(\hat{\bm u})\}_{k\geq1}$ and $\{\hat N_\psi(t_k)\,h(\hat{\bm u})\}_{k\geq1}$ are uniformly integrable. By \eqref{est:martingale-lem-4-step2-4} and \eqref{est:martingale-lem-4-step2-5}, Vitali's convergence theorem gives
    \begin{equation*}
        \hat{\mathbb{E}}((\hat N_\psi(t)-\hat N_\psi(s))h(\hat{\bm u}))=\lim_{k\to\infty}\hat{\mathbb{E}}((\hat N_\psi(t_k)-\hat N_\psi(s_k))h(\hat{\bm u}))=0.
    \end{equation*}

    The identities for all bounded past-measurable tests establish the martingale property for $\{\mathcal H_t\}_{t\geq0}$. The integrable supremum established above permits the augmentation argument in Lemma \ref{lem:martingale-1}, yielding the stated $\hat{\mathbb F}$-martingale property.
\end{proof}

\subsubsection{Identification of the quadratic variation}
Fix $\psi\in V$ and define a real-valued process $\hat M_\psi$ by
\begin{equation}\label{def:hat-M-psi}
    \begin{aligned}
        \hat{M}_\psi(t):={}_{V^\prime}\langle\hat{M}(t),\psi\rangle_V= & {}_{V^\prime}\langle\hat{\bm u}(t),\psi\rangle_V-{}_{V^\prime}\langle\hat{\bm u}(0),\psi\rangle_V \\
        & +\int_0^t{}_{V^\prime}\langle A\hat{\bm u}(s)+\widehat{\mathcal B}(s)+R(\hat{\bm u}(s)),\psi\rangle_V\,\mathrm{d}s,\quad t\in[0,\mathcal T],
    \end{aligned}
\end{equation}
where $\hat M$ is the $V^\prime$-valued martingale defined in \eqref{eqn:def-hat-M}. By Lemma \ref{lem:martingale-3}, $\hat M_\psi$ is a square-integrable real-valued $\hat{\mathbb F}$-martingale.

\begin{lemma}
    For every $\psi\in V$ the predictable quadratic variation of the martingale $\hat M_\psi$ is given by
    \begin{equation}\label{eqn:qv}
        \langle\hat M_\psi\rangle(t)=\int_0^t\int_E {}_{V^\prime}\langle G(\hat{\bm u}(r),\xi),\psi\rangle_V^2\,\nu(\mathrm{d}\xi)\,\mathrm{d}r,\quad t\in[0,\mathcal T],\quad\hat{\mathbb P}\text{-a.s.}
    \end{equation}
\end{lemma}

\begin{proof}
    The proof follows the approach of Kallianpur and Xiong; see \cite[Lemma~6.1.10]{0859.60050} and \cite[Proposition~6.11]{8220046}. The process defined by the right-hand side of \eqref{eqn:qv} is continuous, $\hat{\mathbb F}$-adapted and of finite variation; in particular, it is $\hat{\mathbb F}$-predictable. Hence, by the uniqueness of the predictable compensator in the Doob--Meyer decomposition, it is enough to prove that the process
    \begin{equation}\label{eqn:qv-claim}
        K_\psi(t):=(\hat M_\psi(t))^2-\int_0^t\int_E{}_{V^\prime}\langle G(\hat{\bm u}(r),\xi),\psi\rangle_V^2\,\nu(\mathrm{d}\xi)\,\mathrm{d}r,\quad t\in[0,\mathcal T],
    \end{equation}
    is a martingale. Put
    \begin{equation*}
        \zeta(t):={}_{V^\prime}\langle\hat{\bm u}(t),\psi\rangle_V, \qquad t\in[0,\mathcal T].
    \end{equation*}
    Since the trajectories of $\hat{\bm u}$ are c\`adl\`ag in $V^\prime$, the process $\zeta$ is locally bounded and, by \eqref{def:hat-M-psi},
    \begin{equation}
        \zeta(t)=\zeta(0)-\int_0^t{}_{V^\prime}\langle A\hat{\bm u}(r)+\widehat{\mathcal B}(r)+R(\hat{\bm u}(r)),\psi\rangle_V\,\mathrm{d}r+\hat M_\psi(t).
    \end{equation}
    Applying It\^o's formula to the real-valued semimartingale $\zeta$ and using that the finite variation part of $\zeta$ is continuous, so that $[\zeta,\zeta](t)=[\hat M_\psi,\hat M_\psi](t)$, we obtain
    \begin{equation}\label{eq:qv-ito}
        \begin{aligned}
            \zeta(t)^2 & =\zeta(0)^2-2\int_0^t{}_{V^\prime}\langle A\hat{\bm u}(r)+\widehat{\mathcal B}(r)+R(\hat{\bm u}(r)),\psi\rangle_V\,\zeta(r)\,\mathrm{d}r \\
            & \quad+2\int_0^t\zeta(r-)\,\mathrm{d}\hat M_\psi(r)+[\hat M_\psi,\hat M_\psi](t).
        \end{aligned}
    \end{equation}
    Since $\zeta(r-)$ is predictable and locally bounded, the stochastic integral in \eqref{eq:qv-ito} is a local martingale. Comparing \eqref{eq:qv-ito} with the definition \eqref{eqn:def-hat-N-psi} of $\hat N_\psi$, we infer that the process
    \begin{equation*}
        \hat N_\psi(t)-\left([\hat M_\psi,\hat M_\psi](t)-\int_0^t\int_E {}_{V^\prime}\langle G(\hat{\bm u}(r),\xi),\psi\rangle_V^2 \,\nu(\mathrm{d}\xi)\,\mathrm{d}r\right),\quad t\in[0,\mathcal T],
    \end{equation*}
    is a local martingale. Since $\hat N_\psi$ is a martingale by Lemma \ref{lem:martingale-4}, we deduce that
    \begin{equation*}
        [\hat M_\psi,\hat M_\psi](t)-\int_0^t\int_E{}_{V^\prime}\langle G(\hat{\bm u}(r),\xi),\psi\rangle_V^2\,\nu(\mathrm{d}\xi)\,\mathrm{d}r,\quad t\in[0,\mathcal T],
    \end{equation*}
    is a local martingale as well. Since $(\hat M_\psi)^2-[\hat M_\psi,\hat M_\psi]$ is a local martingale by \cite[Corollary~18.9, p.~125]{0503.60054}, adding the two processes shows that $K_\psi$ defined in \eqref{eqn:qv-claim} is a local martingale. Moreover, by Lemma \ref{lem:martingale-3}, the linear growth condition \eqref{asm:g_linear-H} and the energy estimate, we have
    \begin{equation*}
        \hat{\mathbb E}\sup_{t\in[0,\mathcal T]}|K_\psi(t)|\leq\hat{\mathbb E}\sup_{t\in[0,\mathcal T]}|\hat M_\psi(t)|^2+C\|\psi\|_V^2\,\hat{\mathbb E}\int_0^{\mathcal T}(1+|\hat{\bm u}(t)|_H^2)\,\mathrm{d}t<\infty,
    \end{equation*}
    so that $K_\psi$ is a martingale.
\end{proof}

\subsubsection{Identification of the retained Poisson measure}\label{subsec:representation}
The following argument applies to a single Galerkin sequence and to any finite family of sequences driven by the same noise.

\begin{lemma}\label{lem:retained-noise}
    Let $d\in\mathbb N$ and $r_k^i\to\infty$ for $1\leq i\leq d$. Under Assumptions \ref{asm:martingale} and \ref{asm:A-est}, suppose $\bm u_0\in L^4(\Omega,\mathcal F_0;V)$. Consider a joint Skorokhod--Jakubowski representation
    \begin{equation}\label{eq:retained-joint-convergence}
        (a_k,u_k^1,\ldots,u_k^d,\eta_k) \to(a,u^1,\ldots,u^d,\eta^*) \quad\text{almost surely in }V\times\mathcal Z_{\mathcal T}^d\times\mathcal N,
    \end{equation}
    where the $k$th tuple has the law of $(\bm u_0,\bm u_{r_k^1},\ldots,\bm u_{r_k^d},\eta)$. Let $\mathbb H:=\{\mathcal H_t\}_{t\geq0}$ be the usual augmentation of
    \begin{equation}\label{eq:joint-limit-filtration}
        \mathcal H_t^0=\sigma\{a,u^i(r),\eta^*((0,r]\times A): 1\leq i\leq d,\ r\leq t,\ A\in\mathcal B(E),\ \nu(A)<\infty\}.
    \end{equation}
    Then $\eta^*$ is an $\mathbb H$-Poisson random measure with compensator $\mathrm{d}t\otimes\nu$, and each $u^i$ solves the regularised equation with initial datum $a$ and driver $\eta^*$. Moreover, $u^i\in\mathfrak C([0,\mathcal T])$ and has $V$-valued c\`adl\`ag paths.
\end{lemma}
\begin{proof}
    Put $u_k^{*,i}(t)=e^{-tA}P_{r_k^i}a_k$, $1\leq i\leq d$. Then $P_{r_k^i}a_k\to a$ in $V$, and the background evolutions converge in $C([0,\mathcal T];V)$. Endpoint convergence gives $u^i(0)=a$. The bounds of Proposition \ref{prop:uniform-DA} pass to the limits by lower semicontinuity. Interpolation as in \eqref{eq:strong-L2-V} gives $u_k^i\to u^i$ in $L^2(0,\mathcal T;V)$ almost surely. The inherited uniform $V$ bound and the $V'$-valued c\`adl\`ag paths imply $H$-valued c\`adl\`ag paths by interpolation; their left limits belong to $V$. The Borel structure of $V$ is induced by its continuous injection into $H$; hence the $V$-valued left-limit process is predictable, and so is $G(u^i(\cdot-),\cdot)$.

    \emph{Compatibility of the noise and the joint past.} The marginal law of $\eta^*$ is that of $\eta$. Fix $s<t$ at deterministic continuity times for all limiting paths, and let $F_k$ be a bounded continuous cylinder function of $a_k$, the paths up to $s$, and the noise up to $s$. For nonnegative bounded continuous $f$ supported on finitely many strips, the original independent-increment property and equality of joint laws give
    \begin{equation*}
        \mathbb E_*\left[F_k e^{-\int_{(s,t]\times E}f\,\mathrm{d}\eta_k}\right] =\mathbb E_*F_k\, \exp\left(\int_s^t\int_E(e^{-f}-1)\,\nu(\mathrm{d}\xi)\,\mathrm{d}r\right).
    \end{equation*}
    Here $\mathbb E_*$ denotes expectation on the representation space. Noise restrictions converge at deterministic times, since their limit has no atoms at these times. Passing to the limit and using the monotone class theorem gives the conditional Poisson Laplace functional relative to $\mathcal H_s^0$. Right continuity extends this identity to all $s<t$ and to the usual augmentation, as in the augmentation argument of Lemma \ref{lem:martingale-1}. Thus the future noise is independent of the joint past, with the stated compensator.

    \emph{Joint martingale identities.} For $1\leq i\leq d$ define the drift residual
    \begin{equation*}
        M^i(t)=u^i(t)-a+\int_0^t [Au^i+\theta(\|u^i-e^{-rA}a\|_V)B(u^i)+R(u^i)]\mathrm{d}r,
    \end{equation*}
    and let $M_k^i$ be its Galerkin counterpart. For $\psi\in V$, set $M_\psi^i={}_{V'}\langle M^i,\psi\rangle_V$ and
    \begin{equation*}
        b_\psi^i(r,\xi)={}_{V'}\langle G(u^i(r-),\xi),\psi\rangle_V, \qquad K_f(t)=\int_0^t\int_E f(r,\xi)\,\widetilde\eta^*(\mathrm{d}r,\mathrm{d}\xi).
    \end{equation*}
    Initially take deterministic bounded continuous $f$ on finitely many strips. On the original basis, the Galerkin stochastic integrals satisfy the martingale identities for $M_{k,\psi}^i$, its square minus its compensator, and its product with $K_{k,f}$ minus their mixed compensator. These identities transfer by the joint law, using tests of the joint past.

    We justify passage to the limit. The drift converges as in \eqref{eq:convergence-B}, retaining the cutoff and the converging backgrounds. Uniform fourth moments and almost sure $L^2H$ convergence give $u_k^i\to u^i$ in $L^2(\Omega_*\times[0,\mathcal T];H)$. Proposition \ref{prop:uniform-DA} and the background estimate give uniform fourth $V$-moments, also for the limit by lower semicontinuity. Lemma \ref{lem:local-noise-convergence} therefore gives convergence of $G(u_k^i)$ to $G(u^i)$ in $L^2(\Omega_*\times(0,\mathcal T)\times E;H)$. The contraction property of the projections and dominated convergence for $(P_{r_k^i}-I)G(u^i)$ then yield
    \begin{equation}\label{eq:joint-noise-coefficient-convergence}
        P_{r_k^i}G(u_k^i)-G(u^i) \to0\quad\text{in } L^2(\Omega_*\times[0,\mathcal T]\times E;H).
    \end{equation}
    Time integrals are unchanged by replacing paths with their left limits. Thus both the quadratic and mixed compensators converge in $L^1$. For the martingale products, \eqref{est:uni-hat-M-n} gives uniform fourth moments of the scalar residuals. The $K_{k,f}$ have a fixed law with finite fourth moments. Their products and the squares of the residuals are therefore uniformly integrable. The same is true of the quadratic compensators, by \eqref{asm:g_linear-H} and the fourth solution moment. The integrals $K_{k,f}(t)$ converge at deterministic times because the limiting point measures have no atoms on the time boundaries. For the solution coordinates, take a common dense set of deterministic almost sure continuity times and include $\mathcal T$, where Skorokhod endpoint convergence applies. Vitali convergence, followed by the past-path monotone class argument and right-continuous augmentation, proves that
    \begin{equation}\label{eq:joint-martingale-identities}
        M_\psi^i,\qquad (M_\psi^i)^2-\int_0^\cdot\int_E|b_\psi^i|^2\,\mathrm{d}\nu\,\mathrm{d}r,\qquad M_\psi^iK_f-\int_0^\cdot\int_E b_\psi^i f\,\mathrm{d}\nu\,\mathrm{d}r
    \end{equation}
    are martingales for the same augmented joint filtration. In particular,
    \begin{equation}\label{eq:joint-mixed-covariation}
        \langle M_\psi^i\rangle_t=\int_0^t\int_E|b_\psi^i|^2\,\mathrm{d}\nu\,\mathrm{d}r, \qquad \langle M_\psi^i,K_f\rangle_t =\int_0^t\int_E b_\psi^i f\,\mathrm{d}\nu\,\mathrm{d}r.
    \end{equation}
    Density and the Poisson isometry extend the second identity to deterministic $f\in L^2(\mathrm{d}t\otimes\nu)$.

    \emph{Identification of the retained driver.} Let $Z$ be bounded and measurable at time $s$ for the augmented joint filtration. The product martingale in \eqref{eq:joint-martingale-identities}, together with the martingale properties of $M_\psi^i$ and $K_f$, gives, for $s<t$,
    \begin{equation*}
        \mathbb E_*\big[M_\psi^i(t)Z(K_f(t)-K_f(s))\big] =\mathbb E_*\int_s^t\int_E Zb_\psi^i f\,\mathrm{d}\nu\,\mathrm{d}r.
    \end{equation*}
    Finite sums give this identity for elementary predictable integrands. Their density in the predictable subspace of $L^2(\Omega_*\times[0,\mathcal T]\times E)$, the Poisson isometry, and Cauchy--Schwarz then give, with $I_\psi^i=\int b_\psi^i\,\mathrm{d}\widetilde\eta^*$,
    \begin{equation*}
        \mathbb E_*[M_\psi^i(t)I_\psi^i(t)] =\mathbb E_*\int_0^t\int_E|b_\psi^i|^2\,\mathrm{d}\nu\,\mathrm{d}r.
    \end{equation*}
    The first identity in \eqref{eq:joint-mixed-covariation} and the isometry therefore yield
    \begin{equation}\label{eq:joint-noise-residual-zero}
        \mathbb E_*|M_\psi^i(t)-I_\psi^i(t)|^2=0.
    \end{equation}
    Taking a countable dense set of $\psi\in V$ and a countable dense set of times, then using c\`adl\`ag paths, proves
    \begin{equation}\label{eq:pair-common-noise}
        M^i(t)=\int_0^t\int_E G(u^i(r-),\xi)\, \widetilde\eta^*(\mathrm{d}r,\mathrm{d}\xi), \qquad 1\leq i\leq d,\quad 0\leq t\leq\mathcal T,
    \end{equation}
    almost surely. This identifies every equation with the retained noise. The inherited bounds and Lemma \ref{lem:strong-cadlag} give the asserted strong path regularity and membership in $\mathfrak C([0,\mathcal T])$.
\end{proof}

\begin{proposition}\label{prop:pure-disc}
    For every $\psi\in V$, the martingale $\hat M_\psi$ is purely discontinuous. In particular,
    \begin{equation}\label{eq:pure-disc}
        [\hat M_\psi,\hat M_\psi](t) =\sum_{0<r\leq t}|\Delta\hat M_\psi(r)|^2, \qquad 0\leq t\leq\mathcal T,\quad\hat{\mathbb P}\text{-almost surely}.
    \end{equation}
\end{proposition}
\begin{proof}
    Apply Lemma \ref{lem:retained-noise} with $d=1$ to \eqref{sk-rep-u-law}. It gives
    \begin{equation}\label{eq:representation}
        \hat M(t)=\int_0^t\int_EG(\hat{\bm u}(r-),\xi)\, \widetilde{\hat\eta}(\mathrm{d}r,\mathrm{d}\xi).
    \end{equation}
    A square-integrable compensated Poisson integral is purely discontinuous. Its scalar quadratic variation is the sum of its squared jumps, proving \eqref{eq:pure-disc}.
\end{proof}

\begin{proof}[Proof of Theorem \ref{thm:exi-martingale}]
    The representation \eqref{sk-rep-u-law} and Lemma \ref{lem:retained-noise} with $d=1$ give a solution driven by $\hat\eta$ on the same representation space, with the usual joint filtration \eqref{eq:joint-limit-filtration}. Substituting \eqref{eq:representation} into \eqref{eqn:def-hat-M} proves \eqref{eqn:u-hat-main}. The bounds follow from Corollary \ref{cor:limit-bounds}, and Lemma \ref{lem:strong-cadlag} gives the $V$-valued c\`adl\`ag version. The same tightness and joint identification on all integer time intervals, with a diagonal extraction, yield a single solution on $[0,\infty)$. All estimates hold on each finite interval. Finally, the right continuity in $V$ and $\hat{\bm u}(0)=\hat{\bm u}^*(0)$ imply that the exit time \eqref{def:tau-exit} is strictly positive. Before this exit the cutoff is one. Stopping the integral equation at the exit retains its terminal jump and gives the local solution of \eqref{eqn:abstract}.
\end{proof}

\section{Pathwise solutions and global well-posedness}\label{sec:pathwise}

Theorem \ref{thm:exi-martingale} gives a global martingale solution of the regularised equation and a local martingale solution of \eqref{eqn:abstract}. We first prove pathwise uniqueness, then use the Gyöngy--Krylov lemma to construct a maximal pathwise solution.

\subsection{Pathwise uniqueness}\label{subsec:unique}
Fix $\mathcal T>0$. We use the class $\mathfrak C([0,\mathcal T])$ of Definition \ref{def:unif-bound}.

\begin{theorem}[Pathwise uniqueness]\label{thm:unique}
    Let Assumptions \ref{asm:martingale} and \ref{asm:A-est} hold and let $\bm u^{(1)},\bm u^{(2)}\in \mathfrak C([0,\mathcal T])$ be two martingale solutions of \eqref{eqn:abstract} defined on the same stochastic basis and driven by the same compensated Poisson random measure $\tilde{\eta}$, with $\bm u^{(1)}(0)=\bm u^{(2)}(0)=\bm u_0$.  Then
    \begin{equation}\label{eq:pathwise-uniqueness}
        \mathbb P\big(\bm u^{(1)}(t)=\bm u^{(2)}(t)\ \text{for all } t\in[0,\mathcal T]\big)=1 .
    \end{equation}
\end{theorem}

\begin{remark}\label{rk:unique-regularised}
    For the regularised equation set $\theta_i=\theta(\|\bm u^{(i)}-\bm u^*\|_V)$ and $\bm\rho=\bm u^{(1)}-\bm u^{(2)}$. Its drift difference is $\theta_1[B(\bm u^{(1)})-B(\bm u^{(2)})]+(\theta_1-\theta_2)B(\bm u^{(2)})$. The first term obeys \eqref{est:B-diff_prod}, and
    \begin{equation*}
        |(\theta_1-\theta_2)\langle B(\bm u^{(2)}),\bm\rho\rangle_H| \leq\varepsilon\|\bm\rho\|_V^2+C_{\varepsilon,\theta}|B(\bm u^{(2)})|_H^2|\bm\rho|_H^2.
    \end{equation*}
    By \eqref{est:B-H-norm}, $|B(\bm u)|_H\leq C\|\bm u\|_V\|\bm u\|_{\mathscr D(A)}$, so this additional coefficient is time integrable for paths in $\mathfrak C$. Thus the uniqueness argument below also applies, after including this coefficient in the stopping time.
\end{remark}

\begin{proof}
    Set $\bm\rho:=\bm u^{(1)}-\bm u^{(2)}$.  Subtracting the two copies of \eqref{eqn:weak-form} we obtain, for every $\psi\in V$,
    \begin{equation}
        \begin{aligned}
            \langle\bm\rho(t),\psi\rangle_H & +\int_0^t{}_{V^\prime}\big\langle A\bm\rho +B(\bm u^{(1)})-B(\bm u^{(2)}) +R(\bm\rho),\psi\big\rangle_V\,\mathrm{d}s \\
            & =\int_0^t\int_E{}_{V^\prime}\big\langle G(\bm u^{(1)}(s-),\xi)-G(\bm u^{(2)}(s-),\xi),\psi\big\rangle_V \tilde{\eta}(\mathrm{d}s,\mathrm{d}\xi).
        \end{aligned}
    \end{equation}
    By Lemma \ref{lem:strong-cadlag} the solutions have $V$-valued c\`adl\`ag paths. Estimate \eqref{est:B-H-symmetric} and the local energy bound put their drifts in $L^2(0,\mathcal T;H)$ almost surely. For $r\geq1$, set
    \begin{equation}\label{eq:uniqueness-norm-stop}
        \sigma_r=\inf\{t\in[0,\mathcal T]:\|\bm u^{(1)}(t)\|_V\vee\|\bm u^{(2)}(t)\|_V>r\}\wedge\mathcal T,
    \end{equation}
    with $\inf\varnothing=\infty$. On $(0,\sigma_r]$ the left limits of both processes have $V$ norm at most $r$. All subsequent identities are used on $[0,\sigma_r]$, with the terminal jump retained. The stopped Hilbert-space It\^o formula gives
    \begin{equation}\label{eq:ito-rho}
        \begin{aligned}
            \mathrm{d}|\bm\rho|^2_H+2\|\bm\rho\|^2_V\,\mathrm{d}t & =-2\,{}_{V^\prime}\big\langle B(\bm u^{(1)})-B(\bm u^{(2)}), \bm\rho\big\rangle_V\,\mathrm{d}t -2\,r(\bm\rho,\bm\rho)\,\mathrm{d}t \\
            & \quad+\int_E\big|G(\bm u^{(1)})-G(\bm u^{(2)})\big|^2_H \,\nu(\mathrm{d}\xi)\,\mathrm{d}t+\mathrm{d}M_t,
        \end{aligned}
    \end{equation}
    where $M$ is a real-valued c\`adl\`ag local martingale with $M_0=0$. All identities and estimates are first localized.

    \textbf{Step 1: the drift terms.}
    By \eqref{est:B-diff_prod},
    \begin{equation}\label{est:rho-B}
        \big|{}_{V^\prime}\langle B(\bm u^{(1)})-B(\bm u^{(2)}), \bm\rho\rangle_V\big| \leq\frac{1}{2}\|\bm\rho\|^2_V+C\,\Psi(t)\,|\bm\rho|^2_H,
    \end{equation}
    with
    \begin{equation}
        \Psi(t):=\sum_{i=1}^2(\|\bm u^{(i)}(t)\|_V^4+\|\bm u^{(i)}(t)\|_V^2\|\bm u^{(i)}(t)\|_{\mathscr D(A)}^2),
    \end{equation}
    while \eqref{est:r} and Young's inequality give
    \begin{equation}
        |r(\bm\rho,\bm\rho)| \leq\frac{1}{4}\|\bm\rho\|^2_V+C|\bm\rho|^2_H.
    \end{equation}
    On $(0,\sigma_r]$, \eqref{asm:g_Lip-local} yields
    \begin{equation}\label{est:rho-G}
        \int_E\big|G(\bm u^{(1)}_-,\xi)-G(\bm u^{(2)}_-,\xi)\big|^2_H \,\nu(\mathrm{d}\xi)\leq C_r|\bm\rho_-|^2_H.
    \end{equation}
    Inserting \eqref{est:rho-B}--\eqref{est:rho-G} into \eqref{eq:ito-rho}, we obtain
    \begin{equation}\label{eq:rho-diff}
        \mathrm{d}|\bm\rho|^2_H+\frac{1}{4}\|\bm\rho\|^2_V\,\mathrm{d}t \leq C_r\Lambda(t)|\bm\rho|^2_H\,\mathrm{d}t+\mathrm{d}M_t, \qquad \Lambda:=1+\Psi ,
    \end{equation}
    the inequality holding on $[0,\sigma_r]$. Values at jump times in terms integrated against $\mathrm{d}t$ may be replaced by their left limits.

    \textbf{Step 2: the stopping time.}
    For $n\geq1$ set
    \begin{equation}\label{eq:tau-n}
        \tau_n:=\inf\Big\{t\in[0,\mathcal T]: \int_0^t\Lambda(s)\,\mathrm{d}s\geq n\Big\},
    \end{equation}
    with the convention $\inf\varnothing=\mathcal T$.  We claim that
    \begin{equation}\label{eq:tau-exhausts}
        \tau_n\to\mathcal T\quad\mathbb P\text{-almost surely as }n\to\infty .
    \end{equation}
    Since both $\bm u^{(1)}$ and $\bm u^{(2)}$ satisfy \eqref{eq:class-bound}, $\sup_{t\in[0,\mathcal T]}\|\bm u^{(i)}(t)\|_V$ and $\int_0^{\mathcal T}\|\bm u^{(i)}(t)\|^2_{\mathscr D(A)}\,\mathrm{d}t$ are finite almost surely for $i=1,2$; hence
    \begin{equation}
        \int_0^{\mathcal T}\Lambda(s)\,\mathrm{d}s \leq C\sum_{i=1}^2\Big(\mathcal T+\mathcal T\sup_{t\in[0,\mathcal T]} \|\bm u^{(i)}(t)\|^4_V +\sup_{t\in[0,\mathcal T]}\|\bm u^{(i)}(t)\|^2_V \int_0^{\mathcal T}\|\bm u^{(i)}(t)\|^2_{\mathscr D(A)}\,\mathrm{d}t\Big) <\infty
    \end{equation}
    almost surely, and \eqref{eq:tau-exhausts} follows. Since the paths are c\`adl\`ag in $V$, they are bounded on $[0,\mathcal T]$; hence $\sigma_r\uparrow\mathcal T$ almost surely. Write $\tau_{n,r}=\tau_n\wedge\sigma_r$.

    \textbf{Step 3: a stochastic Gr\"onwall argument.}
       Put $g=G(\bm u^{(1)}_-,\xi)-G(\bm u^{(2)}_-,\xi)$. For stopping times $a\leq b\leq\tau_{n,r}$, the noise terms, with the quadratic remainder written in uncompensated form, are $2\int_a^\cdot\int_E(\bm\rho_-,g)_H\,\tilde\eta$ and $\int_a^\cdot\int_E|g|_H^2\,\eta$. By BDG, Young and compensation, after localization,
    \begin{equation*}
        \begin{split}
            \mathbb E\sup_{a\leq t\leq b} \left|2\int_a^t\int_E(\bm\rho_-,g)_H\,\tilde\eta\right| & \leq\delta\mathbb E\sup_{a\leq t\leq b}|\bm\rho(t)|_H^2 +C_{r,\delta}\mathbb E\int_a^b|\bm\rho|_H^2\,\mathrm{d}s, \\
            \mathbb E\int_a^b\int_E|g|_H^2\,\eta & \leq C_r\mathbb E\int_a^b|\bm\rho|_H^2\,\mathrm{d}s.
        \end{split}
    \end{equation*}
    Here only the second-moment difference bound \eqref{asm:g_Lip-local} is used; the predictable indicator $\mathbbm1_{(0,\sigma_r]}$ preserves the terminal jump. Estimate the supremum and terminal dissipation separately, then absorb with small $\delta$ and remove the auxiliary localization by Fatou's lemma. Fix $n\geq1$ and $0\leq\tau_a\leq\tau_b\leq\tau_{n,r}$. This gives
    \begin{equation}\label{eq:rho-gronwall-input}
        \begin{aligned}
            \mathbb E\Big[\sup_{t\in[\tau_a,\tau_b]}|\bm\rho(t)|^2_H +\frac{1}{4}\int_{\tau_a}^{\tau_b}\|\bm\rho\|^2_V \,\mathrm{d}t\Big] & \leq C_r\mathbb E\Big[|\bm\rho(\tau_a)|^2_H +\int_{\tau_a}^{\tau_b}\Lambda(s)|\bm\rho(s)|^2_H \,\mathrm{d}s\Big] \\
            & \quad+\frac{1}{4}\,\mathbb E\sup_{t\in[\tau_a,\tau_b]} |\bm\rho(t)|^2_H.
        \end{aligned}
    \end{equation}
    After absorbing the last term, \eqref{eq:rho-gronwall-input} has the form \eqref{eq:gron-local} with
    \begin{equation*}
        X=|\bm\rho|^2_H,\quad Y=\frac{1}{4}\|\bm\rho\|^2_V,\quad R=C_r\Lambda,\qquad Z=0.
    \end{equation*}
    Here $\int_0^{\tau_{n,r}}\Lambda(s)\,\mathrm{d}s\leq n$ almost surely, and
    \begin{equation*}
        \mathbb E\int_0^{\tau_{n,r}}\Lambda(s)|\bm\rho(s)|^2_H\,\mathrm{d}s \leq n\,\mathbb E\sup_{t\in[0,\tau_{n,r}]} |\bm\rho(t)|^2_H<\infty ,
    \end{equation*}
    so that the hypotheses \eqref{eq:gron-hyp} are satisfied.  The stochastic Gr\"onwall inequality of Lemma \ref{lem:sto_gron} therefore applies on $[0,\tau_{n,r}]$ and yields
    \begin{equation}
        \mathbb E\Big[\sup_{t\in[0,\tau_{n,r}]}|\bm\rho(t)|^2_H +\frac{1}{4}\int_0^{\tau_{n,r}}\|\bm\rho\|^2_V\,\mathrm{d}t\Big] \leq C_{n,r}\mathbb E|\bm\rho(0)|^2_H=0 ,
    \end{equation}
    since $\bm u^{(1)}(0)=\bm u^{(2)}(0)$.  Thus $\bm\rho=0$ on $[0,\tau_{n,r}]$ almost surely for every $n$.  Letting first $n\to\infty$ and then $r\to\infty$, we conclude that $\bm\rho(t)=0$ for all $t\in[0,\mathcal T]$ almost surely, which is \eqref{eq:pathwise-uniqueness}.
\end{proof}

\begin{corollary}[Uniqueness of local solutions]\label{cor:unique-local}
    Let $(\bm u^{(1)},\tau_1)$ and $(\bm u^{(2)},\tau_2)$ be two local pathwise solutions in the sense of Definition \ref{def:pathwise-sol}, defined on the same stochastic basis, driven by the same compensated Poisson random measure and starting from the same initial datum.  Then
    \begin{equation}
        \mathbb P\big(\bm u^{(1)}(t)=\bm u^{(2)}(t)\ \text{for all } t\in[0,\tau_1\wedge\tau_2]\big)=1 .
    \end{equation}
\end{corollary}

\begin{proof}
    Apply the proof of Theorem \ref{thm:unique} on the common interval $[0,\tau_1\wedge\tau_2]$, with the norm exit time \eqref{eq:uniqueness-norm-stop} restricted to this interval and with $\tau_{n,r}$ replaced by $\tau_{n,r}\wedge\tau_1\wedge\tau_2$. The local bound \eqref{eq:local-bound} gives finiteness of the integrated drift coefficient, and the stopped $V$-valued c\`adl\`ag paths are bounded. Letting $n\to\infty$ and then $r\to\infty$ proves the claim, including the common terminal jump.
\end{proof}

\subsection{Existence of pathwise solutions}\label{subsec:gk}

We pass from martingale solutions to pathwise solutions using the following criterion of Gyöngy and Krylov \cite[Lemma~1.1]{GyongyKrylov1996}.

\begin{lemma}[Gyöngy--Krylov]\label{lem:GK}
    Let $E$ be a Polish space and let $\{Y_n\}_{n\geq1}$ be a tight sequence of $E$-valued random variables.  If, for every pair of subsequences $n_k,m_k\to\infty$, every limit point of the sequence of laws $\{\mathcal L(Y_{n_k},Y_{m_k})\}_{k\geq1}$ is concentrated on the diagonal $\Delta_E:=\{(y,y):y\in E\}$ of $E\times E$, then there exists an $E$-valued random variable $Y$, defined on the underlying probability space, such that $Y_n\to Y$ in probability as $n\to\infty$; in particular the set of limit points of $\{\mathcal L(Y_n)\}_{n\geq1}$ is the singleton $\{\mathcal L(Y)\}$.
\end{lemma}

\begin{lemma}\label{lem:limit-passage}
    Assume Assumptions \ref{asm:martingale} and \ref{asm:A-est}, with $\bm u_0\in L^4(\Omega,\mathcal F_0;V)$. Let $n_k,m_k\to\infty$. Apply Lemma \ref{lem:s-j} to a convergent subsequence of the joint laws of $(\bm u_0,\bm u_{n_k},\bm u_{m_k},\eta)$ on $V\times\mathcal Z_{\mathcal T}^2\times\mathcal N$, and write the resulting convergence as
    \begin{equation}\label{eq:pair-convergence}
        (a_k,u_k^1,u_k^2,\eta_k)\to(a,u^1,u^2,\eta^*) \quad\text{almost surely}.
    \end{equation}
    Each tuple has the corresponding original joint law. With respect to the joint filtration \eqref{eq:joint-limit-filtration} with $d=2$, $\eta^*$ is a Poisson random measure with compensator $\mathrm{d}t\otimes\nu$. Both $u^1$ and $u^2$ belong to $\mathfrak C([0,\mathcal T])$ and solve the regularised equation, with initial datum $a$ and the same compensated measure $\widetilde\eta^*=\eta^*-\mathrm{d}t\otimes\nu$.
\end{lemma}
\begin{proof}
    Apply Lemma \ref{lem:retained-noise} with $d=2$, $r_k^1=n_k$ and $r_k^2=m_k$. It yields both equations with the same Poisson driver relative to the joint filtration, as well as the stated energy bounds and path regularity.
\end{proof}

We may now state the existence result.

\begin{theorem}[Existence of pathwise solutions]\label{thm:pathwise-mod}
    Let Assumptions \ref{asm:martingale} and \ref{asm:A-est} hold and let $\bm u_0\in L^4(\Omega,\mathcal F_0;V)$.  Then there exists a process $\bm u\in\mathfrak C([0,\mathcal T])$, defined on the original stochastic basis and driven by the original Poisson random measure $\eta$, which solves the regularised equation; its stopped version $\bm u(\cdot\wedge\tau)$ solves \eqref{eqn:abstract} on $[0,\tau]$.  Moreover, $\bm u$ is unique among the solutions of the regularised equation in $\mathfrak C([0,\mathcal T])$, by the same argument as in Theorem \ref{thm:unique}; see Remark \ref{rk:unique-regularised}.
\end{theorem}

\begin{proof}
    Let $\mathbb G=\{\mathcal G_t\}_{t\geq0}$ be the usual augmentation of the filtration generated by $\bm u_0$ and $\eta$. Then $\mathcal G_t\subset\mathcal F_t$ for every $t\geq0$, and $\eta$ is also a Poisson random measure relative to $\mathbb G$. Write $\mathcal I=\sigma(\sigma(\bm u_0,\eta)\cup\mathcal N_{\mathbb P})$, where $\mathcal N_{\mathbb P}$ consists of the $\mathbb P$-null sets in $\mathcal F$. In particular, $\mathcal G_t\subset\mathcal I$ for all $t\geq0$. The finite-dimensional construction on this basis gives $\mathbb G$-adapted Galerkin solutions, which coincide with $\bm u_n$ by finite-dimensional pathwise uniqueness. We carry out the limiting construction in this filtration.

    The laws of the Galerkin solutions are tight on $\mathcal Z_{\mathcal T}$ by Lemma \ref{lem:tight-u_n}, and on the Polish space $\mathcal X=L^2(0,\mathcal T;V)$ by the compact embedding
    \begin{equation*}
        L^2(0,\mathcal T;\mathscr D(A))\cap W^{\alpha,2}(0,\mathcal T;H) \hookrightarrow\hookrightarrow L^2(0,\mathcal T;V),\qquad 0<\alpha<\tfrac12.
    \end{equation*}
    Fix a limit point $\mu$ of the laws of any pair of subsequences in $\mathcal X^2$. Adjoin $\bm u_0$ and $\eta$ and extract further in $V\times\mathcal Z_{\mathcal T}^2\times\mathcal N$. Lemma \ref{lem:limit-passage} gives two limits with pair law $\mu$, the same initial datum $a$, and the same driving measure $\eta^*$ on the same filtered space. Remark \ref{rk:unique-regularised} implies $u^1=u^2$ almost surely. Thus every such $\mu$ is supported on the diagonal of $\mathcal X^2$. Lemma \ref{lem:GK} yields $\bm u_n\to\bm u$ in probability in $\mathcal X$ on the original basis.

    Each $\bm u_n$ is $\mathcal I$-measurable as an $\mathcal X$-valued random variable. Passing to an almost surely convergent subsequence shows that the $\mathcal X$-valued limit $\bm u$ has the same measurability.

    It remains to identify the equation on that basis. The uniform fourth $V$-moment in \eqref{est:uniform-DA-p} and uniform integrability upgrade this convergence to
    \begin{equation*}
        \mathbb E\int_0^{\mathcal T}\|\bm u_n-\bm u\|_V^2\,\mathrm{d}s\to0.
    \end{equation*}
    Since $\bm u_n(\cdot-)$ are $\mathbb G$-predictable, the limit has a $\mathbb G$-predictable representative in this space. Moreover, lower semicontinuity of the energy bounds, or the joint representations above, gives $\bm u\in L^\infty(0,\mathcal T;V)\cap L^2(0,\mathcal T;\mathscr D(A))$ almost surely, with the inherited moment bounds. Lemma \ref{lem:local-noise-convergence} and the projection contraction give
    \begin{equation}\label{eq:original-noise-integrand-convergence}
        \mathbb E\int_0^{\mathcal T}\int_E |P_nG(\bm u_n(s-),\xi)-G(\bm u(s),\xi)|_H^2 \,\nu(\mathrm{d}\xi)\,\mathrm{d}s\to0.
    \end{equation}
    Indeed, split the difference into $P_n[G(\bm u_n)-G(\bm u)]+(P_n-I)G(\bm u)$ and apply Lemma \ref{lem:local-noise-convergence}, using the inherited fourth $V$-moments, followed by dominated convergence for the projection error. Doob's inequality and the Poisson isometry imply
    \begin{equation}\label{eq:original-noise-integral-convergence}
        \mathbb E\sup_{t\leq\mathcal T} \left|\int_0^t\int_E [P_nG(\bm u_n(s-),\xi)-G(\bm u(s),\xi)]\, \tilde\eta(\mathrm{d}s,\mathrm{d}\xi)\right|_H^2 \to0.
    \end{equation}
    For each $\psi\in V$, the drift integrals converge uniformly in time in probability by the estimates in \eqref{eq:convergence-B}, using the convergence in $\mathcal X$, the convergence of the backgrounds, and the uniform energy bounds. The linear terms follow from $A:V\to V'$ and \eqref{est:r}. Passing to a subsequence for a countable dense family of tests identifies the time-integrated equation. Its right-hand side provides a $\mathbb G$-adapted c\`adl\`ag representative, equal to the $L^2$ limit almost everywhere in time, with initial value $\bm u_0$. Lemma \ref{lem:strong-cadlag} gives its $V$-valued c\`adl\`ag version. For this version the predictable integrand equals $G(\bm u(s-),\xi)$ almost everywhere in $\mathrm{d}\mathbb P\,\mathrm{d}s\,\mathrm{d}\nu$. Thus the regularised equation holds with the original noise $\eta$. The inherited bounds give membership in $\mathfrak C([0,\mathcal T])$.

    Since the Borel $\sigma$-field of $D([0,\mathcal T];V)$ is generated by coordinate evaluations, $\mathbb G$-adaptedness makes this representative $\mathcal I$-measurable as a $D([0,\mathcal T];V)$-valued random variable. The continuous injection $L^2(0,\mathcal T;\mathscr D(A))\to\mathcal X$ has a Borel inverse on its image, so the same holds with values in $L^2(0,\mathcal T;\mathscr D(A))$. The path therefore admits a $\sigma(\bm u_0,\eta)$-measurable version in their intersection. The exit time \eqref{def:tau-exit} is a $\mathbb G$-stopping time and likewise admits a $\sigma(\bm u_0,\eta)$-measurable version. Stopping at this exit yields the local pathwise solution of \eqref{eqn:abstract}; uniqueness follows from Remark \ref{rk:unique-regularised} and Corollary \ref{cor:unique-local}.
\end{proof}

\begin{proposition}[Square-integrable initial data]\label{prop:initial-localization}
    Let Assumptions \ref{asm:martingale} and \ref{asm:A-est} hold and let $\bm u_0\in L^2(\Omega,\mathcal F_0;V)$. For every $\mathcal T>0$, there is a local pathwise solution $(\bm u,\tau)$ of \eqref{eqn:abstract}, with $0<\tau\leq\mathcal T$ almost surely.
\end{proposition}
\begin{proof}
    We localize the initial datum as in \cite[Proposition~3.6]{BrzezniakSlavik2021}. Put $\Omega_j=\{j-1\leq\|\bm u_0\|_V<j\}$ and $\bm u_0^j=\mathbbm{1}_{\Omega_j}\bm u_0$, $j\geq1$. Theorem \ref{thm:pathwise-mod} gives a regularised solution $\bm w^j\in\mathfrak C([0,\mathcal T])$ with initial datum $\bm u_0^j$ and driver $\eta$. Let $\rho_j$ be its exit time \eqref{def:tau-exit}. Set
    \begin{equation*}
        Y_j(t)=\sup_{s\leq t\wedge\rho_j}\|\bm w^j(s)\|_V^2 +\int_0^{t\wedge\rho_j}|A\bm w^j|_H^2\,\mathrm{d}s.
    \end{equation*}
    Right continuity at zero and the integrable bound $Y_j(\mathcal T)$ give $Y_j(t)\to\|\bm u_0^j\|_V^2$ almost surely and in $L^1$ as $t\downarrow0$. Choose deterministic $t_j\in(0,\mathcal T]$ such that
    \begin{equation*}
        \mathbb E[\mathbbm{1}_{\Omega_j}Y_j(t_j)] \leq\mathbb E[\mathbbm{1}_{\Omega_j}\|\bm u_0\|_V^2]+2^{-j}.
    \end{equation*}
    Since $\Omega_j\in\mathcal F_0$, the pasting
    \begin{equation*}
        \bm u=\sum_{j\geq1}\mathbbm{1}_{\Omega_j}\bm w^j, \qquad \tau=\sum_{j\geq1}\mathbbm{1}_{\Omega_j}(t_j\wedge\rho_j)
    \end{equation*}
    is adapted, $\tau$ is a strictly positive stopping time, and the original equation holds on $[0,\tau]$. Moreover,
    \begin{equation*}
        \mathbb E\left[\sup_{s\leq\tau}\|\bm u(s)\|_V^2 +\int_0^\tau|A\bm u|_H^2\,\mathrm{d}s\right] \leq\mathbb E\|\bm u_0\|_V^2+1.
    \end{equation*}
    The processes are measurable functions of $(\bm u_0,\eta)$, and uniqueness follows from Corollary \ref{cor:unique-local}.
\end{proof}

Local solutions with the same initial datum and noise agree on their common intervals by Corollary \ref{cor:unique-local}. We use this consistency to construct the maximal solution.

\begin{lemma}[Essential supremum of local existence times]\label{lem:ess-sup}
    Let $\mathcal S$ be a nonempty family of stopping times such that $\tau\vee\hat\tau\in\mathcal S$ whenever $\tau,\hat\tau\in\mathcal S$. Then there exist a stopping time $\xi$ and an increasing sequence $\{\tau_N\}_{N\geq1}\subset\mathcal S$ with $\tau_N\uparrow\xi$ almost surely, and $\xi$ dominates almost surely every element of $\mathcal S$.
\end{lemma}

\begin{proof}
    Let $\mathscr L$ be the set of countable subfamilies of $\mathcal S$ and let $c:=\sup_{D\in\mathscr L}\mathbb E[1-e^{-\sup_{\tau\in D}\tau}]$. Choose $D_n\in\mathscr L$ with $\mathbb E[1-e^{-\sup D_n}]\geq c-1/n$ and replace $D_n$ by $D_1\cup\cdots\cup D_n$, so that $\sup D_n$ increases to $\sup D$, where $D:=\bigcup_nD_n\in\mathscr L$; by monotone convergence, $\mathbb E[1-e^{-\sup D}]=c$.  For $\xi:=\sup_{\tau\in D}\tau$ and any $\tau\in\mathcal S$ the family $D\cup\{\tau\}$ belongs to $\mathscr L$, whence $\mathbb E[1-e^{-(\xi\vee\tau)}]\leq c=\mathbb E[1-e^{-\xi}]$, which is possible only if $\xi\vee\tau=\xi$ almost surely.  Finally, enumerating $D=\{\sigma_1,\sigma_2,\dots\}$ and setting $\tau_N:=\sigma_1\vee\cdots\vee\sigma_N$, the lattice property gives $\tau_N\in\mathcal S$ and $\tau_N\uparrow\xi$ almost surely; in particular $\xi$ is a stopping time, since $\{\xi\leq t\}=\bigcap_N\{\tau_N\leq t\}\in\mathcal F_t$ for every $t\geq0$.
\end{proof}

\begin{theorem}[Maximal pathwise solution]\label{thm:wp-maximal}
    Let Assumptions \ref{asm:martingale} and \ref{asm:A-est} hold and let $\bm u_0\in L^2(\Omega,\mathcal F_0;V)$.  Then \eqref{eqn:abstract} possesses a unique maximal pathwise solution $(\bm u,\xi)$, and
    \begin{equation}\label{eq:blow-up}
        \sup_{t\in[0,\xi)}\|\bm u(t)\|_V +\int_0^\xi\|\bm u\|^2_{\mathscr D(A)}\,\mathrm{d}t=\infty \quad\mathbb P\text{-almost surely on }\{\xi<\infty\}.
    \end{equation}
\end{theorem}

\begin{proof}
    Let $\mathbb G:=\{\mathcal G_t\}_{t\geq0}$ be the usual augmentation of the filtration generated by $\bm u_0$ and $\eta$. It satisfies $\mathcal G_t\subset\mathcal F_t$ for every $t\geq0$, and $\eta$ is a Poisson random measure relative to $\mathbb G$. We first work on this basis; all constructed processes and stopping times are then measurable functions of $(\bm u_0,\eta)$, up to completion.

    \emph{Construction of a consistent family.} Let $\mathcal S$ be the set of lifetimes of local pathwise solutions on the $\mathbb G$-basis. Proposition \ref{prop:initial-localization}, applied on this basis and with a finite time horizon, shows that $\mathcal S$ is nonempty. Its initial exit time is strictly positive by right continuity in $V$. If $(\bm u^1,\tau_1)$ and $(\bm u^2,\tau_2)$ are local solutions, they coincide on $[0,\tau_1\wedge\tau_2]$ by Corollary \ref{cor:unique-local}. Following the first solution up to $\tau_1$ and the second afterwards, when $\tau_2>\tau_1$, gives a local solution up to $\tau_1\vee\tau_2$. This is an adapted pasting at stopping times, and its bound \eqref{eq:local-bound} is at most the sum of the two original bounds. Thus $\mathcal S$ is closed under finite maxima.

    Lemma \ref{lem:ess-sup} yields $\tau_N\in\mathcal S$ with $\tau_N\uparrow\xi$, where $\xi$ dominates every element of $\mathcal S$. Let $\bm u^{(N)}$ be the corresponding solutions. Pathwise uniqueness makes the prescription
    \begin{equation}\label{eq:maximal-u}
        \bm u(t)=\bm u^{(N)}(t),\qquad 0\leq t\leq\tau_N,
    \end{equation}
    consistent on $\bigcup_N[0,\tau_N]$, which contains $[0,\xi)$.

    \emph{A closed stopped interval.} Fix $R,T>0$ and put
    \begin{equation}\label{eq:rho-R}
        \begin{split}
            K_N(t) & =\sup_{s\leq t}\|\bm u^{(N)}(s)\|_V +\int_0^t|A\bm u^{(N)}(s)|_H^2\,\mathrm{d}s, \qquad t\leq\tau_N, \\
            \sigma_N & =\tau_N\wedge T\wedge \inf\{t\in[0,\tau_N]:K_N(t)\geq R\}, \qquad \rho_{R,T}=\lim_N\sigma_N.
        \end{split}
    \end{equation}
    Here $\inf\varnothing=\infty$. These are $\mathbb G$-stopping times, and consistency implies $\sigma_N\leq\sigma_{N+1}$. The stochastic interval $D_{R,T}:=\bigcup_N\{(\omega,s):0<s\leq\sigma_N(\omega)\}$ is predictable. On this interval define
    \begin{equation*}
        f_R=-\mathbbm{1}_{D_{R,T}}[B(\bm u)+R(\bm u)],\qquad h_R(s,z)=\mathbbm{1}_{D_{R,T}}G(\bm u(s-),z),
    \end{equation*}
    and set both coefficients to zero elsewhere. The coefficients are well defined by consistency; $f_R$ is progressively measurable and $h_R$ is predictable. The threshold bounds hold before each $\sigma_N$, without imposing a bound on its terminal jump. Hence \eqref{est:B-H-symmetric}, \eqref{est:r} and \eqref{asm:g_moment-V}, with interpolation between its second and fourth moments, give
    \begin{equation}\label{eq:extension-data-bounds}
        \mathbb E\int_0^T|f_R|_H^2\,\mathrm{d}s\leq C_{R,T}, \qquad \int_0^T\int_E\|h_R(s,z)\|_V^p\,\nu(\mathrm{d}z)\,\mathrm{d}s \leq C_{p,R,T}\quad\text{a.s.},\qquad 2\leq p\leq4.
    \end{equation}
    Indeed, $\int_{D_{R,T}}|A\bm u|_H^2\,\mathrm{d}s\leq R$ and $\|\bm u\|_V\leq R$ there for almost every time, while $\|\bm u(s-)\|_V\leq R$ on $D_{R,T}$.

    The linear equation
    \begin{equation}\label{eq:stopped-linear-extension}
        \mathrm{d}Z_R+AZ_R\,\mathrm{d}t=f_R\,\mathrm{d}t +\int_E h_R(t,z)\,\widetilde\eta(\mathrm{d}t,\mathrm{d}z), \qquad Z_R(0)=\bm u_0,
    \end{equation}
    has a unique adapted solution with
    \begin{equation*}
        Z_R\in D([0,T];V)\cap L^2(0,T;\mathscr D(A)),\qquad \mathbb E\left[\sup_{t\leq T}\|Z_R(t)\|_V^2 +\int_0^T|AZ_R|_H^2\,\mathrm{d}t\right]<\infty.
    \end{equation*}
    This follows from the linear spectral estimate in the proof of Lemma \ref{lem:strong-cadlag}. Linear uniqueness gives $Z_R=\bm u$ on every $[0,\sigma_N]$ and hence on $[0,\rho_{R,T})$.

    Write $\rho=\rho_{R,T}$. The preceding identification and the strong left limit give $\|Z_R(\rho-)\|_V\leq R$ on $\{\rho>0\}$. Thus $\mathbbm{1}_{\{0<s\leq\rho\}}G(Z_R(s-),z)$ is predictable and square integrable in $V$. The predictable sets $D_{R,T}$ and $\{0<s\leq\rho\}$ differ only on the graph of $\rho$, a set of zero $\mathrm{d}\mathbb P\,\mathrm{d}s$ measure. Consequently
    \begin{equation}\label{eq:closed-noise-identification}
        \mathbb E\int_0^T\int_E \big\|h_R(s,z)-\mathbbm{1}_{\{0<s\leq\rho\}}G(Z_R(s-),z)\big\|_V^2 \,\nu(\mathrm{d}z)\,\mathrm{d}s=0.
    \end{equation}
    The corresponding Poisson integrals are indistinguishable by the isometry. The drift coefficients likewise agree almost everywhere up to $\rho$. Thus $Z_R(\cdot\wedge\rho)$ satisfies the original equation on the closed interval $[0,\rho]$, including its endpoint jump. In particular, with $M_R=\int h_R\,\mathrm{d}\widetilde\eta$,
    \begin{equation}\label{eq:u-at-xi}
        Z_R(\rho)=Z_R(\rho-)+\Delta M_R(\rho) \quad\text{on }\{\rho>0\}.
    \end{equation}
    The linear energy bound gives
    \begin{equation}\label{eq:restart-moment}
        \mathbb E\|Z_R(\rho)\|_V^2 \leq\mathbb E\sup_{t\leq T}\|Z_R(t)\|_V^2<\infty.
    \end{equation}

    \emph{Continuation and the blow-up criterion.} On the shifted filtration $\{\mathcal G_{\rho+t}\}_{t\geq0}$,
    \begin{equation*}
        \eta_\rho((0,t]\times A)=\eta((\rho,\rho+t]\times A)
    \end{equation*}
    is a Poisson random measure with compensator $\mathrm{d}t\otimes\nu$; the stopping-time version of independent increments follows by approximating the bounded time $\rho$ from above by simple stopping times. Apply Proposition \ref{prop:initial-localization} on this basis with initial value $Z_R(\rho)$ and time horizon $1$. By \eqref{eq:restart-moment}, it gives a local solution $(\bm v,\delta)$ with $0<\delta\leq1$ almost surely. Define
    \begin{equation}\label{eq:glued}
        \bm z(t)= \begin{cases}
            Z_R(t), & 0\leq t\leq\rho, \\
            \bm v(t-\rho), & \rho<t\leq\rho+\delta.
        \end{cases} \qquad \varpi=\rho+\delta.
    \end{equation}
    The time $\varpi$ is a $\mathbb G$-stopping time. The states agree at the joining time, and the shifted noise uses only $(\rho,\rho+t]$, so that the jump at $\rho$ is counted exactly once. The linear bound and the bound for $\bm v$ give \eqref{eq:local-bound}. Therefore $\varpi\in\mathcal S$, and
    \begin{equation}\label{eq:strict-energy-localisation}
        \rho_{R,T}<\varpi\leq\xi\qquad\text{almost surely}.
    \end{equation}
    On the event
    \begin{equation*}
        \left\{\xi\leq T,\quad \sup_{t<\xi}\|\bm u(t)\|_V +\int_0^\xi|A\bm u|_H^2\,\mathrm{d}t<R\right\},
    \end{equation*}
    all thresholds before $\xi$ are unhit and $\lim_N\sigma_N=\xi$, whether or not any $\tau_N$ attains $\xi$. This contradicts \eqref{eq:strict-energy-localisation}. Taking a countable union over integers $R,T\geq1$ proves \eqref{eq:blow-up}.

    \emph{Maximality and uniqueness on the original basis.} Let $(\hat{\bm u},\hat\tau)$ be any local pathwise solution on the original stochastic basis. Uniqueness gives $\hat{\bm u}=\bm u$ on $[0,\hat\tau\wedge\tau_N]$ for every $N$. On $\{\xi<\infty,\hat\tau\geq\xi\}$, its finite local energy bound would contradict \eqref{eq:blow-up}. Hence $\hat\tau\leq\xi$ almost surely. Together with \eqref{eq:maximal-u}, this proves both clauses of Definition \ref{def:pathwise-sol}. For two maximal solutions, applying this domination to their localizing sequences makes their lifetimes equal; local uniqueness then makes their paths equal on $[0,\xi)$. This completes the proof.
\end{proof}

\subsection{Global well-posedness}

Theorem \ref{thm:wp-maximal} provides a maximal pathwise solution $(\bm u,\xi)$ and the blow-up criterion \eqref{eq:blow-up}. We exclude finite-time blow-up by localized energy estimates.

The energy estimates use the barotropic and baroclinic decomposition, followed by estimates for the vertical derivative, the horizontal gradient and the temperature.

Throughout this subsection we write
\begin{equation}\label{eq:tau-N-global}
    \tau_N:=\Big(\inf\Big\{t\in[0,\xi): \|\bm u(t)\|_V+\int_0^t\|\bm u(s)\|^2_{\mathscr D(A)}\,\mathrm{d}s \geq N\Big\}\Big)\wedge\mathcal T,\quad N\geq1 .
\end{equation}
Then $\tau_N\uparrow\xi\wedge\mathcal T$, and, on $\{\tau_N>0\}$,
\begin{equation}\label{eq:tau-N-bound}
    \sup_{t\in[0,\tau_N)}\|\bm u(t)\|_V +\int_0^{\tau_N}\|\bm u\|^2_{\mathscr D(A)}\,\mathrm{d}t\leq 2N \quad\mathbb P\text{-almost surely}.
\end{equation}
The closed stopped solutions have the bound \eqref{eq:local-bound}: indeed, $\tau_N\leq\rho_{2N+1,\mathcal T}$, and the construction \eqref{eq:stopped-linear-extension}--\eqref{eq:restart-moment} controls the endpoint jump as well as the energy before it. For the $L^6$ It\^o formula and the baroclinic energy estimate, we impose the following growth conditions; cf. \cite[(2.27)]{BrzezniakSlavik2021}.

\begin{assumption}\label{asm:global}
    In addition to Assumptions \ref{asm:martingale} and \ref{asm:A-est}, assume that \eqref{asm:g_moment-V} also holds for $p=3$ and that, for $p=1,3$ and $\bm u=(\bm v,T)\in V$,
    \begin{equation}\label{asm:g_L6}
        \int_E|\widetilde{G_1}(\bm u,\xi)|_6^{2p}\,\nu(\mathrm{d}\xi) \leq C_p(1+|\widetilde{\bm v}|_6^{2p}).
    \end{equation}
\end{assumption}

\subsubsection{The Itô formula on the stopped interval}

\begin{lemma}\label{lem:ito-stopped}
    Let $(\bm u,\tau)$ be a local pathwise solution under Assumption \ref{asm:global}. On $[0,\tau]$ it has the regularity
    \begin{equation}\label{eq:stopped-regularity}
        \bm u\in D([0,\tau];V)\cap L^2(0,\tau;\mathscr D(A)), \qquad A\bm u+B(\bm u)+R(\bm u)\in L^2(0,\tau;H) \quad\text{almost surely}.
    \end{equation}
    For the functional
    \begin{equation}\label{eq:psi-L6}
        \psi(\bm u)=|\widetilde{\bm v}|_6^6,
    \end{equation}
    the following identity holds on $[0,\tau]$:
    \begin{align*}
        \psi(\bm u(t))={} & \psi(\bm u_0) -\int_0^tD\psi(\bm u)[A\bm u+B(\bm u)+R(\bm u)]\,\mathrm{d}s \\
        & +\int_0^t\int_E D\psi(\bm u(s-))[G(\bm u(s-),\xi)]\,\tilde\eta(\mathrm{d}s,\mathrm{d}\xi) \\
        & +\int_0^t\int_E\big\{\psi(\bm u(s-)+G)-\psi(\bm u(s-)) -D\psi(\bm u(s-))[G]\big\}\,\eta(\mathrm{d}s,\mathrm{d}\xi).
    \end{align*}
\end{lemma}
\begin{proof}
    The first assertion follows from Lemma \ref{lem:strong-cadlag}, \eqref{est:B-H-symmetric} and the local energy bound. For $\bm h=(\bm h_1,h_2)\in V$,
    \begin{equation*}
        D\psi(\bm u)[\bm h] =6\int_{\mathcal O}|\widetilde{\bm v}|^4\widetilde{\bm v} \cdot\widetilde{\bm h_1}\,\mathrm{d}x\,\mathrm{d}y\,\mathrm{d}z.
    \end{equation*}
    The interpolation bound \eqref{est:GN-d3-10} implies $\|D\psi(\bm u)\|_{H'}\leq C\|\bm u\|_V^4\|\bm u\|_{\mathscr D(A)}$. Thus the drift pairing is integrable. Moreover,
    \begin{equation*}
        |\psi(\bm u+\bm g)-\psi(\bm u)-D\psi(\bm u)[\bm g]| \leq C\big(\|\bm u\|_V^4\|\bm g\|_V^2+\|\bm g\|_V^6\big).
    \end{equation*}
    Write $F=A\bm u+B(\bm u)+R(\bm u)$ and set
    \begin{equation*}
        \sigma_m=\tau\wedge m\wedge\inf\left\{t\leq\tau: \|\bm u(t)\|_V>m\ \text{or}\ \int_0^t(\|\bm u\|_{\mathscr D(A)}^2+|F|_H^2)\,\mathrm{d}s>m\right\}.
    \end{equation*}
    Work first on the $\mathcal F_0$-measurable event $\{\|\bm u_0\|_V\leq m\}$ and on $[0,\sigma_m]$. The pre-jump values have $V$ norm at most $m$. Thus
    \begin{equation*}
        \int_0^{\sigma_m}\!\int_E |D\psi(P_n\bm u_-)[P_nG]|^2\,\mathrm{d}\nu\,\mathrm{d}s \leq C_m\int_0^{\sigma_m}\!\int_E\|G\|_V^2\,\mathrm{d}\nu\,\mathrm{d}s \leq C_m,
    \end{equation*}
    and the compensator of the jump remainder is bounded by $C_m\int_0^{\sigma_m}\int_E(\|G\|_V^2+\|G\|_V^6)\,\mathrm{d}\nu\,\mathrm{d}s \leq C_m$. Apply the finite-dimensional jump It\^o formula to $\psi(P_n\bm u)$ on these stopped intervals. Spectral projections converge uniformly on the compact closure of each c\`adl\`ag $V$-path and in $L^2(0,\tau;\mathscr D(A))$. The interpolation estimate gives convergence of the drift pairings in $L^1(0,\tau)$. After localization, the isometry gives convergence of the linear stochastic integral using the second moment bound, while compensation and dominated convergence give $L^1$ convergence of the jump remainder using the second and sixth moment bounds in Assumption \ref{asm:global}. Passing to the limit in $n$ and then letting $m\to\infty$ proves the identity, since the events and stopping times exhaust the local existence interval by \eqref{eq:stopped-regularity}.
\end{proof}

\subsubsection{The baroclinic \texorpdfstring{$L^6$}{L6} estimate}

By Lemma \ref{lem:ito-stopped}, we may apply the jump It\^o formula to $|\widetilde{\bm v}|_6^6$ along \eqref{eq:baroclinic}. For an $L^p(\mathbb R^d)$ jump formula, see \cite[Theorem~2.1]{7984367}.

Pairings with $|\bm w|^4\bm w$ below denote the duality between $L^6(\mathcal O)$ and $L^{6/5}(\mathcal O)$. In integrals against $\eta$ or $\tilde\eta$, unmarked solution values denote their left limits. The transport terms cancel by integration by parts.

\begin{equation}\label{est:baroclinic-0-1}
    -\int_{\mathcal{O}}( (\widetilde{\bm v}\cdot\nabla)\widetilde{\bm v}+\Phi(\widetilde{\bm v})\partial_z\widetilde{\bm v})\cdot|\widetilde{\bm v}|^4\widetilde{\bm v}\,\mathrm{d}x\,\mathrm{d}y\,\mathrm{d}z=0.
\end{equation}
The Coriolis term vanishes:
\begin{equation}
    \langle -f\bm k\times\widetilde{\bm v},|\widetilde{\bm v}|^4\widetilde{\bm v}\rangle=0.
\end{equation}
By \eqref{eq:div-v-barotropic} and \eqref{bd_cond-barotropic}, the barotropic transport term also cancels:
\begin{equation}\label{est:baroclinic-0-3}
    \int_{\mathcal{O}}( (\overline{\bm v}\cdot\nabla)\widetilde{\bm v})\cdot|\widetilde{\bm v}|^4\widetilde{\bm v}\,\mathrm{d}x\,\mathrm{d}y\,\mathrm{d}z=0.
\end{equation}

Combining these cancellations with the jump It\^o formula gives
\begin{equation}\label{est:L6-ito}
    \begin{aligned}
        |\widetilde{\bm v}|_6^6 & +6\int_0^t\int_{\mathcal{O}}|\nabla\widetilde{\bm v}|^2|\widetilde{\bm v}|^4+\left|\nabla |\widetilde{\bm v}|^2\right|^2\,|\widetilde{\bm v}|^2\,\mathrm{d}x\,\mathrm{d}y\,\mathrm{d}z\,\mathrm{d}s \\
        & +6\int_0^t\int_{\mathcal{O}}|\partial_z\widetilde{\bm v}|^2|\widetilde{\bm v}|^4+\left|\partial_z |\widetilde{\bm v}|^2\right|^2\,|\widetilde{\bm v}|^2\,\mathrm{d}x\,\mathrm{d}y\,\mathrm{d}z\,\mathrm{d}s \\
        = & |\widetilde{\bm v}_0|_6^6-6\int_0^t\int_{\mathcal{O}}( (\widetilde{\bm v}\cdot\nabla)\overline{\bm v}- \overline{(\widetilde{\bm v}\cdot\nabla)\widetilde{\bm v}+(\nabla\cdot\widetilde{\bm v})\widetilde{\bm v}})\cdot|\widetilde{\bm v}|^4\widetilde{\bm v}\,\mathrm{d}x\,\mathrm{d}y\,\mathrm{d}z\,\mathrm{d}s \\
        & +6\int_0^t\int_{\mathcal{O}}\left( \int_{-h}^z\nabla T\,\mathrm{d}z^\prime-\frac{1}{h}\int_{-h}^0\int_{-h}^z\nabla T\,\mathrm{d}z^\prime\,\mathrm{d}z\right)\cdot|\widetilde{\bm v}|^4\widetilde{\bm v}\,\mathrm{d}x\,\mathrm{d}y\,\mathrm{d}z\,\mathrm{d}s \\
        & +\int_0^t\int_E\left(|\widetilde{\bm v}+\widetilde{G_1}(\bm v(s-),T(s-),\xi)|_6^6-|\widetilde{\bm v}|_6^6-6\langle \widetilde{G_1}(\bm v(s-),T(s-),\xi),|\widetilde{\bm v}|^4\widetilde{\bm v}\rangle\right)\eta(\mathrm{d}s,\mathrm{d}\xi) \\
        & +6\int_0^t\int_E\langle \widetilde{G_1}(\bm v(s-),T(s-),\xi),|\widetilde{\bm v}|^4\widetilde{\bm v}\rangle\tilde{\eta}(\mathrm{d}s,\mathrm{d}\xi).
    \end{aligned}
\end{equation}
Integrating by parts and using the boundary conditions \eqref{bd:baroclinic-3} and \eqref{bd_cond-barotropic}, the two deterministic terms on the right-hand side of \eqref{est:L6-ito} can be regrouped as
\begin{equation}
    \begin{aligned}
        & -\int_{\mathcal{O}}((\widetilde{\bm v}\cdot\nabla)\overline{\bm v} -\overline{(\widetilde{\bm v}\cdot\nabla)\widetilde{\bm v} +(\nabla\cdot\widetilde{\bm v})\widetilde{\bm v}})\cdot |\widetilde{\bm v}|^4\widetilde{\bm v}\,\mathrm{d}x\,\mathrm{d}y\,\mathrm{d}z \\
        & \qquad+\int_{\mathcal{O}}\left( \int_{-h}^z\nabla T\,\mathrm{d}z^\prime -\frac{1}{h}\int_{-h}^0\int_{-h}^z\nabla T\,\mathrm{d}z^\prime\,\mathrm{d}z \right)\cdot|\widetilde{\bm v}|^4\widetilde{\bm v}\,\mathrm{d}x\,\mathrm{d}y\,\mathrm{d}z =I_1+I_2+I_3,
    \end{aligned}
\end{equation}
where
\begin{equation}
    \begin{aligned}
        I_1 & :=\int_{\mathcal{O}}\left((\mathrm{div}\,\widetilde{\bm v}) \overline{\bm v}\cdot|\widetilde{\bm v}|^4\widetilde{\bm v} +(\widetilde{\bm v}\cdot\nabla)(|\widetilde{\bm v}|^4\widetilde{\bm v}) \cdot\overline{\bm v}\right)\,\mathrm{d}x\,\mathrm{d}y\,\mathrm{d}z, \\
        I_2 & :=-\int_{\mathcal{O}}\overline{\widetilde{v}_i\widetilde{v}_j}\, \partial_i(|\widetilde{\bm v}|^4\widetilde{v}_j)\,\mathrm{d}x\,\mathrm{d}y\,\mathrm{d}z, \\
        I_3 & :=\int_{\mathcal{O}}\left(-\int_{-h}^zT\,\mathrm{d}z^\prime +\frac{1}{h}\int_{-h}^0\int_{-h}^zT\,\mathrm{d}z^\prime\,\mathrm{d}z\right) \mathrm{div}(|\widetilde{\bm v}|^4\widetilde{\bm v})\,\mathrm{d}x\,\mathrm{d}y\,\mathrm{d}z.
    \end{aligned}
\end{equation}
Each of the three terms is estimated by the same three steps: the Cauchy--Schwarz inequality in $z$, H\"older's inequality in $\mathcal{D}$, and the two-dimensional Sobolev inequality
\begin{equation*}
    |\psi|_{L^4(\mathcal{D})}\le C|\psi|^{1/2}_{L^2(\mathcal D)}\|\psi\|^{1/2}_{H^1(\mathcal D)}.
\end{equation*}
By Minkowski's integral inequality \cite[Theorem~6.19]{0924.28001},
\begin{equation*}
    \left(\int_{\mathcal{D}}\left(\int_{-h}^0|\widetilde{\bm v}|^6\,\mathrm{d}z\right)^2\,\mathrm{d}x\,\mathrm{d}y\right)^{1/2}\leq C\int_{-h}^0\left(\int_{\mathcal{D}}|\widetilde{\bm v}|^{12}\,\mathrm{d}x\,\mathrm{d}y\right)^{1/2}\,\mathrm{d}z.
\end{equation*}
By \eqref{est:L12}, for almost every $z$,
\begin{equation*}
    \begin{split}
        \left(\int_{\mathcal D}|\widetilde{\bm v}|^{12}\,\mathrm{d}x\,\mathrm{d}y\right)^{1/2} & \leq C\left(\int_{\mathcal D}|\widetilde{\bm v}|^6\,\mathrm{d}x\,\mathrm{d}y\right)^{1/2} \left(\int_{\mathcal D}|\widetilde{\bm v}|^4 |\nabla\widetilde{\bm v}|^2\,\mathrm{d}x\,\mathrm{d}y\right)^{1/2} +C\int_{\mathcal D}|\widetilde{\bm v}|^6\,\mathrm{d}x\,\mathrm{d}y.
    \end{split}
\end{equation*}
Thus, by the Cauchy--Schwarz inequality,
\begin{equation}\label{est:v-tilde-6}
    \left(\int_{\mathcal{D}}\left(\int_{-h}^0|\widetilde{\bm v}|^6\,\mathrm{d}z\right)^2\,\mathrm{d}x\,\mathrm{d}y\right)^{1/2}\leq C|\widetilde{\bm v}|_6^3\left(\int_{\mathcal{O}}|\widetilde{\bm v}|^4|\nabla\widetilde{\bm v}|^2\,\mathrm{d}x\,\mathrm{d}y\,\mathrm{d}z\right)^{1/2}+C|\widetilde{\bm v}|_6^6.
\end{equation}
Similarly, by Minkowski's inequality and \eqref{est:GN-d2-8},
\begin{equation}
    \begin{aligned}
        \left(\int_{\mathcal{D}}\left(\int_{-h}^0|\widetilde{\bm v}|^4\,\mathrm{d}z\right)^2\,\mathrm{d}x\,\mathrm{d}y\right)^{1/2} & \leq C\int_{-h}^0\left(\int_{\mathcal{D}}|\widetilde{\bm v}|^8\,\mathrm{d}x\,\mathrm{d}y\right)^{1/2}\,\mathrm{d}z \\
        & \leq C\int_{-h}^0|\widetilde{\bm v}|_{L^6(\mathcal{D})}^3(|\nabla\widetilde{\bm v}|_{L^2(\mathcal{D})}+|\widetilde{\bm v}|_{L^2(\mathcal{D})})\,\mathrm{d}z \\
        & \leq C|\widetilde{\bm v}|_6^3(|\nabla\widetilde{\bm v}|_2+|\widetilde{\bm v}|_2),
    \end{aligned}
\end{equation}
and
\begin{equation}\label{est:v-tilde-2}
    \begin{aligned}
        \left(\int_{\mathcal{D}}\left(\int_{-h}^0|\widetilde{\bm v}|^2\,\mathrm{d}z\right)^4\,\mathrm{d}x\,\mathrm{d}y\right)^{1/4} & \leq C\int_{-h}^0\left(\int_{\mathcal{D}}|\widetilde{\bm v}|^8\,\mathrm{d}x\,\mathrm{d}y\right)^{1/4}\,\mathrm{d}z \\
        & \leq C\int_{-h}^0|\widetilde{\bm v}|_{L^6(\mathcal{D})}^{3/2}(|\nabla\widetilde{\bm v}|_{L^2(\mathcal{D})}^{1/2}+|\widetilde{\bm v}|_{L^2(\mathcal{D})}^{1/2})\,\mathrm{d}z \\
        & \leq C|\widetilde{\bm v}|_6^{3/2}(|\nabla\widetilde{\bm v}|_2^{1/2}+|\widetilde{\bm v}|_2^{1/2}).
    \end{aligned}
\end{equation}
We write $X_z:=\int_{-h}^0|\nabla\widetilde{\bm v}|^2 |\widetilde{\bm v}|^4\,\mathrm{d}z$ and $X:=\int_{\mathcal{D}}X_z\,\mathrm{d}x\,\mathrm{d}y$.  Applying to each of $I_1,I_2,I_3$ the three steps described above, and using \eqref{est:v-tilde-6}--\eqref{est:v-tilde-2}, we obtain
\begin{equation}\label{est:I-1}
    \begin{aligned}
        |I_1| & \leq C\int_{\mathcal{D}}|\overline{\bm v}| \left(\int_{-h}^0|\nabla\widetilde{\bm v}||\widetilde{\bm v}|^5\,\mathrm{d}z\right)\,\mathrm{d}x\,\mathrm{d}y \\
        & \leq C\int_{\mathcal{D}}|\overline{\bm v}|X_z^{1/2} \left(\int_{-h}^0|\widetilde{\bm v}|^6\,\mathrm{d}z\right)^{1/2}\,\mathrm{d}x\,\mathrm{d}y \\
        & \leq C|\overline{\bm v}|_{L^4(\mathcal{D})}X^{1/2} \left(\int_{\mathcal{D}}\left(\int_{-h}^0|\widetilde{\bm v}|^6\,\mathrm{d}z\right)^2\,\mathrm{d}x\,\mathrm{d}y\right)^{1/4} \\
        & \leq C|\overline{\bm v}|^{1/2}_2|\nabla\overline{\bm v}|^{1/2}_2 (|\widetilde{\bm v}|^{3/2}_6X^{3/4}+|\widetilde{\bm v}|^3_6X^{1/2}),
    \end{aligned}
\end{equation}
\begin{equation}\label{est:I-2}
    \begin{aligned}
        |I_2| & \leq C\int_{\mathcal{D}}\left(\int_{-h}^0|\widetilde{\bm v}|^2\,\mathrm{d}z\right) \left(\int_{-h}^0|\nabla\widetilde{\bm v}||\widetilde{\bm v}|^4\,\mathrm{d}z\right)\,\mathrm{d}x\,\mathrm{d}y \\
        & \leq C\int_{\mathcal{D}}\left(\int_{-h}^0|\widetilde{\bm v}|^2\,\mathrm{d}z\right)X_z^{1/2} \left(\int_{-h}^0|\widetilde{\bm v}|^4\,\mathrm{d}z\right)^{1/2}\,\mathrm{d}x\,\mathrm{d}y \\
        & \leq C\left(\int_{\mathcal{D}}\left(\int_{-h}^0|\widetilde{\bm v}|^2\,\mathrm{d}z\right)^4\,\mathrm{d}x\,\mathrm{d}y\right)^{1/4}X^{1/2} \left(\int_{\mathcal{D}}\left(\int_{-h}^0|\widetilde{\bm v}|^4\,\mathrm{d}z\right)^2\,\mathrm{d}x\,\mathrm{d}y\right)^{1/4} \\
        & \leq C|\widetilde{\bm v}|^3_6(|\nabla\widetilde{\bm v}|_2+|\widetilde{\bm v}|_2)X^{1/2}.
    \end{aligned}
\end{equation}
Here $\|\overline{|T|}\|_{L^4(\mathcal D)}\leq C|T|_2^{1/2}\|T\|_{H^1}^{1/2}\leq C|T|_2^{1/2}\|T\|_{V_2}^{1/2}$, by Minkowski and the two-dimensional Ladyzhenskaya inequality.
\begin{equation}\label{est:I-3}
    \begin{aligned}
        |I_3| & \leq C\int_{\mathcal{D}}\overline{|T|} \left(\int_{-h}^0|\nabla\widetilde{\bm v}||\widetilde{\bm v}|^4\,\mathrm{d}z\right)\,\mathrm{d}x\,\mathrm{d}y \\
        & \leq C\int_{\mathcal{D}}\overline{|T|}X_z^{1/2} \left(\int_{-h}^0|\widetilde{\bm v}|^4\,\mathrm{d}z\right)^{1/2}\,\mathrm{d}x\,\mathrm{d}y \\
        & \leq C\left|\,\overline{|T|}\,\right|_{L^4(\mathcal{D})}X^{1/2} \left(\int_{\mathcal{D}}\left(\int_{-h}^0|\widetilde{\bm v}|^4\,\mathrm{d}z\right)^2\,\mathrm{d}x\,\mathrm{d}y\right)^{1/4} \\
        & \leq C|T|_2^{1/2}\|T\|_{V_2}^{1/2}|\widetilde{\bm v}|^{3/2}_6 (|\nabla\widetilde{\bm v}|^{1/2}_2+|\widetilde{\bm v}|^{1/2}_2)X^{1/2}.
    \end{aligned}
\end{equation}
Put $U=|\widetilde{\bm v}|_6$, $H_v=\|\widetilde{\bm v}\|_{H^1}$ and $A_T=|T|_2\|T\|_{V_2}$. Young's inequality gives
\begin{equation*}
    \begin{aligned}
        |I_1| & \leq\varepsilon X+C_\varepsilon (1+|\overline{\bm v}|_2^2|\nabla\overline{\bm v}|_2^2)U^6, \\
        |I_2| & \leq\varepsilon X+C_\varepsilon U^6H_v^2, \\
        |I_3| & \leq\varepsilon X+C_\varepsilon A_TU^3H_v \leq\varepsilon X+C_\varepsilon(A_T^2+U^6H_v^2).
    \end{aligned}
\end{equation*}
Here the factors involving $\nabla\widetilde{\bm v}$ in \eqref{est:I-2}--\eqref{est:I-3} are $L^2$ norms. Consequently,
\begin{equation}\label{est:L6-determined}
    |I_1+I_2+I_3| \leq\varepsilon X +C_\varepsilon(1+|\overline{\bm v}|^2_2 |\nabla\overline{\bm v}|^2_2 +|\nabla\widetilde{\bm v}|^2_2+|\widetilde{\bm v}|^2_2) |\widetilde{\bm v}|^6_6 +C_\varepsilon|T|_2^2\|T\|_{V_2}^2,
\end{equation}
the term $\varepsilon X$ being absorbed by the dissipation $6\int_0^t\int_{\mathcal{O}}|\nabla\widetilde{\bm v}|^2 |\widetilde{\bm v}|^4\,\mathrm{d}x\,\mathrm{d}y\,\mathrm{d}z\,\mathrm{d}s$ in \eqref{est:L6-ito} upon choosing $\varepsilon$ small enough. All subsequent estimates are stopped: time $t$ denotes $t\wedge\tau_N$, and the same argument applies on stopping subintervals. Set
\begin{equation*}
    R_6=1+|\overline{\bm v}|_{L^2(\mathcal D)}^2 |\nabla\overline{\bm v}|_{L^2(\mathcal D)}^2 +\|\widetilde{\bm v}\|_{H^1}^2, \qquad \mathcal X_6=\int_{\mathcal O}|\widetilde{\bm v}|^4 |\nabla_3\widetilde{\bm v}|^2\,\mathrm{d}x\,\mathrm{d}y\,\mathrm{d}z.
\end{equation*}
After absorbing the deterministic terms, estimate the supremum and the terminal dissipation separately in \eqref{est:L6-ito}. Their sum is bounded by a constant times the initial energy, the integral of $R_6|\widetilde{\bm v}|_6^6+|T|_2^2\|T\|_{V_2}^2$, the nonnegative jump remainder and the supremum of the linear martingale. The stopped version of \eqref{est:T-L4L2} gives
\begin{equation}\label{est:Z-in-L6}
    \mathbb E\int_0^{\tau_N}|T|_2^2\|T\|_{V_2}^2\,\mathrm{d}s \leq C_{\mathcal T}(1+\mathbb E|T_0|_{H_2}^4).
\end{equation}
The joint low-order estimate, proved as in \eqref{eq:low-energy-ito} on the stopped interval, also gives
\begin{equation}\label{est:barotropic-low-order}
    \mathbb E\int_0^{\tau_N}|\overline{\bm v}|_{L^2(\mathcal D)}^2 |\nabla\overline{\bm v}|_{L^2(\mathcal D)}^2\,\mathrm{d}s \leq C_{\mathcal T}(1+\mathbb E|\bm u_0|_H^4).
\end{equation}
Indeed, Jensen's inequality bounds its integrand by $C|\bm u|_H^2\|\bm u\|_V^2$; the stopped $m=4$ estimate is proved explicitly in \eqref{est:global-low-energy} below.

For the stochastic terms, write $\bm w(s)=\widetilde{\bm v}(s)$ and $\bm g(s,\xi)=\widetilde{G_1}(\bm u(s-),\xi)$. All integrands against $\eta$ or $\tilde\eta$ below use $\bm w(s-)$, and the estimates are understood on stopped intervals. Taylor's formula gives
\begin{equation}\label{est:uniform-taylor}
    \big||\bm w+\bm g|_6^6-|\bm w|_6^6-6\langle \bm g,|\bm w|^4\bm w\rangle\big| \leq C\big(|\bm w|_6^4|\bm g|_6^2+|\bm g|_6^6\big).
\end{equation}
By compensation and \eqref{asm:g_L6},
\begin{equation}\label{est:L6-sto-1}
    \begin{aligned}
        & \mathbb E\int_0^t\int_E \big||\bm w+\bm g|_6^6-|\bm w|_6^6-6\langle \bm g,|\bm w|^4\bm w\rangle\big| \,\eta(\mathrm{d}s,\mathrm{d}\xi) \\
        & \qquad\leq C\mathbb E\int_0^t\int_E \big(|\bm w|_6^4|\bm g|_6^2+|\bm g|_6^6\big)\,\nu(\mathrm{d}\xi)\,\mathrm{d}s \leq C\int_0^t\big(1+\mathbb E|\bm w(s)|_6^6\big)\,\mathrm{d}s.
    \end{aligned}
\end{equation}
The Burkholder--Davis--Gundy inequality, Young's inequality and compensation yield, for every $\delta>0$,
\begin{equation}\label{est:L6-sto-2}
    \begin{aligned}
        & \mathbb E\sup_{r\leq t}\left| \int_0^r\int_E\langle \bm g,|\bm w|^4\bm w\rangle \,\tilde\eta(\mathrm{d}s,\mathrm{d}\xi)\right| \\
        & \quad\leq C\mathbb E\left[ \sup_{r\leq t}|\bm w(r)|_6^3 \left(\int_0^t\int_E|\bm w(s-)|_6^4|\bm g(s,\xi)|_6^2 \,\eta(\mathrm{d}s,\mathrm{d}\xi)\right)^{1/2}\right] \\
        & \quad\leq\delta\mathbb E\sup_{r\leq t}|\bm w(r)|_6^6 +C_\delta\mathbb E\int_0^t\int_E |\bm w(s-)|_6^4|\bm g(s,\xi)|_6^2\,\nu(\mathrm{d}\xi)\,\mathrm{d}s \\
        & \quad\leq\delta\mathbb E\sup_{r\leq t}|\bm w(r)|_6^6 +C_\delta\int_0^t\big(1+\mathbb E|\bm w(s)|_6^6\big)\,\mathrm{d}s.
    \end{aligned}
\end{equation}

Choosing $\delta$ sufficiently small in \eqref{est:L6-sto-2} yields
\begin{equation}\label{est:L6-before-gronwall}
    \begin{split}
        \mathbb E\left[\sup_{s\leq t}|\widetilde{\bm v}(s)|_6^6 +\int_0^t\mathcal X_6\,\mathrm{d}s\right] & \leq C\mathbb E|\widetilde{\bm v}_0|_6^6 +C_{\mathcal T}(1+\mathbb E|T_0|_{H_2}^4) \\
        & \quad+C\mathbb E\int_0^t R_6(s) |\widetilde{\bm v}(s)|_6^6\,\mathrm{d}s.
    \end{split}
\end{equation}
Set
\begin{equation*}
    \sigma_{N,K}=\tau_N\wedge\inf\left\{t:\int_0^tR_6\,\mathrm{d}s\geq K\right\}.
\end{equation*}
Applying the preceding estimates on stopped subintervals and using Lemma \ref{lem:sto_gron} gives
\begin{equation}\label{est:v-tilde-L6}
    \mathbb E\left[\sup_{t\leq\sigma_{N,K}}|\widetilde{\bm v}(t)|_6^6 +\int_0^{\sigma_{N,K}}\int_{\mathcal O} |\widetilde{\bm v}|^4|\nabla_3\widetilde{\bm v}|^2\,\mathrm{d}x\,\mathrm{d}y\,\mathrm{d}z\,\mathrm{d}s\right] \leq C_{\mathcal T,K},
\end{equation}
where the constant also depends on the initial moments and is independent of $N$.

For each fixed $K$, estimate \eqref{est:v-tilde-L6} is uniform in $N$.

\subsubsection{The \texorpdfstring{$H^1$}{H1} estimate}

\begin{proposition}\label{prop:H1-global}
    Under Assumption \ref{asm:global}, let $(\bm u,\xi)$ be a maximal pathwise solution. Assume $\mathbb E\|\bm u_0\|_V^6<\infty$ and set $\mathcal E_N=\sup_{t\leq\tau_N}\|\bm u(t)\|_V^2+ \int_0^{\tau_N}|A\bm u|_H^2\,\mathrm{d}t$. Then
    \begin{equation}\label{est:H1-global}
        \lim_{M\to\infty}\sup_{N\geq1}\mathbb P(\mathcal E_N>M)=0.
    \end{equation}
\end{proposition}

The proof combines four component energy estimates, with dissipations $|S\overline{\bm v}|_2^2$, $|\nabla_3\bm v_z|_2^2$, $|A_h\bm v|_2^2+|\nabla\bm v_z|_2^2$ and $|A_2T|_2^2$. The temperature energy includes the Robin surface term.

Let $P_\sigma$ be the Leray projection on $\mathcal D$ and let $S=-P_\sigma\Delta$ be the Stokes operator on divergence-free fields on $\mathcal D$, with $\bm v\cdot\bm n=0$ and $(\partial_{\bm n}\bm v)\times\bm n=0$. Its form is $\int_{\mathcal D}\nabla\bm v:\nabla\bm w\,\mathrm{d}x\,\mathrm{d}y$. Define
\begin{equation*}
    A_h\bm v=S\overline{\bm v}-\Delta\widetilde{\bm v}.
\end{equation*}
Then $A_1=A_h-\partial_{zz}$, $\langle A_h\bm v,\bm v\rangle_{H_1}=|\nabla\bm v|_2^2$, and
\begin{equation*}
    \langle A_1\bm v,A_h\bm v\rangle_{H_1} =|A_h\bm v|_2^2+|\nabla\bm v_z|_2^2, \qquad \bm v\in\mathscr D(A_1).
\end{equation*}
The horizontal and vertical realizations commute on the cylinder. Since $\overline{A_h\bm v}=S\overline{\bm v}$ is divergence free and tangent to $\partial\mathcal D$, $\int_{\mathcal O}\nabla\phi(x,y)\cdot A_h\bm v\,\mathrm{d}x\,\mathrm{d}y\,\mathrm{d}z =h\int_{\mathcal D}\nabla\phi\cdot S\overline{\bm v}\,\mathrm{d}x\,\mathrm{d}y=0$. The $z$-dependent buoyancy is estimated separately below.

All energy inequalities below are stopped at a common bounded stopping time within the local existence interval; we omit this stopping from the notation. Stochastic integrands are evaluated at left limits. For each quadratic energy $Y=|L\bm u|_2^2$ used below, the noise terms are
\begin{equation*}
    2\int_0^t\int_E\langle L\bm u(s-),LG(\bm u(s-),\xi)\rangle\,\tilde\eta(\mathrm{d}s,\mathrm{d}\xi) +\int_0^t\int_E|LG(\bm u(s-),\xi)|_2^2\,\eta(\mathrm{d}s,\mathrm{d}\xi).
\end{equation*}
Here the energy maps are $L\bm u=\nabla\overline{\bm v}$, $\partial_z\bm v$, $\nabla\bm v$ and $A_2^{1/2}T$, respectively. Each is bounded from $V$ into its corresponding $L^2$ space. In particular, testing the drift with $A_h\bm v$ or $A_2T$ gives the jump remainders $|\nabla G_1|_2^2$ and $|A_2^{1/2}G_2|_2^2=a_2(G_2,G_2)$. By \eqref{asm:g_moment-V} with $p=1$, the Burkholder--Davis--Gundy inequality and compensation bound their expected supremum by
\begin{equation}\label{est:h1-noise}
    \delta\mathbb E\sup_{s\leq t}Y(s) +C_\delta\mathbb E\int_0^t(1+\|\bm u\|_V^2)\,\mathrm{d}s, \qquad\delta>0.
\end{equation}
The quadratic identities follow by spectral approximation, as in Lemma \ref{lem:strong-cadlag}.

\subsubsection{The barotropic gradient}
Set $Y_b=|\nabla\overline{\bm v}|_2^2$ and $\mathcal X=\int_{\mathcal O}|\widetilde{\bm v}|^4|\nabla\widetilde{\bm v}|^2\,\mathrm{d}x\,\mathrm{d}y\,\mathrm{d}z$. Testing the projected barotropic equation with $S\overline{\bm v}$ gives the dissipation $|S\overline{\bm v}|_2^2$. The two-dimensional interpolation inequality yields
\begin{equation}\label{id:nabla-v-bar-ito}
    |\langle(\overline{\bm v}\cdot\nabla)\overline{\bm v},S\overline{\bm v}\rangle| \leq C|\overline{\bm v}|_2^{1/2}|\nabla\overline{\bm v}|_2|S\overline{\bm v}|_2^{3/2} \leq\varepsilon|S\overline{\bm v}|_2^2+C_\varepsilon|\overline{\bm v}|_2^2Y_b^2.
\end{equation}
For the Reynolds stress, Cauchy--Schwarz gives
\begin{equation*}
    \left\|\int_{-h}^0|\widetilde{\bm v}||\nabla\widetilde{\bm v}|\,\mathrm{d}z\right\|_{L^2(\mathcal D)}^2 \leq C|\nabla\widetilde{\bm v}|_2\mathcal X^{1/2}.
\end{equation*}
Its pairing with $S\overline{\bm v}$ is therefore bounded by $\varepsilon|S\overline{\bm v}|_2^2+C_\varepsilon(|\nabla\widetilde{\bm v}|_2^2+\mathcal X)$. Using \eqref{est:h1-noise}, we obtain
\begin{equation}\label{est:barotropic-local}
    \mathbb E\sup_{s\leq t}Y_b(s)+\mathbb E\int_0^t|S\overline{\bm v}|_2^2\,\mathrm{d}s \leq C\mathbb EY_b(0) +C\mathbb E\int_0^t\big[(1+|\overline{\bm v}|_2^2Y_b)Y_b +1+\|\bm u\|_V^2+\mathcal X\big]\,\mathrm{d}s.
\end{equation}

\subsubsection{The vertical derivative}
Writing $\bm v_z=\partial_z\bm v$, differentiation of \eqref{eqn:r-1} gives
\begin{equation}\label{eqn:v-z}
    \begin{aligned}
        \mathrm{d}\bm v_z & +\big[(\bm v\cdot\nabla)\bm v_z+\Phi(\bm v)\partial_z\bm v_z\big]\,\mathrm{d}t \\
        & +\big[(\bm v_z\cdot\nabla)\bm v-(\mathrm{div}\bm v)\bm v_z +f\bm k\times\bm v_z-\nabla T+L_1\bm v_z\big]\,\mathrm{d}t \\
        & =\int_E\partial_zG_1\,\tilde\eta(\mathrm{d}t,\mathrm{d}\xi).
    \end{aligned}
\end{equation}
The transport terms cancel when tested against $\bm v_z$. Integration by parts and the zero normal flux give
\begin{equation*}
    \left|\int_{\mathcal O}\big[-(\bm v_z\cdot\nabla)\bm v+(\mathrm{div}\bm v)\bm v_z\big]\cdot\bm v_z\,\mathrm{d}x\,\mathrm{d}y\,\mathrm{d}z\right| \leq C|\bm v|_6|\bm v_z|_2^{1/2}|\nabla_3\bm v_z|_2^{3/2} \leq\varepsilon|\nabla_3\bm v_z|_2^2+C_\varepsilon|\bm v|_6^4|\bm v_z|_2^2.
\end{equation*}
The temperature term is bounded by $C|T|_2|\nabla\bm v_z|_2$. Thus
\begin{equation}\label{id:v-z-ito}
    \mathbb E\sup_{s\leq t}|\bm v_z(s)|_2^2+\mathbb E\int_0^t|\nabla_3\bm v_z|_2^2\,\mathrm{d}s \leq C\mathbb E|\bm v_z(0)|_2^2 +C\mathbb E\int_0^t\big[|\bm v|_6^4|\bm v_z|_2^2+1+\|\bm u\|_V^2\big]\,\mathrm{d}s.
\end{equation}

\subsubsection{The horizontal gradient}
Put $X=|\nabla\bm v|_2$, $Q=|A_h\bm v|_2$ and
\begin{equation*}
    D_h^2=Q^2+|\nabla\bm v_z|_2^2,\qquad R_v=1+|\bm v|_6^4+(1+|\bm v_z|_2^2)(|\bm v_z|_2^2+|\nabla\bm v_z|_2^2).
\end{equation*}
The anisotropic estimate used in \eqref{est:B-H-norm}, together with horizontal elliptic regularity, gives
\begin{equation}\label{est:nabla-v-determined}
    \begin{aligned}
        |\langle B_1(\bm v,\bm v),A_h\bm v\rangle| & \leq C|\bm v|_6X^{1/2}(X+D_h)^{1/2}Q \\
        & \quad+CX^{1/2}(X+Q)^{1/2}|\bm v_z|_2^{1/2}(|\bm v_z|_2+|\nabla\bm v_z|_2)^{1/2}Q \\
        & \leq\varepsilon D_h^2+C_\varepsilon R_vX^2.
    \end{aligned}
\end{equation}
To verify the first line, use $|\nabla\bm v|_3\leq C X^{1/2}(X+D_h)^{1/2}$. For the vertical term, the functions $F_f$ used in the proof of Lemma \ref{lem:Euler_bi_op} give
\begin{equation*}
    \left|\int_{\mathcal O}\Phi(\bm v)\bm v_z\cdot A_h\bm v\,\mathrm{d}x\,\mathrm{d}y\,\mathrm{d}z\right| \leq C\|F_{\operatorname{div}\bm v}\|_{L^4(\mathcal D)} \|F_{\bm v_z}\|_{L^4(\mathcal D)}Q.
\end{equation*}
Horizontal elliptic regularity gives $\|\nabla\bm v\|_{L^2_zH^1_{x,y}}\leq C(X+Q)$, which proves the second line. With $Z=|\bm v_z|_2$ and $W=|\nabla\bm v_z|_2$, its two contributions are bounded by $CX[Z(Z+W)]^{1/2}Q$ and $CX^{1/2}[Z(Z+W)]^{1/2}Q^{3/2}$. Young's inequality with exponents $(2,2)$ and $(4,4/3)$ gives $\varepsilon Q^2+C_\varepsilon[Z(Z+W)+Z^2(Z+W)^2]X^2$. Since the bracket is at most $C(1+Z^2)(Z^2+W^2)$, this proves the last line of \eqref{est:nabla-v-determined}. The Coriolis and buoyancy terms are bounded by $\varepsilon Q^2+C_\varepsilon(|\bm v|_2^2+|\nabla T|_2^2)$. Testing with $A_h\bm v$ and using \eqref{est:h1-noise} yields
\begin{equation}\label{est:horizontal-local}
    \mathbb E\sup_{s\leq t}X(s)^2+\mathbb E\int_0^tD_h^2\,\mathrm{d}s \leq C\mathbb EX(0)^2+C\mathbb E\int_0^t\big[R_vX^2+1+\|\bm u\|_V^2\big]\,\mathrm{d}s.
\end{equation}

\subsubsection{The temperature gradient}
The energy is $\|T\|_{V_2}^2=a_2(T,T)$, including the Robin surface term, and its derivative on $\mathscr D(A_2)$ is $2A_2T$. By the component version of \eqref{est:B-H-norm},
\begin{equation}\label{est:T-V-determined}
    |\langle B_2(\bm v,T),A_2T\rangle| \leq\varepsilon|A_2T|_2^2+C_\varepsilon R_T\|T\|_{V_2}^2, \qquad R_T=|\bm v|_6^4+\|\bm v\|_{V_1}^2|A_1\bm v|_2^2.
\end{equation}
Indeed, the vertical contribution is bounded by $C\|\bm v\|_{V_1}^{1/2}|A_1\bm v|_2^{1/2}\|T\|_{V_2}^{1/2}|A_2T|_2^{3/2}$. Consequently,
\begin{equation}\label{est:temperature-local}
    \mathbb E\sup_{s\leq t}\|T(s)\|_{V_2}^2+\mathbb E\int_0^t|A_2T|_2^2\,\mathrm{d}s \leq C\mathbb E\|T_0\|_{V_2}^2 +C\mathbb E\int_0^t\big[R_T\|T\|_{V_2}^2+1+\|\bm u\|_V^2\big]\,\mathrm{d}s.
\end{equation}

\begin{proof}[Proof of Proposition \ref{prop:H1-global}]
    Cancellation of $B$, the bound on $R$, and the $H$-moment assumptions, applied to $|\bm u|_H^2$ and $|\bm u|_H^4$, give
    \begin{equation}\label{est:global-low-energy}
        \sup_N\mathbb E\left[\sup_{t\leq\tau_N}|\bm u(t)|_H^4 +\int_0^{\tau_N}(1+|\bm u|_H^2)\|\bm u\|_V^2\,\mathrm{d}t\right] \leq C_{\mathcal T}(1+\mathbb E|\bm u_0|_H^4).
    \end{equation}
    Indeed, the drift contribution of $R$ is absorbed into half the dissipation plus $C(|\bm u|_H^2+|\bm u|_H^4)$. The quadratic and quartic jump remainders are bounded by the second and fourth noise moments. More explicitly, applying \eqref{eq:low-energy-ito} to the original solution, with $m=2,4$, gives, for any further localization $\sigma$,
    \begin{equation*}
        \begin{split}
            & \mathbb E\left[\sup_{s\leq t\wedge\tau_N\wedge\sigma}|\bm u(s)|_H^m +\int_0^{t\wedge\tau_N\wedge\sigma}|\bm u|_H^{m-2}\|\bm u\|_V^2\,\mathrm{d}s\right] \\
            & \qquad\leq C_m\mathbb E|\bm u_0|_H^m +C_m\int_0^t\left(1+\mathbb E\sup_{r\leq s\wedge\tau_N\wedge\sigma}|\bm u(r)|_H^m\right)\mathrm{d}s.
        \end{split}
    \end{equation*}
    The constants do not depend on $N$ or $\sigma$. Deterministic Gr\"onwall and Fatou's lemma remove $\sigma$; adding the cases $m=2,4$ proves \eqref{est:global-low-energy}. Applying the same argument to the temperature equation proves the stopped version of \eqref{est:T-L4L2} used in \eqref{est:Z-in-L6}.

    We use the following localization argument repeatedly. Suppose a local energy inequality has coefficient $r\geq0$ and additive term $z\geq0$, and both $\int_0^{\tau_N}r$ and $\int_0^{\tau_N}z$ are bounded in probability, uniformly in $N$. Stop also when their sum first reaches $K$. Lemma \ref{lem:sto_gron} bounds the expected stopped energy by $C_K(1+\mathbb E Y(0)+K)$. Consequently its unstopped energy $E_N$ satisfies
    \begin{equation}\label{est:localization-transfer}
        \mathbb P(E_N>M)\leq \mathbb P\left(\int_0^{\tau_N}(r+z)\,\mathrm{d}t\geq K\right) +\frac{C_K(1+\mathbb E Y(0)+K)}{M}.
    \end{equation}
    First let $M\to\infty$, then $K\to\infty$. The local inequalities hold on arbitrary stopping subintervals, as required by that lemma. The additional stops use continuous accumulated integrals; terminal jumps remain in the energy estimate.

    The integral of the coefficient in \eqref{est:L6-before-gronwall} is bounded by $C\int_0^{\tau_N}[1+(1+|\bm u|_H^2)\|\bm u\|_V^2]\,\mathrm{d}s$. Thus \eqref{est:global-low-energy}, \eqref{est:v-tilde-L6} and \eqref{est:localization-transfer} imply that
    \begin{equation*}
        \sup_{t\leq\tau_N}|\widetilde{\bm v}(t)|_6^6 +\int_0^{\tau_N}\!\int_{\mathcal O} |\widetilde{\bm v}|^4|\nabla_3\widetilde{\bm v}|^2 \,\mathrm{d}x\,\mathrm{d}y\,\mathrm{d}z\,\mathrm{d}s
    \end{equation*}
    is bounded in probability uniformly in $N$. In \eqref{est:barotropic-local}, the integral of $1+|\overline{\bm v}|_2^2Y_b$ is controlled by \eqref{est:global-low-energy}, and the additive integral is controlled by the preceding step. Hence $\sup Y_b+\int|S\overline{\bm v}|_2^2$ is bounded in probability. The two-dimensional Sobolev inequality gives
    \begin{equation*}
        \sup_{t\leq\tau_N}|\bm v(t)|_6^4 \leq C\sup_{t\leq\tau_N} \big(\|\overline{\bm v}(t)\|_{H^1(\mathcal D)}^4 +|\widetilde{\bm v}(t)|_6^4\big).
    \end{equation*}
    Applying the same argument to \eqref{id:v-z-ito} now controls $\sup|\bm v_z|_2^2+\int|\nabla_3\bm v_z|_2^2$ in probability. In particular,
    \begin{equation*}
        \int_0^{\tau_N}R_v\leq \mathcal T(1+\sup|\bm v|_6^4) +(1+\sup|\bm v_z|_2^2) \left(\mathcal T\sup|\bm v_z|_2^2+\int_0^{\tau_N}|\nabla\bm v_z|_2^2\right)
    \end{equation*}
    is bounded in probability. Estimate \eqref{est:horizontal-local} therefore controls $\sup\|\bm v\|_{V_1}^2+\int|A_1\bm v|_2^2$; here $A_1=A_h-\partial_{zz}$ and the vertical dissipation controls $|\partial_{zz}\bm v|_2^2$. Finally,
    \begin{equation*}
        \int_0^{\tau_N}R_T\leq\mathcal T\sup|\bm v|_6^4 +\sup\|\bm v\|_{V_1}^2\int_0^{\tau_N}|A_1\bm v|_2^2
    \end{equation*}
    is bounded in probability. Applying \eqref{est:localization-transfer} to \eqref{est:temperature-local} proves \eqref{est:H1-global}. All suprema and integrals in this argument are over $[0,\tau_N]$.
\end{proof}

\subsubsection{Proof of the global theorem}

\begin{theorem}[Global well-posedness]\label{thm:global}
    Let Assumptions \ref{asm:martingale}, \ref{asm:A-est} and \ref{asm:global} hold and let $\bm u_0\in L^2(\Omega,\mathcal F_0;V)$.  Then the unique maximal pathwise solution $(\bm u,\xi)$ of \eqref{eqn:abstract} given by Theorem \ref{thm:wp-maximal} is global:
    \begin{equation}\label{eq:global}
        \xi=+\infty\quad\mathbb P\text{-almost surely},
    \end{equation}
    and for every $\mathcal T>0$
    \begin{equation}\label{eq:global-bound}
        \sup_{t\in[0,\mathcal T]}\|\bm u(t)\|_V^2 +\int_0^{\mathcal T}|A\bm u|_H^2\,\mathrm{d}t<\infty \quad\mathbb P\text{-almost surely}.
    \end{equation}
\end{theorem}

\begin{proof}
    First suppose $\bm u_0\in L^6(\Omega,\mathcal F_0;V)$. The energies $\mathcal E_N$ increase with $N$. By \eqref{est:H1-global}, their limit is finite almost surely:
    \begin{equation}\label{eq:bound-on-xi}
        \mathbb P\left(\lim_N\mathcal E_N=\infty\right) \leq\lim_{M\to\infty}\sup_N\mathbb P(\mathcal E_N>M)=0.
    \end{equation}
    On $\{\xi\leq\mathcal T\}$, the blow-up criterion \eqref{eq:blow-up} and $\tau_N\uparrow\xi\wedge\mathcal T$ force $\lim_N\mathcal E_N=\infty$. Therefore $\mathbb P(\xi\leq\mathcal T)=0$ for every $\mathcal T>0$. This proves \eqref{eq:global} and \eqref{eq:global-bound} in this case. For general square-integrable initial data, condition on $\Omega_j=\{j-1\leq\|\bm u_0\|_V<j\}$ whenever $\mathbb P(\Omega_j)>0$. Since $\Omega_j\in\mathcal F_0$, the independent-increment property of $\eta$ and its compensator are unchanged under $\mathbb P(\cdot\mid \Omega_j)$. The initial datum is bounded under this law, and the local energy bounds and blow-up criterion remain valid. Applying the preceding argument under each conditional law proves both assertions on $\bigcup_{j\geq1}\Omega_j=\Omega$.
\end{proof}

    \section{Large deviation principle}\label{sec:ldp}
We establish the Freidlin--Wentzell large deviation principle for the family \eqref{eqn:spe} by the weak convergence approach of \cite{0936.60059,0994.60028,0904.60001}.

Let $\{X^\varepsilon\}$ be a family of random variables defined on a probability space $(\Omega,\mathcal{F},\mathbb{P})$ taking values in some Polish space $\mathbb U$.

\begin{definition}[Rate function]
    A function $I:\mathbb U\to[0,\infty]$ is called a rate function if it is lower semicontinuous and not identically $+\infty$. It is good if every sublevel set $\{x\in\mathbb U:I(x)\leq M\}$, $M<\infty$, is compact.
\end{definition}

\begin{definition}[Large deviation principle]
    The family $\{X^\varepsilon\}$ is said to satisfy the large deviation principle with rate function $I:\mathbb U\to[0,\infty]$ if the following conditions hold:
    \begin{enumerate}
        \item (Upper bound) For each closed subset $F$ of $\mathbb U$,
            \begin{equation*}
                \limsup_{\varepsilon\to 0}\varepsilon\log\mathbb{P}(X^\varepsilon\in F)\leq -\inf_{x\in F}I(x).
            \end{equation*}
        \item (Lower bound) For each open subset $G$ of $\mathbb U$,
            \begin{equation*}
                \liminf_{\varepsilon\to 0}\varepsilon\log\mathbb{P}(X^\varepsilon\in G)\geq -\inf_{x\in G}I(x).
            \end{equation*}
    \end{enumerate}
\end{definition}

On a Polish space, the large deviation principle with a good rate function is equivalent to the Laplace principle with the same rate function \cite[Section~1.2]{0904.60001}.

\begin{definition}
    The family $\{X^\varepsilon\}$ is said to satisfy the Laplace principle in $\mathbb U$ with rate function $I:\mathbb U\to[0,\infty]$ if for all bounded continuous functions $h:\mathbb U\to\mathbb{R}$,
    \begin{equation*}
        \lim_{\varepsilon\to 0}\varepsilon\log\mathbb{E}\left(\exp\left(-\frac{h(X^\varepsilon)}{\varepsilon}\right)\right)=-\inf_{x\in\mathbb U}\{h(x)+I(x)\}.
    \end{equation*}
\end{definition}

\subsection{Poisson random measure and a variational representation}
Throughout, $E$ is a locally compact Polish space. Let $\mathcal{M}_{FC}(E)$ be the space of all measures $\nu$ on the Borel $\sigma$-algebra $\mathcal{B}(E)$ such that $\nu(K)<\infty$ for all compact subsets $K$ of $E$. We endow $\mathcal{M}_{FC}(E)$ with the weakest topology such that for all continuous functions $f$ on $E$ with compact support, the map $\nu\mapsto\int_Ef(x)\nu(\mathrm{d}x)$ is continuous. It is well known that with this topology, $\mathcal{M}_{FC}(E)$ is a Polish space; see \cite{1231.60018}.

Fix $\mathcal{T}\in(0,\infty)$, and let $E_{\mathcal{T}}=[0,\mathcal{T}]\times E$. Fix a measure $\nu\in\mathcal{M}_{FC}(E)$, and let $\nu_{\mathcal{T}}=\mathrm{Leb}_{[0,\mathcal{T}]}\otimes\nu$, where $\mathrm{Leb}_{[0,\mathcal{T}]}$ is Lebesgue measure on $[0,\mathcal{T}]$. We recall the definition of a Poisson random measure; see, e.g., \cite[Definition I.8.1]{0684.60040}.

\begin{definition}[Poisson random measure]
    A Poisson random measure $\eta$ on $E_{\mathcal{T}}$ with intensity measure $\nu_{\mathcal{T}}$ is a $\mathcal{M}_{FC}(E_{\mathcal{T}})$-valued random variable such that
    \begin{enumerate}
        \item for each $B\in\mathcal{B}(E_{\mathcal{T}})$ with $\nu_{\mathcal{T}}(B)<\infty$, $\eta(B)$ has a Poisson distribution with mean $\nu_{\mathcal{T}}(B)$;
        \item for disjoint $B_1,B_2, \cdots, B_k\in\mathcal{B}(E_{\mathcal{T}})$, $\eta(B_1),\eta(B_2),\cdots,\eta(B_k)$ are mutually independent random variables.
    \end{enumerate}
\end{definition}

Denote by $\mathbb{P}$ the probability measure induced by $\eta$ on $(\mathcal{M}_{FC}(E_{\mathcal{T}}),\mathcal{B}(\mathcal{M}_{FC}(E_{\mathcal{T}})))$. Let $\mathbb{M}=\mathcal{M}_{FC}(E_{\mathcal{T}})$. Then $\mathbb{P}$ is the unique probability measure on $(\mathbb{M},\mathcal{B}(\mathbb{M}))$ under which the canonical map $N:\mathbb{M}\to\mathbb{M}$, $N(m):=m$, is a Poisson random measure with intensity measure $\nu_{\mathcal{T}}$. With applications to large deviations in mind, we also consider, for $\theta>0$, probability measures $\mathbb{P}_\theta$ on $(\mathbb{M}, \mathcal{B}(\mathbb{M}))$ under which $N$ is a Poisson random measure with intensity $\theta\nu_{\mathcal{T}}$. The corresponding expectation operators will be denoted by $\mathbb{E}$ and $\mathbb{E}_\theta$, respectively.

Set
\begin{equation*}
    Y:=E\times[0,\infty),\quad Y_{\mathcal{T}}=[0,\mathcal{T}]\times Y.
\end{equation*}

Let $\bar{\mathbb M}=\mathcal M_{FC}(Y_{\mathcal T})$ carry the law $\bar{\mathbb P}$ under which the canonical measure $\bar N$ is Poisson with intensity $\bar\nu_{\mathcal T}=\mathrm{d}t\otimes\nu\otimes\mathrm{d}r$. Write $\bar{\mathbb E}$ for expectation under $\bar{\mathbb P}$, and let $\bar{\mathbb F}:=\{\bar{\mathcal F}_t\}_{0\leq t\leq\mathcal T}$ be the usual augmentation of the natural filtration of $\bar N$. Set $\bar{\mathcal F}=\bar{\mathcal F}_{\mathcal T}$. Denote its predictable $\sigma$-field by $\bar{\mathcal P}$. Let $\bar{\mathcal{A}}$ be the class of all $(\bar{\mathcal{P}}\otimes\mathcal{B}(E))$-measurable maps $\varphi: E_{\mathcal{T}}\times\bar{\mathbb{M}}\to[0,\infty)$. For $\varphi\in\bar{\mathcal{A}}$, define the controlled random measure $N^\varphi$ on $E_{\mathcal{T}}$ by
\begin{equation*}
    N^\varphi((0,t]\times U) =\int_{(0,t]\times U\times(0,\infty)} \mathbbm{1}_{\{r\leq\varphi(s,\xi)\}} \,\bar N(\mathrm{d}s,\mathrm{d}\xi,\mathrm{d}r), \qquad U\in\mathcal B(E).
\end{equation*}
When $\varphi(s,\xi,\bar{m})\equiv\theta\in(0,\infty)$, we write $N^\varphi=N^\theta$. Note that $N^\theta$ has the same distribution with respect to $\bar{\mathbb{P}}$ as $N$ has with respect to $\mathbb{P}_\theta$. Define the local rate function $l:[0,\infty)\to [0, \infty)$ for a standard Poisson process by
\begin{equation*}
    l(r)=r\log r-r+1,\quad r>0,\qquad l(0)=1.
\end{equation*}
For any $\varphi\in\bar{\mathcal{A}}$ the quantity
\begin{equation}
    L_{\mathcal{T}}(\varphi)=\int_{E_{\mathcal{T}}}l(\varphi(s,\xi))\nu_{\mathcal{T}}(\mathrm{d}s,\mathrm{d}\xi)
\end{equation}
is well defined as a $[0, \infty]$-valued random variable. The following is a representation formula proved in \cite{1231.60018}.
\begin{theorem}
    Let $F$ be a bounded measurable function on $\mathbb{M}$. Then for $\theta>0$,
    \begin{equation*}
        -\log\mathbb{E}_\theta(e^{-F(N)})=-\log\bar{\mathbb{E}}(e^{-F(N^\theta)})=\inf_{\varphi\in\bar{\mathcal{A}}}\bar{\mathbb{E}}(\theta L_{\mathcal{T}}(\varphi)+F(N^{\theta\varphi})).
    \end{equation*}
\end{theorem}

In the variational formula, controls with infinite expected entropy are assigned cost $+\infty$. Finite entropy implies local finiteness of the controlled intensity, by Young's entropy inequality on compact sets.

Write $\tilde N=\bar N-\bar\nu_{\mathcal T}$. For a predictable control $\varphi$, the predictable compensator of $N^\varphi$ is $\varphi(s,\xi)\,\mathrm{d}s\,\nu(\mathrm{d}\xi)$, on sets where it is locally integrable. In particular, if $\int_0^{\mathcal T}\int_U\varphi\,\mathrm{d}\nu\,\mathrm{d}s<\infty$ almost surely, then
\begin{equation*}
    \tilde N^\varphi((0,t]\times U) =N^\varphi((0,t]\times U)-\int_0^t\int_U\varphi\,\mathrm{d}\nu\,\mathrm{d}s.
\end{equation*}
For a predictable Hilbert-space-valued integrand $g$ satisfying $\int_0^{\mathcal T}\int_E\|g\|^2\varphi\,\mathrm{d}\nu\,\mathrm{d}s<\infty$ almost surely, define the locally square-integrable integral by
\begin{equation*}
    \int_0^t\int_E g\,\mathrm{d}\tilde N^\varphi :=\int_0^t\int_E\int_0^\infty g(s,\xi)\mathbbm1_{[0,\varphi(s,\xi)]}(r) \,\tilde N(\mathrm{d}s,\mathrm{d}\xi,\mathrm{d}r).
\end{equation*}
The controls and stopped integrands used below satisfy these conditions. For constant $\varphi=\theta>0$, $N^\theta$ is a Poisson random measure with intensity measure $\theta\nu_{\mathcal T}$, and the integral above is taken with respect to its compensated version.

Throughout this section the initial datum $\bm u_0\in V$ is fixed and deterministic. The notation $\Upsilon_{\mathcal{T}}(V)$ denotes the path space
\begin{equation}
    \Upsilon_{\mathcal{T}}(V):=D([0,\mathcal{T}];V)\cap L^2(0,\mathcal{T};\mathscr D(A)),
\end{equation}
where $D([0,\mathcal{T}];V)$ is the space of c\`adl\`ag $V$-valued paths, endowed with the Skorokhod topology, and $L^2(0,\mathcal{T};\mathscr D(A))$ carries its usual norm. Then $\Upsilon_{\mathcal{T}}(V)$ is a closed subspace of the product $D([0,\mathcal{T}];V)\times L^2(0,\mathcal{T};\mathscr D(A))$, and hence a Polish space for the metric
\begin{equation*}
    d_\Upsilon(\bm u,\tilde{\bm u}):=d_{J_1}(\bm u,\tilde{\bm u})+\left(\int_0^{\mathcal{T}}\|\bm u(t)-\tilde{\bm u}(t)\|_{\mathscr D(A)}^2\,\mathrm{d}t\right)^{1/2},
\end{equation*}
where $d_{J_1}$ denotes the Skorokhod metric on $D([0,\mathcal{T}];V)$. The large deviation principle for the primitive equations is formulated on $\Upsilon_{\mathcal{T}}(V)$. By Lemma \ref{lem:strong-cadlag} and Theorem \ref{thm:global}, the solutions of \eqref{eqn:spe} belong almost surely to $\Upsilon_{\mathcal T}(V)$.

Fix $\varepsilon>0$. On the canonical stochastic basis introduced above, consider the SPE
\begin{equation}\label{eqn:spe}
    \begin{cases}
        \displaystyle\,\mathrm{d}\bm u^\varepsilon(t)+(A(\bm u^\varepsilon(t))+B(\bm u^\varepsilon(t))+R(\bm u^\varepsilon(t)))\,\mathrm{d}t=\varepsilon\int_EG(\bm u^\varepsilon(t-),\xi)\tilde{N}^{1/\varepsilon}(\mathrm{d}t,\mathrm{d}\xi), \\
        \bm u^\varepsilon(0)=\bm u_0\in V.
    \end{cases}
\end{equation}

By Theorem \ref{thm:global}, applied with the noise scaled by $\varepsilon$, for each $\varepsilon\in(0,1]$ there exists a unique global strong solution $\bm u^\varepsilon$ of the above SPE on $[0,\mathcal{T}]$.

\subsection{A general large deviation result}
We use the variational criterion of \cite{1259.60065,1231.60018}. Let $\mathbb U$ be a Polish space and let $\mathcal G^\varepsilon:\mathbb M\to\mathbb U$ be measurable maps. Set $Z^\varepsilon=\mathcal G^\varepsilon(\varepsilon N^{\varepsilon^{-1}})$. For nonnegative measurable controls $g$, define
\begin{equation}
    S^N = \{g:E_{\mathcal{T}}\to [0, \infty) : L_{\mathcal{T}}(g) \leq N \}.
\end{equation}
A function $g\in S^N$ can be identified with a measure $\nu_{\mathcal{T}}^g \in \mathbb{M}$, defined by
\begin{equation}\label{def:nu_T^g}
    \nu_{\mathcal{T}}^g(B) = \int_B g(s, x) \nu_{\mathcal{T}}(\mathrm{d}s,\,\mathrm{d}x), \quad B \in \mathcal{B}(E_{\mathcal{T}}).
\end{equation}
This identification induces a topology on $S^N$ under which $S^N$ is a compact space; see the Appendix of \cite{1259.60065} for details. Throughout we use this topology on $S^N$. Define $\mathbb{S}=\bigcup_{N=1}^\infty S^N$, and let
\begin{equation*}
    \mathcal{U}^N:=\{\varphi\in \bar{\mathcal{A}} : \varphi(w) \in S^N , \bar{\mathbb{P}}\text{-a.e. }w\}.
\end{equation*}
Choose increasing compact sets $K_m$ with $\bigcup_mK_m=E$. Set
\begin{equation}\label{def:restricted-controls}
    \begin{split}
        \bar{\mathcal A}_{b,m} & =\{\varphi\in\bar{\mathcal A}:m^{-1}\leq\varphi\leq m \text{ on }[0,\mathcal T]\times K_m,\quad \varphi=1\text{ on }[0,\mathcal T]\times K_m^c\}, \\
        \bar{\mathcal A}_b & =\bigcup_{m\geq1}\bar{\mathcal A}_{b,m}, \qquad \widetilde{\mathcal U}^N =\mathcal U^N\cap\bar{\mathcal A}_b.
    \end{split}
\end{equation}
The bounds hold almost surely. We use the restricted-control criterion of \cite[Theorem~2.4]{1259.60065}; the limit control need only belong to $\mathcal U^N$.

The criterion involves a limiting map $\mathcal G^0$ and the following two convergence properties.

\begin{condition}\label{cond:ldp}
    There exists a measurable map $\mathcal{G}^0:\mathbb{M}\to\mathbb{U}$ such that, for every $N\in\mathbb N$, the following hold:
    \begin{enumerate}
        \item If $g_n,g\in S^N$ and $g_n\to g$, then
            \begin{equation*}
                \mathcal{G}^0(\nu_{\mathcal{T}}^{g_n})\to\mathcal{G}^0(\nu_{\mathcal{T}}^g).
            \end{equation*}
        \item If $\varphi_\varepsilon\in\widetilde{\mathcal U}^N$ and
            $\varphi\in\mathcal U^N$ satisfy $\varphi_\varepsilon\Rightarrow\varphi$ as $\varepsilon\to0$, then
            \begin{equation*}
                \mathcal{G}^\varepsilon(\varepsilon N^{\varepsilon^{-1}\varphi_\varepsilon})\Rightarrow \mathcal{G}^0(\nu_{\mathcal{T}}^\varphi)
            \end{equation*}
            in distribution as $\varepsilon\to 0$.
    \end{enumerate}
\end{condition}

For $\phi\in\mathbb{U}$, define $\mathbb{S}_\phi=\{g\in\mathbb{S} :\phi=\mathcal{G}^0(\nu_{\mathcal{T}}^g)\}$. Let $I:\mathbb{U}\to [0,\infty]$ be defined by
\begin{equation}\label{def:rate_fun}
    I(\phi) = \inf\{L_{\mathcal{T}}(g) : g\in\mathbb{S}_\phi\},\quad\phi\in\mathbb{U}.
\end{equation}
By convention, $I(\phi) = \infty$ if $\mathbb{S}_\phi = \varnothing$.

The following restricted-control form of the variational criterion is \cite[Theorem~2.4]{1259.60065}, based on \cite[Theorem~4.2]{1231.60018}.

\begin{theorem}\label{thm:ldp-abstract}
    For $\varepsilon > 0$, let $Z^\varepsilon$ be defined by $Z^\varepsilon = \mathcal{G}^\varepsilon(\varepsilon N^{\varepsilon^{-1}})$, and suppose that Condition \ref{cond:ldp} holds. Then $I$ defined as in \eqref{def:rate_fun} is a good rate function on $\mathbb{U}$ and the family $\{Z^\varepsilon\}_{\varepsilon>0}$ satisfies a large deviation principle with rate function $I$.
\end{theorem}

\subsection{Main result of LDP}
We apply this criterion to \eqref{eqn:spe} under the following assumptions on the noise coefficient.

\begin{assumption}\label{assump:LDP}
    For $i=1,2$, assume that $G_i:V\times E\to V_i$ is a measurable mapping. Suppose that there exist nonnegative functions $L\in L^2(\nu)\cap\mathscr H^2$ and $\tilde L\in L^2(\nu)\cap\mathscr H^1$ such that
    \begin{enumerate}
        \item (Lipschitz) For $\bm u=(\bm v,T),\hat{\bm u}=(\hat{\bm v},\hat{T})\in V$,
            \begin{equation}\label{asm:g_Lip-LDP}
                |G(\bm u,\xi)-G(\hat{\bm u},\xi)|_H \leq L(\xi)|\bm u-\hat{\bm u}|_H.
            \end{equation}
        \item (Linear growth in $H$) For $\bm u=(\bm v,T)\in V$,
            \begin{equation}\label{asm:g_growth-H-LDP}
                \begin{aligned}
                    & |G_1(\bm v,T,\xi)|_{H_1}\leq L(\xi)(1+|\bm v|_{H_1}), \\
                    & |G_2(\bm v,T,\xi)|_{H_2}\leq L(\xi)(1+|T|_{H_2}).
                \end{aligned}
            \end{equation}
        \item (Linear growth in $V$) For $\bm u=(\bm v,T)\in V$,
            \begin{equation}\label{asm:g_growth-V-LDP}
                \begin{aligned}
                    & \|G_1(\bm v,T,\xi)\|_{V_1}\leq L(\xi)(1+\|\bm v\|_{V_1}), \\
                    & \|G_2(\bm v,T,\xi)\|_{V_2}\leq L(\xi)(1+\|T\|_{V_2}).
                \end{aligned}
            \end{equation}
        \item (Linear growth in $L^6$) For $\bm u=(\bm v,T)\in V$,
            \begin{equation}\label{asm:g_growth-L6-LDP}
                |\widetilde{G_1}(\bm v,T,\xi)|_6 \leq\tilde L(\xi)(1+|\widetilde{\bm v}|_6).
            \end{equation}
        \item (Joint barotropic and vertical energy growth) For $\bm v\in V_1$, set
            \begin{equation*}
                Y(\bm v)=h|\nabla\overline{\bm v}|_{L^2(\mathcal D)}^2 +|\partial_z\bm v|_2^2.
            \end{equation*}
            For $\bm u=(\bm v,T)\in V$ and $\nu$-almost every $\xi$, assume
            \begin{equation}\label{asm:g_joint-energy-LDP}
                \begin{split}
                    \big|h\langle\nabla\overline{G_1}(\bm u,\xi), \nabla\overline{\bm v}\rangle_{L^2(\mathcal D)} & +\langle\partial_zG_1(\bm u,\xi),\partial_z\bm v\rangle_{L^2(\mathcal O;\mathbb R^2)}\big| \\
                    & \leq L(\xi)(1+Y(\bm v)).
                \end{split}
            \end{equation}
    \end{enumerate}
    Here $L^2(\nu)$ denotes the space of Borel measurable functions $h:E\to\mathbb R$, identified up to $\nu$-almost everywhere equality, such that
    \begin{equation}\label{cond:L1-nu}
        \int_E |h(\xi)|^2\nu(\mathrm{d}\xi)<\infty.
    \end{equation}
    For $1\leq n<\infty$, the class $\mathscr H^n$ consists of Borel measurable functions $h:E\to\mathbb R$ for which there exists $\delta>0$ such that, for every $\Gamma\in\mathcal B(E)$ with $\nu(\Gamma)<\infty$,
    \begin{equation}\label{cond:exp-int}
        \int_\Gamma e^{\delta |h(\xi)|^n}\nu(\mathrm{d}\xi)<\infty.
    \end{equation}
\end{assumption}
\begin{remark}\label{rk:ldp-assumptions-suffice}
    Assumption \ref{assump:LDP} implies all noise hypotheses of Theorem \ref{thm:global}. Both $L$ and $\tilde L$ belong to $L^2(\nu)\cap\mathscr H^1$. If $h\geq0$ belongs to this class, then $h\in L^p(\nu)$ for every finite $p\geq2$: on $\{h\leq1\}$ use $h^p\leq h^2$, and on the finite-measure set $\{h>1\}$ use $h^p\leq C_{p,\delta}e^{\delta h}$. Apply this to $L$ and $\tilde L$ and integrate the pointwise bounds.
\end{remark}

\begin{example}\label{ex:coupled-noise}
    Let $E=\mathbb R$, let $\nu$ be the standard Gaussian probability measure, and fix $\psi\in V_1$ with $\widetilde\psi\ne0$ and $0\ne\chi\in H_2$. The coefficients
    \begin{equation*}
        G_1(\bm v,T,\xi)=\xi\big[\bm v+ \tanh(\langle T,\chi\rangle_{H_2})\psi\big], \qquad G_2(\bm v,T,\xi)=\xi T
    \end{equation*}
    satisfy Assumption \ref{assump:LDP} with $L=\tilde L=C(1+|\xi|)$, for sufficiently large $C$. Indeed, $\tanh$ is bounded and Lipschitz, $V\hookrightarrow L^6$, and the left side of \eqref{asm:g_joint-energy-LDP} is at most $C|\xi|(Y(\bm v)+Y(\bm v)^{1/2})$. The Gaussian measure admits the required quadratic exponential moment for sufficiently small $\delta>0$.
\end{example}

\begin{remark}\label{rk:est-int-L}
    By \cite[Lemma~3.4]{1259.60065}, the assumption $L\in L^2(\nu)\cap\mathscr H^2$ implies that, for every $N\in\mathbb N$,
    \begin{align}
        & \sup_{g\in S^N}\int_{E_{\mathcal{T}}}L(\xi)|g(t,\xi)-1|\nu(\mathrm{d}\xi)\,\mathrm{d}t<\infty,\label{est:int-L-1} \\
        & \sup_{g\in S^N}\int_{E_{\mathcal{T}}}L^2(\xi)(g(t,\xi)+1)\nu(\mathrm{d}\xi)\,\mathrm{d}t<\infty.\label{est:int-L-2}
    \end{align}
    The same entropy truncation argument, applied on measurable time sets, gives: for $A_n\subset [0,\mathcal{T}]$ satisfying $\mathrm{Leb}_{[0,\mathcal{T}]}(A_n)\to 0$ as $n\to\infty$,
    \begin{equation}\label{est:int-L-3}
        \lim_{n\to\infty}\sup_{g\in S^N}\int_{A_n}\int_EL(\xi)|g(t,\xi)-1|\nu(\mathrm{d}\xi)\,\mathrm{d}t=0,
    \end{equation}
    where $\mathrm{Leb}_{[0,\mathcal{T}]}$ denotes the Lebesgue measure on $[0,\mathcal{T}]$. For $\tilde L$, choose $\delta>0$ as in \eqref{cond:exp-int} with $n=1$. The entropy inequality
    \begin{equation*}
        a|r-1|\leq l(r)+e^a-1-a,\qquad a,r\geq0,
    \end{equation*}
    gives
    \begin{equation*}
        \sup_{g\in S^N}\int_0^{\mathcal T}\int_E \tilde L(\xi)|g(t,\xi)-1|\,\nu(\mathrm{d}\xi)\,\mathrm{d}t \leq\frac{1}{\delta}\left[N+\mathcal T\int_E \big(e^{\delta\tilde L}-1-\delta\tilde L\big)\,\mathrm{d}\nu\right] <\infty.
    \end{equation*}
    The integral on the right is finite by $\tilde L\in L^2(\nu)\cap\mathscr H^1$, splitting $E$ into $\{\tilde L\leq1\}$ and $\{\tilde L>1\}$.
\end{remark}
\begin{theorem}\label{thm:LDP}
    Assume that Assumption \ref{assump:LDP} holds. For every $\varepsilon>0$, let $\bm u^\varepsilon$ be the unique global pathwise solution of \eqref{eqn:spe} on the stochastic basis $(\bar{\mathbb{M}},\bar{\mathcal{F}},\bar{\mathbb F},\bar{\mathbb{P}})$ by Theorem \ref{thm:global}. Then the family $\{\bm u^\varepsilon\}_{\varepsilon>0}$ satisfies the large deviation principle on $\Upsilon_{\mathcal{T}}(V)$ as $\varepsilon\to0$, with a good rate function $I$ defined by
    \begin{equation}
        I(\phi)=\inf_{g\in\mathbb{S},\bm{u}^g=\phi}L_{\mathcal{T}}(g),\quad \phi\in \Upsilon_{\mathcal{T}}(V),
    \end{equation}
    where $\bm{u}^g=(\bm v^g,T^g)$ is the unique solution to the following skeleton equation
    \begin{equation}\label{eqn:skeleton}
        \begin{cases}
            \displaystyle\frac{\mathrm{d}\bm{u}^g(t)}{\mathrm{d}t}+A\bm{u}^g(t)+B(\bm{u}^g(t))+R(\bm{u}^g(t))=\int_EG(\bm{u}^g(t),\xi)(g(t,\xi)-1)\nu(\mathrm{d}\xi), \\
            \bm{u}^g(0)=\bm u_0=(\bm v_0,T_0)\in V.
        \end{cases}
    \end{equation}
\end{theorem}

\begin{proof}
    For the fixed deterministic initial datum, Theorem \ref{thm:global} and Definition \ref{def:pathwise-sol} give Borel maps $\mathcal G^\varepsilon:\mathbb M\to\Upsilon_{\mathcal T}(V)$ with
    \begin{equation}
        \bm u^\varepsilon=\mathcal G^\varepsilon(\varepsilon N^{\varepsilon^{-1}}).
    \end{equation}
    The maps are nonanticipating up to null sets. Equality of joint laws transfers the solution identity to any Poisson random measure with intensity $\varepsilon^{-1}\mathrm{d}t\otimes\nu$, as in \cite[proof of Theorem~5.1 and Corollary~5.1]{1519.35220}. By Lemma~\ref{lem:wp-skeleton}, define $\mathcal G^0:\mathbb M\to\Upsilon_{\mathcal T}(V)$ on finite-entropy intensities by
    \begin{equation}
        \mathcal{G}^0(\nu_{\mathcal{T}}^g)=\bm u^g,
    \end{equation}
    where $\nu_{\mathcal{T}}^g$ is defined in \eqref{def:nu_T^g} and controls equal $\nu_{\mathcal T}$-almost everywhere are identified. The set $\mathbb M_0=\bigcup_{N\geq1}\{\nu_{\mathcal T}^g:g\in S^N\}$ is Borel, since each set in the union is compact. Define $\mathcal G^0$ to be the zero path on $\mathbb M\setminus\mathbb M_0$. Lemma \ref{lem:cond-1} proves continuity on each entropy level, hence Borel measurability of this extension. Condition \ref{cond:ldp} is verified in Lemmas \ref{lem:cond-1} and \ref{lem:cond-2}. Theorem \ref{thm:ldp-abstract} then gives the conclusion.
\end{proof}

\subsection{Global well-posedness of the skeleton equation}\label{wp-skeleton}
We first establish the well-posedness of the skeleton equation \eqref{eqn:skeleton}.

\begin{lemma}\label{lem:wp-skeleton}
    Let Assumption \ref{assump:LDP} hold and $N\in\mathbb N$. Then for each $g\in S^N$, there exists a unique strong solution $\bm{u}^g\in C([0,\mathcal{T}];V)\cap L^2([0,\mathcal{T}];\mathscr D(A))$ to the skeleton equation \eqref{eqn:skeleton}, which depends continuously in this space on the initial datum in $V$. Moreover, for any $C_0>0$, there exists a constant $C_N>0$, depending also on $C_0$, such that for every $g\in S^N$ and $\bm u_0\in V$ such that $\|\bm u_0\|_V\leq C_0$, the following estimate holds:
    \begin{equation}\label{est:skeleton-bound}
        \sup_{t\in[0,\mathcal{T}]}\|\bm{u}^g(t)\|_V^2+\int_0^{\mathcal{T}}\|\bm{u}^g(t)\|_{\mathscr D(A)}^2\,\mathrm{d}t\leq C_N.
    \end{equation}
\end{lemma}

\begin{proof}
    We use the barotropic--baroclinic estimates of \cite[Sections~3.2--3.3]{1151.35074}. Write $\bm u=\bm u^g=(\bm v,T)$ and
    \begin{equation*}
        F_g(\bm u)(t)=\int_EG(\bm u,\xi)(g(t,\xi)-1)\,\nu(\mathrm{d}\xi), \qquad a_g(t)=\int_E L(\xi)|g(t,\xi)-1|\,\nu(\mathrm{d}\xi).
    \end{equation*}
    Define $\widetilde a_g$ with $\widetilde L$ in place of $L$. Remark \ref{rk:est-int-L} gives
    \begin{equation}\label{def:J}
        \sup_{g\in S^N}\int_0^{\mathcal T}(a_g+\widetilde a_g)\,\mathrm{d}t \leq C_N.
    \end{equation}
    All constants below are uniform for $g\in S^N$ and $\|\bm u_0\|_V\leq C_0$.

    \emph{Local solutions.} The growth and Lipschitz assumptions imply
    \begin{equation}\label{est:skeleton-force}
        \|F_g(\bm u)\|_V\leq C a_g(1+\|\bm u\|_V),\qquad |F_g(\bm u)-F_g(\hat{\bm u})|_H \leq C a_g|\bm u-\hat{\bm u}|_H.
    \end{equation}
    The local spectral construction applies with this $L^1V$ forcing. Indeed, the spectral theorem and Minkowski's inequality give, on any interval $I=[s,t]$,
    \begin{equation}\label{est:skeleton-linear}
        \left\|\int_s^r e^{-(r-\ell)A}f(\ell)\,\mathrm{d}\ell\right\|_{C(I;V)} +\left\|A\int_s^r e^{-(r-\ell)A}f(\ell)\,\mathrm{d}\ell\right\|_{L^2(I;H)} \leq C\|f\|_{L^1(I;V)}.
    \end{equation}
    For the second term use $\int_0^\infty|Ae^{-rA}v|_H^2\,\mathrm{d}r=\|v\|_V^2/2$. Put $\bm u^*(t)=e^{-tA}\bm u_0$ and $w_n=\bm u_n-P_n\bm u^*$ in the spectral approximation. Before $\|w_n\|_V$ reaches a fixed sufficiently small $\rho$, \eqref{est:B-H-symmetric} absorbs $B(w_n,w_n)$ into the dissipation. The mixed terms and \eqref{est:skeleton-force} give
    \begin{equation*}
        \sup_{r\leq t}\|w_n(r)\|_V^2+\int_0^t|Aw_n|_H^2\,\mathrm{d}s \leq C b_g(t)\exp(Cb_g(t)+Ct),
    \end{equation*}
    where
    \begin{equation*}
        b_g(t)=\int_0^t(1+\|\bm u^*\|_V^2) (|A\bm u^*|_H^2+a_g)\,\mathrm{d}s.
    \end{equation*}
    Since $b_g(t)\to0$ as $t\downarrow0$, the approximations remain inside that ball on a common interval. Truncation outside the ball therefore gives a local solution in $CV\cap L^2\mathscr D(A)$. In passing to the limit, compactness gives strong $L^2V$ convergence; the second inequality in \eqref{est:skeleton-force}, the local $V$ bound, and the absolute continuity of $\int a_g$ give convergence of the forcing in $L^1H$. The limiting forcing belongs to $L^1V$, and \eqref{est:skeleton-linear} gives continuity in $V$. The same construction at any initial time gives a maximal solution; see \cite[Section~5.1.4]{1370.60103} for the continuation principle. The estimates that follow hold on every compact subinterval of its existence interval. The energy pairings are justified by spectral approximation for forcing in $L^1V$.

    \emph{Low-order and $L^6$ estimates.} Testing with $\bm u$, using the cancellation of $B$ and \eqref{est:skeleton-force}, gives
    \begin{equation*}
        \frac{\mathrm{d}}{\mathrm{d}t}|\bm u|_H^2+c\|\bm u\|_V^2 \leq C(1+a_g)(1+|\bm u|_H^2).
    \end{equation*}
    Hence, with $q=\|\bm u\|_V^2$,
    \begin{equation}\label{est:skeleton-low}
        \sup_{s\leq t}|\bm u(s)|_H^2+\int_0^t q(s)\,\mathrm{d}s\leq C_N.
    \end{equation}
    Put $\mathcal X=\int_{\mathcal O}|\widetilde{\bm v}|^4 |\nabla_3\widetilde{\bm v}|^2\,\mathrm{d}x\,\mathrm{d}y\,\mathrm{d}z$. The deterministic estimates \eqref{est:L6-determined} and \eqref{asm:g_growth-L6-LDP} give
    \begin{equation*}
        \frac{\mathrm{d}}{\mathrm{d}t}|\widetilde{\bm v}|_6^6+c\mathcal X \leq C\big[(1+|\bm u|_H^2)q+\widetilde a_g\big] (1+|\widetilde{\bm v}|_6^6).
    \end{equation*}
    For example, the control contribution is at most $C\widetilde a_g(1+|\widetilde{\bm v}|_6^6)$. Gr\"onwall's inequality and \eqref{est:skeleton-low} imply
    \begin{equation}\label{est:skeleton-L6}
        \sup_{s\leq t}|\widetilde{\bm v}(s)|_6^6 +\int_0^t\mathcal X(s)\,\mathrm{d}s\leq C_N.
    \end{equation}

    \emph{The barotropic gradient and vertical derivative.} Let $S=-P_\sigma\Delta$ be the barotropic Stokes operator and $A_h\bm v=S\overline{\bm v}-\Delta\widetilde{\bm v}$, as in the estimates preceding \eqref{est:barotropic-local}, and set
    \begin{equation*}
        Y_b=h|\nabla\overline{\bm v}|_{L^2(\mathcal D)}^2,\qquad Y_z=|\bm v_z|_2^2,\qquad Y=Y_b+Y_z, \qquad D_z=|\nabla_3\bm v_z|_2^2.
    \end{equation*}
    The estimates preceding \eqref{est:barotropic-local} give
    \begin{equation*}
        Y_b'+c h|S\overline{\bm v}|_{L^2(\mathcal D)}^2 \leq C(1+|\bm v|_2^2Y_b)Y_b+C(1+q+\mathcal X) +2h\langle\nabla\overline{F_{g,1}}, \nabla\overline{\bm v}\rangle_{L^2(\mathcal D)}.
    \end{equation*}
    For the vertical equation, separate $\bm v=\overline{\bm v} +\widetilde{\bm v}$ in its stretching term. Since $\operatorname{div}\overline{\bm v}=0$, the barotropic contribution obeys
    \begin{equation}\label{est:skeleton-vertical-bar}
        \begin{split}
            \left|\int_{\mathcal O}(\bm v_z\cdot\nabla) \overline{\bm v}\cdot\bm v_z\,\mathrm{d}x\,\mathrm{d}y\,\mathrm{d}z\right| & \leq C|\nabla\overline{\bm v}|_{L^2(\mathcal D)} |\bm v_z|_2\big(|\bm v_z|_2+|\nabla\bm v_z|_2\big) \\
            & \leq\delta D_z+C_\delta(1+Y_b)Y_z.
        \end{split}
    \end{equation}
    Here the two-dimensional interpolation inequality is applied to $(\int_{-h}^0|\bm v_z|^2\,\mathrm{d}z)^{1/2}$. Integration by parts in the baroclinic contribution yields $\delta D_z+C_\delta(1+|\widetilde{\bm v}|_6^4)Y_z$. The temperature term is bounded by $\delta D_z+C_\delta|T|_2^2$. Consequently,
    \begin{equation*}
        Y_z'+cD_z\leq C(1+Y_b+|\widetilde{\bm v}|_6^4)Y_z +C|T|_2^2+2\langle\partial_zF_{g,1},\bm v_z\rangle_{L^2(\mathcal O;\mathbb R^2)}.
    \end{equation*}
    The sum of the two control terms is at most $2a_g(1+Y)$ by \eqref{asm:g_joint-energy-LDP}. Since $Y_b\leq Cq$, addition gives
    \begin{equation}\label{est:skeleton-joint}
        Y'+c\big(h|S\overline{\bm v}|_{L^2(\mathcal D)}^2+D_z\big) \leq C\big[1+(1+|\bm u|_H^2)q +|\widetilde{\bm v}|_6^4+a_g\big](1+Y)+C\mathcal X.
    \end{equation}
    Every coefficient on the right has a uniformly bounded time integral by \eqref{def:J}, \eqref{est:skeleton-low} and \eqref{est:skeleton-L6}. Thus
    \begin{equation}\label{est:H1-skeleton-1}
        \sup_{s\leq t}Y(s)+\int_0^t \big(h|S\overline{\bm v}|_{L^2(\mathcal D)}^2+D_z\big)\,\mathrm{d}s \leq C_N.
    \end{equation}
    In particular, $\sup_{s\leq t}|\bm v(s)|_6\leq C_N$.

    \emph{The remaining gradients.} Set $X=|\nabla\bm v|_2^2$ and $D_h=|A_h\bm v|_2^2+|\nabla\bm v_z|_2^2$. The form pairing and \eqref{asm:g_growth-V-LDP} give
    \begin{equation*}
        |\langle F_{g,1},A_h\bm v\rangle_{H_1}| =|\langle\nabla F_{g,1},\nabla\bm v\rangle| \leq C a_g(1+X+Y_z).
    \end{equation*}
    Using \eqref{est:nabla-v-determined}, we obtain
    \begin{equation}\label{est:skeleton-horizontal}
        X'+cD_h\leq C(R_v+a_g)X+C(q+1)+Ca_g(1+Y_z), \qquad R_v=1+|\bm v|_6^4+(1+Y_z)(Y_z+D_z).
    \end{equation}
    Estimate~\eqref{est:H1-skeleton-1} bounds $\int_0^tR_v$. Gr\"onwall's inequality, followed by $A_1=A_h-\partial_{zz}$, yields
    \begin{equation}\label{est:skeleton-velocity}
        \sup_{s\leq t}\|\bm v(s)\|_{V_1}^2 +\int_0^t|A_1\bm v|_2^2\,\mathrm{d}s\leq C_N.
    \end{equation}
    Finally, \eqref{est:T-V-determined} and the $V_2$ form pairing give
    \begin{equation*}
        \frac{\mathrm{d}}{\mathrm{d}t}\|T\|_{V_2}^2+c|A_2T|_2^2 \leq C\big[1+|\bm v|_6^4+\|\bm v\|_{V_1}^2|A_1\bm v|_2^2+a_g\big] (1+\|T\|_{V_2}^2).
    \end{equation*}
    Here $a_2$ includes the Robin boundary term. By \eqref{est:skeleton-velocity}, the coefficient on the right-hand side is integrable. Hence
    \begin{equation}\label{est:H1-skeleton-4}
        \sup_{s\leq t}\|T(s)\|_{V_2}^2 +\int_0^t|A_2T|_2^2\,\mathrm{d}s\leq C_N.
    \end{equation}

    \emph{Continuation and uniqueness.} All preceding bounds are independent of the endpoint $t$ of the local interval. If its maximal endpoint $\mathcal T^*$ were at most $\mathcal T$, they would give $\bm u\in L^\infty(0,\mathcal T^*;V) \cap L^2(0,\mathcal T^*;\mathscr D(A))$ and $F_g(\bm u)\in L^1(0,\mathcal T^*;V)$. Moreover, \eqref{est:B-H-symmetric} implies $B(\bm u)+R(\bm u)\in L^2(0,\mathcal T^*;H)$. The variation-of-constants formula and \eqref{est:skeleton-linear} therefore provide a limit in $V$ at $\mathcal T^*$. Applying the local existence result with this limit as initial datum extends the solution when $\mathcal T^*<\mathcal T$; at $\mathcal T$ the limit defines the endpoint value. This proves global existence on $[0,\mathcal T]$ and \eqref{est:skeleton-bound}.

    For two solutions with the same control put $d=\bm u_1-\bm u_2$. By \eqref{id:B-skew}, \eqref{est:B-diff_prod} and \eqref{est:skeleton-force},
    \begin{equation}\label{est:v^g-minus}
        \frac{\mathrm{d}}{\mathrm{d}t}|d|_H^2+\|d\|_V^2 \leq C\big(1+a_g+\|\bm u_2\|_V^4 +\|\bm u_2\|_V^2|A\bm u_2|_H^2\big)|d|_H^2.
    \end{equation}
    The coefficient has a uniformly bounded integral on bounded sets of initial data in $V$. Gr\"onwall's inequality proves uniqueness and continuous dependence in $CH\cap L^2V$. If $\bm u_{0,n}\to\bm u_0$ in $V$, write $d_n=\bm u_n-\bm u$. Testing the difference equation with $Ad_n$, the symmetric estimate \eqref{est:B-H-symmetric} and Young's inequality give
    \begin{equation*}
        \|d_n\|_{CV}^2+\|Ad_n\|_{L^2H}^2 \leq C_N\left(\|d_n(0)\|_V^2+ \int_0^{\mathcal T}a_g\|d_n\|_V\,\mathrm{d}t\right).
    \end{equation*}
    Here the integrated nonlinear coefficient is bounded by $C\int(1+\|\bm u_n\|_V^2|A\bm u_n|_H^2+ \|\bm u\|_V^2|A\bm u|_H^2)$, and the forcing difference has $V$ norm at most $C_Na_g$. Since $d_n\to0$ in $L^2V$ and is bounded in $L^\infty V$, the last integral tends to zero by absolute continuity of $\int a_g$. This proves continuity in $CV\cap L^2\mathscr D(A)$.
\end{proof}

\subsection{The first continuity lemma}\label{sec:cond-1}
The first part of Condition~\ref{cond:ldp} concerns continuity of the skeleton solution map
\begin{equation*}
    \mathcal{G}^0(\nu_{\mathcal{T}}^g)=\bm u^g,
\end{equation*}
with respect to the control $g\in S^N$.

\begin{lemma}\label{lem:cond-1}
    Let Assumption~\ref{assump:LDP} hold and fix $N\in\mathbb N$. If $g_n\to g$ in $S^N$, then
    \begin{equation*}
        \mathcal{G}^0(\nu_{\mathcal{T}}^{g_n})\to\mathcal{G}^0(\nu_{\mathcal{T}}^g), \quad \text{i.e., }\bm u^{g_n}\to\bm u^g\ \text{in }\Upsilon_{\mathcal{T}}(V).
    \end{equation*}
\end{lemma}

The proof uses fractional Sobolev compactness and a time-continuity criterion. Let $K$ be a separable Hilbert space. For $p>1$ and $\alpha\in(0,1)$, the space $W^{\alpha,p}([0,\mathcal T];K)$ consists of $u\in L^p([0,\mathcal T];K)$ such that
\begin{equation*}
    \int_0^{\mathcal{T}}\int_0^{\mathcal{T}}\frac{\|u(t)-u(s)\|_K^p}{|t-s|^{1+\alpha p}}\,\mathrm{d}t\,\mathrm{d}s<\infty,
\end{equation*}
endowed with the norm
\begin{equation*}
    \|u\|_{W^{\alpha,p}([0,\mathcal{T}]; K)}^p=\int_0^{\mathcal{T}}\|u(t)\|_K^p\,\mathrm{d}t+\int_0^{\mathcal{T}}\int_0^{\mathcal{T}}\frac{\|u(t)-u(s)\|_K^p}{|t-s|^{1+\alpha p}}\,\mathrm{d}t\,\mathrm{d}s.
\end{equation*}
We use the following compact embedding \cite[Theorem~2.1, p.~370]{0831.60072}.
\begin{lemma}\label{lem:compactness-fractional-sobolev}
    Let $B_0\hookrightarrow B\hookrightarrow B_1$ be continuous embeddings of Banach spaces, with $B_0,B_1$ reflexive and $B_0\hookrightarrow B$ compact. Let $p\in(1,\infty)$ and $\alpha\in(0,1)$ be given. Let $\Lambda$ be the space
    \begin{equation*}
        \Lambda:=L^p([0,\mathcal{T}]; B_0)\cap W^{\alpha,p}([0,\mathcal{T}]; B_1),
    \end{equation*}
    endowed with the natural norm. Then the embedding $\Lambda\hookrightarrow L^p([0,\mathcal{T}];B)$ is compact.
\end{lemma}

For continuity in time, we use the Lions--Magenes lemma \cite[Chapter~III, Lemma~1.2]{0981.35001}.

\begin{lemma}[Lions--Magenes]\label{lem:lions-magenes}
    Let $V\subset\subset H\cong H'\subset V'$ be a Hilbert triple, with $V'$ the dual of $V$. If $u\in L^2(0,\mathcal T;V)$ and $\frac{\mathrm{d}u}{\mathrm{d}t}\in L^2(0,\mathcal T;V')$, then $u$ has a representative in $C([0,\mathcal T];H)$.
\end{lemma}

\begin{proof}[Proof of Lemma~\ref{lem:cond-1}]
    Write $u_n=\bm u^{g_n}$ and $u=\bm u^g$, and let $M$ be a common bound in \eqref{est:skeleton-bound}. Put $a_g(t)=\int_E L(\xi)|g(t,\xi)-1|\,\nu(\mathrm{d}\xi)$. By \eqref{est:int-L-1} and \eqref{est:int-L-3}, the functions $a_g$, $g\in S^N$, are bounded and uniformly integrable in $L^1(0,\mathcal T)$.

    \emph{Compactness and identification.} The control primitives $J_{5,n}(t)=\int_0^t\int_E G(u_n(s),\xi)(g_n(s,\xi)-1)\,\nu(\mathrm{d}\xi)\,\mathrm{d}s$ are uniformly bounded in $BV([0,\mathcal T];V)$ and in $L^\infty V$. Indeed, their total variation is at most $C_M\int_0^{\mathcal T}a_{g_n}$. For a bounded $V$-valued function $J$ of bounded variation,
    \begin{equation*}
        \int_0^{\mathcal T-h}\|J(t+h)-J(t)\|_V^2\,\mathrm{d}t \leq 2\|J\|_{L^\infty V}\,h\,\operatorname{Var}(J).
    \end{equation*}
    Consequently,
    \begin{equation}\label{est:j5}
        \sup_n\|J_{5,n}\|_{W^{\alpha,2}(0,\mathcal T;V)}<\infty, \qquad 0<\alpha<\tfrac12.
    \end{equation}
    Estimate~\eqref{est:B-H-symmetric} bounds $Au_n+B(u_n)+R(u_n)$ in $L^2(0,\mathcal T;H)$. Lemma~\ref{lem:compactness-fractional-sobolev} therefore gives, along a subsequence, $u_n\to\hat u$ in $L^2(0,\mathcal T;V)$, with $\hat u\in L^\infty V\cap L^2\mathscr D(A)$. The Lipschitz assumption and uniform integrability give
    \begin{equation*}
        \int_0^{\mathcal T}\!\int_E |G(u_n,\xi)-G(\hat u,\xi)|_H|g_n-1|\,\nu(\mathrm{d}\xi)\,\mathrm{d}t \leq\int_0^{\mathcal T}a_{g_n}|u_n-\hat u|_H\,\mathrm{d}t\to0.
    \end{equation*}
    For the frozen kernel $G(\hat u(t),\xi)$, apply \cite[Lemma 3.11]{1259.60065} to each scalar test. Its hypotheses follow from $|G(\hat u,\xi)|_H\leq C_M L(\xi)$: quadratic exponential integrability implies every linear exponential integrability on sets of finite measure. Thus the signed control integrals converge. The remaining terms pass to the limit using \eqref{est:B-VDA} and strong $L^2V$ convergence. The resulting weak equation has initial datum $\bm u_0$.

    If $F\in L^1(0,\mathcal T;V)$, the spectral theorem and Minkowski's inequality give
    \begin{equation}\label{est:linear-L1V}
        \left\|\int_0^t e^{-(t-s)A}F(s)\,\mathrm{d}s\right\|_{C V} +\left\|A\int_0^t e^{-(t-s)A}F(s)\,\mathrm{d}s\right\|_{L^2 H} \leq C\|F\|_{L^1 V}.
    \end{equation}
    Indeed, $\int_0^\infty|Ae^{-tA}f|_H^2\,\mathrm{d}t=\|f\|_V^2/2$. The $L^2H$ part of the linear equation is continuous in $V$ by its energy estimate and Lemma~\ref{lem:lions-magenes}, applied to the pair $(\mathscr D(A),V)$. Together with \eqref{est:linear-L1V}, this gives $\hat u\in CV$. Uniqueness, proved by \eqref{est:v^g-minus}, identifies $\hat u=u$. Hence the whole sequence converges to $u$ in $L^2V$.

    \emph{A linear corrector.} Define
    \begin{equation}\label{eq:control-corrector}
        h_n(t)=\int_E G(u(t),\xi)(g_n-g)(t,\xi)\,\nu(\mathrm{d}\xi), \qquad z_n'+Az_n=h_n,\quad z_n(0)=0.
    \end{equation}
    We claim that
    \begin{equation}\label{est:control-corrector}
        z_n\to0\quad\text{in }CV\cap L^2\mathscr D(A).
    \end{equation}
    Let $P_m$ be the spectral projections of $A$. For each fixed $m$, \cite[Lemma 3.11]{1259.60065} gives pointwise convergence to zero of $\int_0^tP_mh_n$. Uniform integrability of $a_{g_n}+a_g$ makes these primitives equicontinuous. They converge uniformly, and the finite-dimensional variation-of-constants formula implies $P_mz_n\to0$ in $CV\cap L^2\mathscr D(A)$.

    To control the remaining modes, put $q_m(t,\xi)=\|(I-P_m)G(u(t),\xi)\|_V$. Then $q_m\to0$ almost everywhere and $q_m\leq C_M L$. The entropy bound implies
    \begin{equation}\label{est:entropy-spectral-tail}
        \lim_{m\to\infty}\sup_{g\in S^N} \int_{E_{\mathcal T}}q_m|g-1|\,\mathrm{d}\nu_{\mathcal T}=0.
    \end{equation}
    Fix $\beta>1$. On $\{|g-1|\leq\beta\}$, the inequality $|g-1|^2\leq C_\beta l(g)$ and Cauchy--Schwarz bound the integral by $(C_\beta N)^{1/2}\|q_m\|_{L^2(\nu_{\mathcal T})}$, which tends to zero as $m\to\infty$. On $\{|g-1|>\beta,q_m<1\}$ use $|g-1|\leq c_\beta l(g)$, where $c_\beta\to0$. On $\{q_m\geq1\}$, Young's entropy inequality bounds the integral by $N/a+C_a\int_{\{q_m\geq1\}}e^{a q_m}\,\mathrm{d}\nu_{\mathcal T}$. For fixed $a$, this last integral tends to zero by dominated convergence: $\{q_m\geq1\}\subset\{C_M L\geq1\}$ has finite measure and $e^{a C_M L}$ is integrable there. Let $m\to\infty$, then $\beta,a\to\infty$. This proves \eqref{est:entropy-spectral-tail}. Apply \eqref{est:linear-L1V} to $(I-P_m)h_n$ and then let $m\to\infty$ to conclude \eqref{est:control-corrector}.

    \emph{The strong energy estimate.} Set $w_n=u_n-u-z_n$ and $f_n=\int_E[G(u_n,\xi)-G(u,\xi)](g_n-1)\,\nu(\mathrm{d}\xi)$. Then
    \begin{equation}\label{eqn:v_n-minus-v^g}
        w_n'+Aw_n+B(u_n)-B(u)+R(w_n+z_n)=f_n,\qquad w_n(0)=0.
    \end{equation}
    The $V$-growth assumption gives $\|f_n(t)\|_V\leq C_Ma_{g_n}(t)$. Use the decomposition
    \begin{equation*}
        B(u_n)-B(u)=B(u_n,w_n)+B(w_n,u)+B(u_n,z_n)+B(z_n,u).
    \end{equation*}
    Estimate~\eqref{est:B-H-symmetric} gives, for $a,d\in\mathscr D(A)$,
    \begin{equation}\label{est:B-strong-difference}
        |(B(a,d),Ad)_H|+|(B(d,a),Ad)_H| \leq\varepsilon|Ad|_H^2 +C_\varepsilon\|a\|_V^2|Aa|_H^2\|d\|_V^2.
    \end{equation}
    For the control term use the form pairing instead:
    \begin{equation}\label{est:r_n-4}
        (f_n,Aw_n)_H=(f_n,w_n)_V \leq C_Ma_{g_n}\|w_n\|_V.
    \end{equation}
    The energy identity follows by spectral approximation, since its forcing lies in $L^2H+L^1V$. Thus
    \begin{equation}\label{est:rho_n}
        \frac{\mathrm{d}}{\mathrm{d}t}\|w_n\|_V^2+|Aw_n|_H^2 \leq C(1+k_n)\|w_n\|_V^2 +C_Ma_{g_n}\|w_n\|_V+C e_n,
    \end{equation}
    where
    \begin{equation*}
        k_n=\|u_n\|_V^2|Au_n|_H^2+\|u\|_V^2|Au|_H^2,\qquad e_n=|B(u_n,z_n)|_H^2+|B(z_n,u)|_H^2+\|z_n\|_V^2.
    \end{equation*}
    The integrals of $k_n$ are uniformly bounded. Moreover,
    \begin{equation*}
        \int_0^{\mathcal T}|B(u_n,z_n)|_H^2\,\mathrm{d}t \leq C\|u_n\|_{L^\infty V}\|z_n\|_{CV} \|Au_n\|_{L^2H}\|Az_n\|_{L^2H}\to0,
    \end{equation*}
    and the other terms in $e_n$ converge to zero in the same way. Since $w_n\to0$ in $L^2V$ and is bounded in $L^\infty V$, uniform integrability yields $\int a_{g_n}\|w_n\|_V\to0$. Gr\"onwall's inequality in \eqref{est:rho_n} now gives
    \begin{equation}\label{est:rho_n-gronwall}
        \|w_n\|_{CV}^2+\|Aw_n\|_{L^2H}^2 \leq C_M\left(\int_0^{\mathcal T}a_{g_n}\|w_n\|_V\,\mathrm{d}t +\int_0^{\mathcal T}e_n\,\mathrm{d}t\right)\to0.
    \end{equation}
    Together with \eqref{est:control-corrector}, this proves
    \begin{equation}\label{est:rho_n-H1}
        u_n\to u\quad\text{in }CV\cap L^2(0,\mathcal T;\mathscr D(A)).\qedhere
    \end{equation}
\end{proof}

\subsection{The second continuity lemma}\label{sec:cond-2}
For the second part of Condition~\ref{cond:ldp}, we compare the stochastic solution with the skeleton driven by the same control. The linear corrector follows the approach of \cite[Lemma~7.6]{1519.35220} and \cite[eqs.~(6.75)--(6.78)]{1370.60103}.

\begin{lemma}\label{lem:cond-2}
    Let Assumption~\ref{assump:LDP} hold. Let $\varphi_\varepsilon\in\widetilde{\mathcal U}^N$ converge in distribution on $S^N$ to $\varphi\in\mathcal U^N$. Then
    \begin{equation*}
        \mathcal G^\varepsilon(\varepsilon N^{\varepsilon^{-1}\varphi_\varepsilon}) \ \Rightarrow\ \mathcal G^0(\nu_{\mathcal T}^\varphi) \quad\text{in }\Upsilon_{\mathcal T}(V).
    \end{equation*}
\end{lemma}
\begin{proof}
    Throughout, the initial datum $\bm u_0\in V$ is fixed and deterministic. Write $F_g(v)(t)=\int_EG(v,\xi)(g(t,\xi)-1)\,\nu(\mathrm{d}\xi)$ and $a_g(t)=\int_E L(\xi)|g(t,\xi)-1|\,\nu(\mathrm{d}\xi)$. Remark~\ref{rk:est-int-L} gives finite constants $J_N,Q_N$ and a modulus $\omega_N(r)\downarrow0$ as $r\downarrow0$ such that
    \begin{equation}\label{eq:controlled-entropy-modulus}
        \begin{gathered}
            \sup_{g\in S^N}\int_0^{\mathcal T}a_g\,\mathrm{d}t\leq J_N, \qquad \sup_{g\in S^N}\int_0^{\mathcal T}\int_E L^2g\,\mathrm{d}\nu\,\mathrm{d}t \leq Q_N, \\
            \sup_{g\in S^N}\sup_{|D|\leq r}\int_D a_g(t)\,\mathrm{d}t \leq\omega_N(r),
        \end{gathered}
    \end{equation}
    where $D$ ranges over measurable time sets. These bounds apply pathwise to random controls and random time sets.

    \emph{The controlled equation.} Fix $\varepsilon$ and choose $m$ with $\varphi_\varepsilon\in\bar{\mathcal A}_{b,m}$. The bounded change of intensity in \cite[Theorem~6.1]{1519.35220} applies: on the marked Poisson space its density uses the predictable factor
    \begin{equation*}
        \vartheta(s,\xi,r)= \begin{cases}
            \varphi_\varepsilon(s,\xi)^{-1}, & r\leq\varepsilon^{-1}\varphi_\varepsilon(s,\xi), \\
            1, & \text{otherwise}.
        \end{cases}
    \end{equation*}
    It is bounded away from zero and infinity and differs from $1$ only on $[0,\mathcal T]\times K_m\times[0,m/\varepsilon]$, a set of finite intensity. The resulting probability measure is equivalent to $\bar{\mathbb P}$, and under this measure $N^{\varepsilon^{-1}\varphi_\varepsilon}$ is a Poisson random measure with intensity measure $\varepsilon^{-1}\mathrm{d}t\otimes\nu$. Thus, by the solution map and the change-of-measure argument of \cite[Lemma~7.1]{1519.35220}, $X^\varepsilon=\mathcal G^\varepsilon (\varepsilon N^{\varepsilon^{-1}\varphi_\varepsilon})$ solves
    \begin{equation}\label{eqn:X^epsilon}
        \begin{split}
            \mathrm{d}X^\varepsilon+[AX^\varepsilon+B(X^\varepsilon) +R(X^\varepsilon)]\,\mathrm{d}t & =F_{\varphi_\varepsilon}(X^\varepsilon)\,\mathrm{d}t +\varepsilon\int_EG(X^\varepsilon(t-),\xi) \,\widetilde N^{\varepsilon^{-1}\varphi_\varepsilon} (\mathrm{d}t,\mathrm{d}\xi), \\
            X^\varepsilon(0) & =\bm u_0.
        \end{split}
    \end{equation}
    Here the compensator under $\bar{\mathbb P}$ is $\varepsilon^{-1}\varphi_\varepsilon\nu\,\mathrm{d}t$. After localization in the energy of $X^\varepsilon$, predictable simple-integrand approximations converge in probability under either of the two equivalent measures. Their finite-sum identities are the same, and the compensator difference is $F_{\varphi_\varepsilon}(X^\varepsilon)\,\mathrm{d}t$. Removing the localization proves \eqref{eqn:X^epsilon}. In particular $X^\varepsilon$ has paths in $\Upsilon_{\mathcal T}(V)$.

    Set $U^\varepsilon=\bm u^{\varphi_\varepsilon}$. By Lemmas~\ref{lem:wp-skeleton} and~\ref{lem:cond-1}, this is a measurable process, and uniqueness on time restrictions makes it adapted. For a deterministic $K_N<\infty$,
    \begin{equation}\label{eq:comparison-skeleton-bound}
        \sup_{t\leq\mathcal T}\|U^\varepsilon(t)\|_V^2 +\int_0^{\mathcal T}|AU^\varepsilon|_H^2\,\mathrm{d}t\leq K_N \quad\text{almost surely}.
    \end{equation}

    \emph{The stopped corrector.} For a path $v$, write $\mathcal E_t(v)=\sup_{s\leq t}\|v(s)\|_V^2+ \int_0^t|Av(s)|_H^2\,\mathrm{d}s$. Fix $M>K_N+1$ and define
    \begin{equation}\label{def:tau-epsilon-M}
        \tau=\tau^{\varepsilon,M} =\inf\{t\in[0,\mathcal T]:\mathcal E_t(X^\varepsilon)\geq M\} \wedge\mathcal T,\qquad\inf\varnothing=\infty.
    \end{equation}
    Before $\tau$ the $V$ norm is bounded by $\sqrt M$, and the dissipation integral up to $\tau$ is at most $M$. Let $Y=Y^{\varepsilon,M}$ solve
    \begin{equation}\label{eqn:Langevin}
        \mathrm{d}Y+AY\,\mathrm{d}t =\varepsilon\int_E\mathbbm{1}_{\{s\leq\tau\}} G(X^\varepsilon(s-),\xi)\, \widetilde N^{\varepsilon^{-1}\varphi_\varepsilon}(\mathrm{d}s,\mathrm{d}\xi), \qquad Y(0)=0.
    \end{equation}
    The indicator is predictable. The linear spectral estimate used in Lemma~\ref{lem:strong-cadlag}, now with the controlled compensator, and \eqref{asm:g_growth-V-LDP} imply
    \begin{equation}\label{est:uniform-stopped}
        \begin{split}
            \bar{\mathbb E}\mathcal E_{\mathcal T}(Y) & \leq C\varepsilon\bar{\mathbb E}\int_0^{\mathcal T}\int_E \mathbbm{1}_{\{s\leq\tau\}}\|G(X^\varepsilon(s-),\xi)\|_V^2 \varphi_\varepsilon\,\mathrm{d}\nu\,\mathrm{d}s \\
            & \leq C\varepsilon(1+M)Q_N.
        \end{split}
    \end{equation}

    \emph{Comparison in $H$ and $V$.} Suppress $\varepsilon$ and put $w=X-U-Y$ on $[0,\tau]$. The stochastic integrals cancel, including at $\tau$, and
    \begin{equation}\label{eqn:Z-epsilon}
        w'+Aw+B(X)-B(U)+R(w+Y)=f,\qquad f=F_{\varphi_\varepsilon}(X)-F_{\varphi_\varepsilon}(U), \quad w(0)=0.
    \end{equation}
    In particular $w$ is continuous in $V$ on this interval. Write
    \begin{equation*}
        B(X)-B(U)=B(X,w)+B(w,U)+E_B, \qquad E_B=B(X,Y)+B(Y,U).
    \end{equation*}
    By \eqref{est:B-H-symmetric},
    \begin{equation}\label{eq:controlled-cross-errors}
        \int_0^\tau\big(|E_B|_H^2+|R(Y)|_H^2\big)\,\mathrm{d}t \leq C_{N,M}\mathcal E_{\mathcal T}(Y).
    \end{equation}
    For example, the first term is bounded using
    \begin{equation*}
        \int_0^\tau|B(X,Y)|_H^2\,\mathrm{d}t \leq C\|X\|_{L^\infty(0,\tau;V)}\|Y\|_{L^\infty(0,\tau;V)} \|AX\|_{L^2(0,\tau;H)}\|AY\|_{L^2(0,\tau;H)}.
    \end{equation*}
    All time norms of $X$ here exclude any terminal jump in the essential supremum. The bounds on $U$ give the same estimate for $B(Y,U)$. The pointwise $H$-Lipschitz assumption gives $|f|_H\leq a_{\varphi_\varepsilon}(|w|_H+|Y|_H)$. Testing \eqref{eqn:Z-epsilon} with $w$, using $b(X,w,w)=0$, \eqref{est:B-diff_prod}, and Gr\"onwall's inequality yields
    \begin{equation}\label{eq:controlled-H-comparison}
        \sup_{t\leq\tau}|w(t)|_H^2+\int_0^\tau\|w\|_V^2\,\mathrm{d}t \leq C_{N,M}\mathcal E_{\mathcal T}(Y).
    \end{equation}
    Indeed, \eqref{eq:controlled-entropy-modulus} and \eqref{eq:comparison-skeleton-bound} bound the integrated coefficient:
    \begin{equation*}
        \int_0^{\mathcal T}(1+a_{\varphi_\varepsilon} +\|U\|_V^4+\|U\|_V^2|AU|_H^2)\,\mathrm{d}t\leq C_N.
    \end{equation*}

    For the $V$ estimate use \eqref{est:B-strong-difference}. The growth assumption gives $\|f\|_V\leq C_{N,M}a_{\varphi_\varepsilon}$ almost everywhere before $\tau$, so that $(f,Aw)_H=(f,w)_V\leq C_{N,M}a_{\varphi_\varepsilon}\|w\|_V$. Spectral approximation justifies this energy pairing for forcing in $L^2H+L^1V$. The coefficient satisfies
    \begin{equation*}
        \int_0^\tau(1+\|X\|_V^2|AX|_H^2 +\|U\|_V^2|AU|_H^2)\,\mathrm{d}t \leq C_{N,M}.
    \end{equation*}
    Hence
    \begin{equation}\label{eq:controlled-V-comparison}
        \mathcal E_\tau(w)\leq C_{N,M}\left[ \mathcal E_{\mathcal T}(Y) +\int_0^\tau a_{\varphi_\varepsilon}\|w\|_V\,\mathrm{d}t\right].
    \end{equation}
    To estimate the last term put $s=\mathcal E_{\mathcal T}(Y)$. On $\{s\leq1\}$, $\|w\|_V\leq D_{N,M}$ almost everywhere before $\tau$. For every $\delta>0$, \eqref{eq:controlled-H-comparison} gives
    \begin{equation}\label{eq:controlled-weighted-error}
        \int_0^\tau a_{\varphi_\varepsilon}\|w\|_V\,\mathrm{d}t \leq J_N\delta+D_{N,M}\omega_N(C_{N,M}s/\delta^2).
    \end{equation}
    First let $s\downarrow0$, then $\delta\downarrow0$. Equations \eqref{est:uniform-stopped}--\eqref{eq:controlled-weighted-error} prove, for each fixed $M>K_N+1$ and $r>0$,
    \begin{equation}\label{eq:controlled-stopped-convergence}
        \lim_{\varepsilon\downarrow0} \sup_{\varphi_\varepsilon\in\widetilde{\mathcal U}^N} \bar{\mathbb P}\{\mathcal E_\tau(X^\varepsilon-U^\varepsilon)>r\}=0.
    \end{equation}
    The suprema include the endpoint $\tau$: the jump of $X$ there is exactly the jump of $Y$, while $w$ and $U$ are continuous.

    \emph{Removal of the stop and convergence in law.} Choose $M>2K_N+2$. On $\{\mathcal E_{\mathcal T}(X^\varepsilon)\geq M\}$, the exit is attained and $\mathcal E_\tau(X^\varepsilon)\geq M$. But
    \begin{equation*}
        \mathcal E_\tau(X^\varepsilon) \leq2K_N+2\mathcal E_\tau(X^\varepsilon-U^\varepsilon).
    \end{equation*}
    Thus \eqref{eq:controlled-stopped-convergence} implies
    \begin{equation}\label{est:chebyshev-tau}
        \lim_{\varepsilon\downarrow0} \sup_{\varphi_\varepsilon\in\widetilde{\mathcal U}^N} \bar{\mathbb P}\{\mathcal E_{\mathcal T}(X^\varepsilon)\geq M\}=0.
    \end{equation}
    Outside this event $\tau=\mathcal T$, so
    \begin{equation}\label{eq:controlled-strong-convergence}
        \mathcal E_{\mathcal T}(X^\varepsilon-U^\varepsilon) \to0\quad\text{in probability}.
    \end{equation}
    Finally, Lemma~\ref{lem:cond-1} and the continuous mapping theorem give $U^\varepsilon\Rightarrow\bm u^\varphi$ in $\Upsilon_{\mathcal T}(V)$. The quantity in \eqref{eq:controlled-strong-convergence} controls the uniform $V$ distance and the $L^2\mathscr D(A)$ distance, hence also $d_\Upsilon(X^\varepsilon,U^\varepsilon)$. The convergence-together principle proves the assertion.
\end{proof}

\vspace{0.3cm}
\noindent{\bf Acknowledgments}
 This work is partially supported by the NSFC (Grant Nos. 12271227, 12571258 and 12671286).\\

\noindent{\bf Conflict of Interest}
The authors declare there is no conflict of interest.\\

 \noindent{\bf Data availability statement}
Data sharing not applicable to this article as no data sets were generated or analyzed during the current study.

    {\emergencystretch=1em

}
\end{document}